\documentclass[12pt, reqno]{amsart}
\usepackage{amsmath, amstext, amsbsy, amssymb, amscd}
\usepackage{amsmath}
\usepackage{amsxtra}
\usepackage{amscd}
\usepackage{amsthm}
\usepackage{amsfonts}
\usepackage{amssymb}
\usepackage{eucal}
\usepackage{color}

\input prepictex
\input pictex
\input postpictex

\newtheorem{theorem}{Theorem}[section]
\newtheorem{procedure}[theorem]{Procedure}
\newtheorem{lemma}[theorem]{Lemma}
\newtheorem{proposition}[theorem]{Proposition}
\newtheorem{corollary}[theorem]{Corollary}
\theoremstyle{definition}
\newtheorem{definition}[theorem]{Definition}
\newtheorem{example}[theorem]{Example}

\theoremstyle{remark}
\newtheorem{remark}[theorem]{Remark}
\newtheorem{conjecture}[theorem]{Conjecture}
\definecolor{A}{rgb}{.75,1,.75}

\numberwithin{equation}{section}

\newcommand{\C}{ \mathbb C }

\newcommand{\Fi}{\mathcal O_{m|\infty}^{++}}

\newcommand{\FPn}{\mathcal O_{m|n}^{+}}
\newcommand{\FPPn}{\mathcal O_{m|n}^{++}}

\newcommand{\Oi}{\mathcal O_{m+\infty}^{++}}

\newcommand{\OPn}{\mathcal O_{m+n}^{+}}
\newcommand{\OPPn}{\mathcal O_{m+n}^{++}}

\newcommand{\Xmn}{X_{m|n}}
\newcommand{\Xmi}{X_{m|\infty}}

\newcommand{\Zmn}{ \Z_+^{m|n}}

\newcommand{\glmpn}{{\mathfrak g\mathfrak l}(m+n)}
\newcommand{\g}{\gamma}
\newcommand{\gl}{{\mathfrak g\mathfrak l} (m+\infty)}
\newcommand{\glmn}{{\mathfrak g\mathfrak l}(m|n)}
\newcommand{\glsuper}{{\mathfrak g\mathfrak l}(m|\infty)}
\newcommand{\hf}{\frac12}
\newcommand{\la}{\lambda}
\newcommand{\La}{\Lambda}
\newcommand{\N}{\mathbb N}

\newcommand{\vac}{|0\rangle}
\newcommand{\vacc}{|0_* \rangle}
\newcommand{\wt}{ \text{wt} }
\newcommand{\wgv}{\Lambda^{\infty} \mathbb V}
\newcommand{\wgw}{\Lambda^{\infty} \mathbb V^*}
\newcommand{\Z}{ \mathbb Z }

{\vskip-\lastskip\medskip
  \noindent
  {\em #1.}\enspace
  }%
{\qed\par\medskip
  }

\begin{document}
\title[Super duality and Kazhdan-Lusztig polynomials] {Super duality and Kazhdan-Lusztig polynomials}

\author[Shun-Jen Cheng]{Shun-Jen Cheng}
\address{Institue of Mathematics, Academia Sinica, Taipei, Taiwan
11529}\email{chengsj@math.sinica.edu.tw}

\author[Weiqiang Wang]{Weiqiang Wang}
\address{Department of Mathematics, University of Virginia,
Charlottesville, VA 22904} \email{ww9c@virginia.edu}

\author[R.B.~Zhang]{R.B.~Zhang}
\address{School of Mathematics and Statistics, University of Sydney, New South Wales
2006, Australia}\email{rzhang@maths.usyd.edu.au}

\begin{abstract}
We establish a direct connection between the representation theories
of Lie algebras and Lie superalgebras (of type $A$) via Fock space
reformulations of their Kazhdan-Lusztig theories. As a consequence,
the characters of finite-dimensional irreducible modules of the
general linear Lie superalgebra are computed by the usual parabolic
Kazhdan-Lusztig polynomials of type $A$. In addition, we establish
closed formulas for canonical and dual canonical bases for the
tensor product of any two fundamental representations of type $A$.
\end{abstract}

\maketitle
\date{}
\tableofcontents

\section{Introduction}
\subsection{Background}

Since the classification of finite-dimensional complex simple Lie
superalgebras by Kac \cite{K1}, there has been an enormous amount
of literature devoted to the representation theory of Lie
superalgebras, in particular that of the general linear
superalgebras $\glmn$ over $\C$ (cf. \cite{BL, BR, Br, Br2, CW,
CZ, JHKT, JZ, K2, Ku, PS, Se1, Se2, Sv, Zou} and the references
therein). While having many similarities with the representation
theory of Lie algebras, it was soon realized however that the
study of representations of Lie superalgebras is in general very
different and difficult: often a given method borrowed from
semisimple Lie algebras (e.g. Borel-Weil construction, Weyl
characters, etc.) can only cover limited classes of
representations. In the case when a weight $\la$ is typical, the
irreducible highest weight module $L_n(\la)$ of $\glmn$ coincides
with the so-called Kac module, whose character was obtained in
\cite{K2}.

For the study of representation theory of reductive Lie algebras,
the Kazhdan-Lusztig (KL) theory \cite{KL} has proved to be an
extremely powerful tool. The Kazhdan-Lusztig conjecture which
expresses the characters of irreducible modules of a simple Lie
algebra in terms of KL polynomials has been established
independently in \cite{BB, BK}. A decade ago, the idea of KL
theory was applied by Serganova \cite{Se1, Se2} (with a mixture of
algebraic and geometric technique) to offer a complete solution to
the longstanding problem of finding the characters of
finite-dimensional irreducible modules of $\glmn$.

Recently, partly inspired by the ideas of Lascoux, Leclerc and
Thibon \cite{LLT}, Brundan \cite{Br} offered a very different,
purely algebraic, elegant Fock space approach to the irreducible
character problem for $\glmn$. Brundan's Fock space\footnote{In
\cite{Br}, the roles of $\mathbb V$ and $\mathbb V^*$ are
switched; but this is not essential.}
$$\mathcal E^{m|n}
:= \Lambda^m \mathbb V \otimes \Lambda^n \mathbb V^*,$$
where $\mathbb V$ (respectively $\mathbb V^*$) is the natural
(respectively dual) module of the quantum group $\mathcal U =
U_q(\mathfrak s\mathfrak l(\infty))$, affords several bases
such as the monomial basis and the Kashiwara-Lusztig (dual)
canonical bases (\cite{Lu1, Kas, Lu}). Brundan computed the
transition matrices of KL polynomials between the (dual) canonical
basis and the monomial basis for $\mathcal E^{m|n}$, and further
established in terms of the KL polynomials the character formulas
of the irreducible modules and of the tilting modules in the
category $\FPn$ of finite-dimensional modules of $\glmn$. (We will
denote by $\Fi$ an analogous category, see
Section~\ref{sec:superRep}). His definition of KL polynomials was
also shown to coincide with the one given originally in \cite{Se1,
Se2} and it has a cohomological interpretation in the sense of
Vogan \cite{Vo} (also cf. \cite{CPS}, \cite{Zou}).

The category $\FPn$ is reminiscent of the parabolic category
$\OPn$ of $\glmpn$-modules which has been well studied in the
literature (cf. e.g. \cite{CC, Deo, ES}; see
Section~\ref{sec:reductiveRep} for a precise definition, where an
analogous category $\Oi$ of $\mathfrak g\mathfrak
l(m+\infty)$-modules is also defined). Kac modules of $\glmn$ can
be viewed as super analogues of generalized Verma modules. The
differences between the categories $\FPn$ and $\OPn$ are however
rather significant: a block in $\FPn$ usually contains infinitely
many simple objects, and it is controlled by the group $S_{m|n}
=S_m \times S_n$ and odd reflections (as demonstrated in
\cite{Se1}). On the other hand, a block in $\OPn$ contains
finitely many objects and it is controlled by the Weyl group
$S_{m+n}$. In the category $\OPn$ tilting modules are usually
different from the projective modules \cite{CoI}. In contrast, the
tilting modules in $\FPn$ coincide with the projective covers of
the simple modules \cite{Br}.
\subsection{The main results}

This paper provides for the first time a direct connection between
the representation theories of Lie algebras and Lie superalgebras
of type $A$ via a Fock space reformulation of the Kazhdan-Lusztig
theory. Motivated by \cite{Br} (also see \cite{LLT}), we
reformulate the KL theory for the category $\OPn$ in terms of the
canonical and dual canonical bases on the Fock space
$$\mathcal E^{m+n} := \Lambda^m \mathbb V \otimes
\Lambda^n \mathbb V$$
via the Schur-Jimbo duality \cite{Jim}. Brundan, via personal
communication, has informed us that he was aware of this. Our key
starting point here is that one has to take the limit $n\to\infty$
in order to make a precise connection with the category $\FPn$.
Let us summarize our main results.

\begin{itemize}
\item There exists a $\mathcal U$-module isomorphism of Fock
spaces of semi-infinite $q$-wedges $\Lambda^\infty \mathbb V \cong
\Lambda^\infty \mathbb V^*$. This induces a $\mathcal U$-module
isomorphism
$$ \natural: \widehat{\mathcal E}^{m+\infty}  \stackrel{\cong}{\longrightarrow} \widehat{\mathcal E}^{m|\infty},$$
which commutes with the bar involutions on $\widehat{\mathcal
E}^{m+\infty}$ and $\widehat{\mathcal E}^{m|\infty}$\footnote{The
hat here stands for a certain topological completion. Here we will
be slightly imprecise by ignoring various topological completions
for the sake of simplicity. These will be made rigorous later on.}.
It follows that the map $\natural$ sends the monomial, canonical,
and dual canonical bases for $\widehat{\mathcal E}^{m+\infty}$ to
the corresponding bases for $\widehat{\mathcal E}^{m|\infty}$
(Theorem~\ref{correspondence}).

\item We obtain closed formulas for the transition matrices among
the canonical, dual canonical, and the monomial bases in
${\mathcal E}^{m+n}$. Our formulas generalize those of
Leclerc-Miyachi \cite{LM} and our method is inspired by \cite{Br,
LM}. We introduce the truncation maps and use them to show that
the combinatorics of the (dual) canonical bases for
$\widehat{\mathcal E}^{m+\infty}$ (respectively $\widehat{\mathcal
E}^{m|\infty}$) is equivalent to those of ${\mathcal E}^{m+n}$
(respectively ${\mathcal E}^{m|n}$) for all $n$. As a consequence,
the Brundan-Serganova KL polynomials for $\glmn$ are exactly the
usual parabolic KL polynomials of type $A$. (See
Remark~\ref{rem:KLn=N}).

\item We establish an isomorphism of $\mathcal U_{q=1}$-modules
between ${\mathcal E}^{m+n}|_{q=1}$ and the rational Grothendieck
group of $\OPn$, where the Chevalley generators act via the
translation functors and the canonical, dual canonical, and
monomial bases are identified with the tilting, irreducible, and
the generalized Verma modules, respectively.

\item The tilting modules for the categories $\FPn$ are shown to
be compatible under the truncation functor (see Definition
\ref{def:trunc}), and then are `glued' together into a module in
the category $\Fi$ which is shown to be a tilting module. There
exists a natural isomorphism of the Grothendieck groups of $\Oi$
and $\Fi$ that matches the tilting, generalized Verma, and
irreducible modules with the tilting, Kac, and irreducible
modules, respectively. This isomorphism (called {\em Super
Duality}) is shown to be compatible with tensor products. The
categories $\Oi$ and $\Fi$ are conjectured to be equivalent.
\end{itemize}

\subsection{Further relations to other works}

The results of this paper in addition provide conceptual
clarification of various results and empirical observations in the
literature. Here are a few examples.

For a given positive integer $n$, let $\la =(\la_{-m}, \cdots,
\la_{-1}, \la_1, \cdots, \la_N)$ be a partition with $\la_1\le n$,
and let $\la^\natural =(\la_{-m}, \cdots, \la_{-1}, \la_1',
\cdots, \la_n')$, where $(\la_1', \cdots, \la_n')$ denotes the
conjugate partition of $(\la_1, \cdots, \la_N)$. A classical
result of Sergeev \cite{Sv} says that the irreducible
$\glmn$-modules with highest weights $\la^\natural$ associated to
such partitions are exactly those appearing in the tensor powers
of the natural module $\C^{m|n}$. Our duality implies that the
character of the irreducible $\glmn$-module $L_n(\la^\natural)$ is
given by the hook Schur polynomial, which were obtained via
different approaches in \cite{Sv} and \cite{BR}.

For a given partition $\la$ as above, since the irreducible
$\mathfrak g\mathfrak l(m+N)$-module $L_N(\la)$ admits the Weyl
character formula, our duality results imply immediately that the
irreducible $\glmn$-module $L_n(\la^\natural)$ affords a Weyl-type
character formula. This recovers and explains the results of
\cite{CZ}.

Our results also provide a conceptual explanation of the
similarity, observed in \cite{LM}, of the formulas {\em
loc.~cit.}~with the ones in \cite{Br} that the number of monomial
basis elements appearing in a canonical basis is always a power of
$2$. The combination of our stability and duality results explain
when and why these two seemingly unrelated calculations yield the
same results, and the truncation maps explain how to make sense of
the difference when some terms are missing.

We believe that the duality principle between Lie algebras and Lie
superalgebras, formulated and established here for type $A$,
provides a new approach to the representation theory of Lie
superalgebras and should be applicable to more general module
categories and other types of Lie (super)algebras. It would also
be interesting to extend this to the positive characteristic case
(c.f.~\cite{Ku}). We intend to address these issues in sequels to
this paper.

\subsection{Organization and Acknowledgments}

The paper is organized as follows. In Section \ref{sec:basic} we
set up the basics on (dual) canonical bases on Fock spaces. We
introduce the truncation maps between Fock spaces and show the bar
involution commutes with the truncations. Section
\ref{sec:superRep} addresses the representation theory of $\glmn$
and $\glsuper$. In Section~\ref{sec:bar}, we formulate and
establish closed formulas for (dual) canonical bases on Fock
spaces $\mathcal E^{m+n}$ and $\widehat{\mathcal E}^{m+\infty}$.
In Section~\ref{sec:reductiveRep} we reformulate a parabolic KL
theory for the representations of $\glmpn$ in terms of the Fock
space $\mathcal E^{m+n}$. In Section \ref{sec:isom} we establish
the isomorphism of the Fock spaces and the isomorphism of
Grothendieck groups with favorable properties.

The first and third authors thank University of Virginia for
hospitality and support. The first author is supported by an NSC
grant of the R.O.C. and is a member of the NCTS Taipei Office and
the Tai-Da Institute for Mathematical Sciences, while the second
author is supported by NSF. We thank Jon Brundan for helpful
comments and references, and Bernard Leclerc for the very
interesting and helpful reference \cite{LM}. We thank Leonard Scott
for helpful consultations on the KL theory. Finally we are greatly
indebted to an anonymous expert for numerous thoughtful suggestions,
criticism, and corrections that have resulted in this greatly
improved version. Notation: $\N =\{0,1,2,\cdots\}$.

\section{Basic constructions of Fock spaces}
\label{sec:basic}
\subsection{Basics on quantum groups and the Iwahori-Hecke algebra}

The quantum group $U_q({\mathfrak g\mathfrak l}(\infty))$ is the
$\mathbb Q(q)$-algebra generated by $E_a, F_a, K^{\pm}_a, a \in
\Z$, subject to the relations

\begin{eqnarray*}
 K_a K_a^{-1} =K_a^{-1} K_a =1, &&
 K_a K_b = K_b K_a, \\
 K_a E_b K_a^{-1} = q^{\delta_{a,b} -\delta_{a,b+1}} E_b, &&
 K_a F_b K_a^{-1} = q^{\delta_{a,b+1}-\delta_{a,b}}
 F_b, \\
 E_a F_b -F_b E_a &=& \delta_{a,b} (K_{a,a+1}
 -K_{a+1,a})/(q-q^{-1}), \\
 E_a E_b = E_b E_a,   &&
 F_a F_b = F_b F_a,  \qquad\qquad \text{if } |a-b|>1,\\
 E_a^2 E_b +E_b E_a^2 &=& (q+q^{-1}) E_a E_b E_a,  \qquad\qquad \text{if } |a-b|=1, \\
 F_a^2 F_b +F_b F_a^2 &=& (q+q^{-1}) F_a F_b F_a,  \qquad\qquad \text{if } |a-b|=1.
\end{eqnarray*}
Here $K_{a,a+1} :=K_aK_{a+1}^{-1}, a\in \Z$.
Define the
bar involution on $U_q({\mathfrak g\mathfrak l}(\infty))$
to be the anti-linear automorphism
$ ^-:  E_a\mapsto E_a, \quad F_a\mapsto F_a, \quad K_a\mapsto
K_a^{-1}.$
Here by {\em anti-linear} we mean the automorphism of $\mathbb
Q(q)$ given by $q \mapsto q^{-1}$.

Let $\mathbb V$ be the natural $U_q({\mathfrak g\mathfrak
l}(\infty))$-module with basis $\{v_a\}_{a\in\Z}$ and $\mathbb W
:=\mathbb V^*$ the dual module with basis $\{w_a\}_{a\in\Z}$ such
that $w_a (v_b) = (-q)^{-a} \delta_{a,b}$.  We have
\begin{align*}
K_av_b=q^{\delta_{ab}}v_b,\quad E_av_b=\delta_{a+1,b}v_a,\quad
F_av_b=\delta_{a,b}v_{a+1},\\
K_aw_b=q^{-\delta_{ab}}w_b,\quad E_aw_b=\delta_{a,b}w_{a+1},\quad
F_aw_b=\delta_{a+1,b}w_{a}.
\end{align*}
Following \cite{Br} we shall use the co-multiplication $\Delta$ on
$U_q({\mathfrak g\mathfrak l}(\infty))$ defined by:
\begin{eqnarray*}
 \Delta (E_a) &=& 1 \otimes E_a + E_a \otimes K_{a+1, a}, \\
 \Delta (F_a) &=& F_a \otimes 1 +  K_{a, a+1} \otimes F_a, \quad
 \Delta (K_a) = K_a \otimes K_{ a}.
\end{eqnarray*}
We let $\mathcal U = U_q({\mathfrak s\mathfrak l}(\infty))$ denote
the subalgebra with generators $E_a$, $F_a$, $K_{a,a+1}, a\in \Z$.

For $m \in \N, n\in\N \cup \infty$, we let
$I(m|n):=\{-m,-m+1,\cdots,-1\}\cup\{1,2,\cdots,n\}.$
Denote by $S_{m+n}$ the symmetric group of (finite) permutations
on $I(m|n)$, and $S_{m|n}$ its subgroup $S_m \times S_n$.
Then $S_{m+n}$ is generated by the simple
transpositions
$
s_{-m+1} =(-m, -m+1), \ldots, s_{-1} =(-2, -1),
s_0 =(-1, 1),
s_1=(1,2), \ldots, s_{n-1} =(n-1, n)
$
and $S_{m|n}$ is generated by those $s_i$ with $i \neq 0$.

The Iwahori-Hecke algebra $\mathcal H_{m+n}$ is the $\mathbb
Q(q)$-algebra with generators $H_i$, where $-m+1 \le i \le n-1$,
subject to the relations
%
$ (H_i -q^{-1})(H_i +q) = 0,
 H_i H_{i+1} H_i = H_{i+1} H_i H_{i+1},
 H_i H_j = H_j H_i, \text{for } |i-j| >1$.
Associated to $x \in S_{m+n}$ with a reduced
expression $x=s_{i_1} \cdots s_{i_r}$, we define $H_x :=H_{i_1}
\cdots H_{i_r}$. The bar involution on $\mathcal H_{m+n}$ is the
unique anti-linear automorphism defined by $\overline{H_x}
=H_{x^{-1}}^{-1}$ for all $x \in S_{m+n}$.
We denote by $\mathcal H_{m|n}$ the subalgebra generated by those
$H_i$ with $i \neq 0$, that is $\mathcal H_{m|n}$ is the Hecke
algebra corresponding to $S_{m|n}$. Denote by $w_0$ the longest
element in $S_{m|n}$, and let
$\widehat{H}_0:= \sum_{x \in S_{m|n}} (-q)^{\ell(x) - \ell (w_0)}
H_x.$

\subsection{The spaces $\wgv$, $\mathcal E^{m+n}$ and $\mathcal
E^{m+\infty}$}\label{tauij}

Let $P$ be the free abelian group with basis $\{\epsilon_a\vert
a\in\Z\}$ equipped with a bilinear form $(\cdot | \cdot)$ for
which the $\epsilon_a$'s are orthonormal. We define a partial
order on $P$ by declaring $\nu\ge\mu$ for $\nu,\mu\in P$ if
$\nu-\mu$ is a non-negative integral linear combination of
$\epsilon_a-\epsilon_{a+1}$, $a \in \Z$.
For $n \in \N \cup \infty$, we let $\Z^{m+n}$ or $\Z^{m|n}$ be the
set of integer-valued functions on $I(m|n)$.   Define
\begin{equation}\label{nonsuperweight}
\wt^{\epsilon} (f):=\sum_{i\in I(m|n)}\epsilon_{f(i)}, \qquad
\text{for } f \in \Z^{m+n}.
\end{equation}

For a finite $n$ consider $\mathbb T^{m+n} :=\mathbb V^{\otimes
(m+n)},$ where we adopt the convention that the $m+n$ tensor
factors are indexed by $I(m|n)$. Similar conventions will apply to
similar situations below. For $f \in \Z^{m+n}$, let
$$\mathcal V_f :=v_{f(-m)} \otimes \cdots  \otimes v_{f(-1)} \otimes
v_{f(1)} \otimes \cdots \otimes v_{f(n)}.$$

Let $\mathcal H_{m+n}$ act on  $\mathbb T^{m+n}$ on the right by:
\begin{eqnarray} \label{eq:heckeaction}
 \mathcal V_f H_i = \left\{
 \begin{array}{ll}
 \mathcal V_{f\cdot s_i}, & \text{if } f < f \cdot s_i,  \\
 q^{-1} \mathcal V_{f}, & \text{if } f = f \cdot s_i, \\
 \mathcal V_{f \cdot s_i} - (q-q^{-1}) \mathcal V_f, & \text{if } f > f \cdot
 s_i.
 \end{array}
 \right.
\end{eqnarray}
The {\em Bruhat ordering} $\leq$ on $\Z^{m+n}$ comes from the
Bruhat ordering on $S_{m+n}$, which is the transitive closure of
the relation $f< f \cdot \tau_{ij}$ if $f(i) < f(j)$, for $i,j \in
I(m|n)$ with $i<j$. Here and further $\tau_{ij}$ denotes the
transposition interchanging $i,j$. The algebra $\mathcal U$ acts
on $\mathbb T^{m+n}$ via the co-multiplication $\Delta$.

\begin{lemma} \cite{Jim} \label{lem:commute}
The actions of $\mathcal U$ and $\mathcal H_{m+n}$ on $\mathbb
T^{m+n}$ commute with each other.
\end{lemma}

Different commuting actions of $\mathcal U$ and the Hecke
algebra on a tensor power of $\mathbb V$ were used in \cite{KMS}
to construct the space $\Lambda^n \mathbb V$ of finite $q$-wedges
and then the space of infinite-wedges by taking the limit $n
\rightarrow \infty$ appropriately. These
spaces carry the action of quantum affine algebras as well as the
action of $\mathcal U$ (as a limiting case). The constructions in
{\em loc.~cit.}~carry over using the above actions of
$\mathcal U$ and the Hecke algebra as formulated below, and we refer to {\em loc.~cit.}~for
details.

One can think of $\Lambda^n \mathbb V$ either as the subspace
${\rm im}\widehat{H}_0$ or the quotient $\mathbb T^n/\ker
\widehat{H}_0$ of $\mathbb T^n$. Here $\ker \widehat{H}_0$ equals
the sum of the kernels of the operators $H_i -q^{-1}$, $1 \le i
\le n-1$ (note that the Hecke algebra generator $T_i$ used in {\em
loc.~cit.}~corresponds to our $-q H_i$). The {\em $q$-wedge}
$v_{a_1} \wedge \cdots \wedge v_{a_n}$ is an element of $\Lambda^n
\mathbb V$, which is the image  of $v_{a_1} \otimes \cdots \otimes
v_{a_n}$ under the canonical map when $\Lambda^n \mathbb V$ is
regarded as the quotient space $\mathbb T^n/\ker \widehat{H}_0$.
We have

\begin{eqnarray}
\cdots \wedge v_{a_i} \wedge v_{a_{i+1}} \wedge \cdots
 &=& -q^{-1} (\cdots \wedge v_{a_{i+1}} \wedge v_{a_i} \wedge
\cdots), \quad \text{if } a_i < a_{i+1}; \nonumber\\
\cdots \wedge v_{a_i} \wedge v_{a_{i+1}} \wedge \cdots
 &=& 0, \quad \text{if } a_i = a_{i+1}.  \label{eq:straighten}
\end{eqnarray}
It follows that the elements $v_{a_1} \wedge \cdots \wedge
v_{a_n}$, where $a_1 > \cdots > a_n$ form a basis for $\Lambda^n
\mathbb V$. By Lemma~\ref{lem:commute}, $\mathcal U$ acts
naturally on $\Lambda^n \mathbb V$.

By taking $n \rightarrow \infty$, one
defines $\wgv$ with a $\mathcal U$-action, with a
basis given by the infinite $q$-wedges
$v_{m_1} \wedge v_{m_2} \wedge v_{m_3} \wedge \cdots,$
where  $m_1 > m_2  >m_3> \cdots,$ and $m_i = 1-i$ for $i \gg0$
(our $\wgv$ is $F_{(0)}$ in  \cite{KMS}).
Alternatively, $\wgv$ has a basis
$$|\la\rangle := v_{\la_1} \wedge v_{\la_2 -1} \wedge v_{\la_3 -2}\wedge \cdots $$
where $\la =(\la_1, \la_2, \cdots)$ runs over the set of all
partitions. Set (for a finite $n$)
\begin{eqnarray*}
\Z_+^{m+n} &:=& \{f \in \Z^{m+n} \mid f(-m) > \cdots > f(-1), f(1)
> \cdots > f(n) \},   \\
\Z^{m+n}_{++} &:=& \{f \in \Z_+^{m+n} \mid  f(n) \ge 1-n \}.\\
\Z_+^{m+\infty} &:=& \{f \in \Z^{m+\infty} \mid f(-m) > \cdots >
f(-1), \\
&& \qquad\qquad \quad f(1)  > f(2) > \cdots; f(i) =1-i \text{ for
} i\gg0 \}.
\end{eqnarray*}
We will call an element in $\Z^{m+n}_+$ or
$\Z^{m+\infty}_+$ {\em dominant}.

For $n \in \N \cup \infty$, we let
$$\mathcal E^{m+n} := \Lambda^m \mathbb V \otimes \Lambda^n \mathbb V,$$
where the factors are indexed by $I(m|n)$. It
has the {\em monomial basis}
$$\mathcal K_f :
= \left\{
 \begin{array}{ll}
 v_{f(-m)} \wedge \cdots  \wedge v_{f(-1)} \otimes
v_{f(1)} \wedge \cdots \wedge v_{f(n)}, & \text{for finite } n, \\
 v_{f(-m)} \wedge \cdots  \wedge v_{f(-1)} \otimes
v_{f(1)} \wedge v_{f(2)} \wedge \cdots , & \text{for } n=\infty,
 \end{array}
 \right.$$
where $f$ runs over the set $\Z_+^{m+n}$. Recalling $\mathbb
T^{m+n} = \ker \widehat{H}_0 \oplus \text{Im} \widehat{H}_0$,
cf.~\cite{KMS}, we may regard equivalently $\mathcal E^{m+n}$ as
the subspace $\text{Im} \widehat{H}_0$ of $\mathbb T^{m+n}$ for
$n$ finite.

%
%
%

Let $g \in \Z_+^{m+0}$ or $g \in \Z_+^{0+n}$. For $a\in\Z$ define
$\tilde{e}_{a,a+1}g$ (respectively $\tilde{e}_{a,a-1}g$) to be the
strictly decreasing sequence of integers obtained from $g$ by
replacing the value that is $a+1$ (respectively $a-1$) by $a$, if
$a+1$ (respectively $a-1$) appears in $g$ and $a$ does not appear
in $g$. Otherwise set $\tilde{e}_{a,a\pm1}g =\emptyset$. For
$g\in\Z^{m+n}$ we let $g^{<0}$ and $g^{>0}$ denote the
restrictions of $g$ to $\{-m,\cdots,-1\}$ and $\{1,\cdots,n\}$,
respectively, and write $g=(g^{<0}|g^{>0})$. Now for $m,n>0$ let
$g\in\Z^{m+n}_+$. We set $\mathcal K_g=0$ if $g$ is of the form
$(g^{<0}|\emptyset)$ or $(\emptyset|g^{>0})$. The formulas for
$\Delta$ on $\Lambda^m\mathbb V\otimes\Lambda^n\mathbb W$ and the
straightening relations (\ref{eq:straightensuper}) give us the
following formula.

\begin{lemma}\label{lem:uactsonk}
Let $n \in \N \cup \infty$. For $f=(f^{<0}|f^{>0})\in\Z^{m+n}_+$,
$\mathcal U =U_q(\mathfrak s\mathfrak l(\infty))$ acts on
$\mathcal E^{m+n}$ as follows:
\begin{align*}
&E_a (\mathcal K_{(f^{<0}|f^{>0})}) =\mathcal
K_{(f^{<0}|\tilde{e}_{a,a+1}f^{>0})} +q^{-(\wt^\epsilon
(f^{>0})|\epsilon_a)+(\wt^\epsilon
(f^{>0})|\epsilon_{a+1})}\mathcal K_{(\tilde{e}_{a,a+1}f^{<0}|f^{>0})},\\
&F_a (\mathcal K_{(f^{<0}|f^{>0})}) =\mathcal
K_{(\tilde{e}_{a+1,a}f^{<0}|f^{>0})} +q^{(\wt^\epsilon (f^{<0})
|\epsilon_a) -(\wt^\epsilon (f^{<0}) |\epsilon_{a+1})}\mathcal
K_{(f^{<0}|\tilde{e}_{a+1,a}f^{>0})}.
\end{align*}
\end{lemma}
\subsection{Super Bruhat ordering on $\Z^{m|n}$}\label{def:di}
 Let $n\in\N\cup\infty$. For $i \in I(m|n)$ we define $d_i \in \Z^{m|n}$ by $j \mapsto -
\text{sgn} (i) \delta_{ij}$. For $f, g\in\Z^{m|n}$, we write $f
\downarrow g$ if one of the following holds:
\begin{enumerate}
\item
$g=f-d_i +d_j$ for some $i<0 <j$ such that $f(i) =f(j)$;
\item
$g =f \cdot \tau_{ij}$ for some $i<j<0$ such that $f(i) >f(j)$;
\item
$g =f \cdot \tau_{ij}$ for some $0<i<j$ such that $f(i) <f(j)$.
\end{enumerate}
The {\em super Bruhat ordering} on $\Z^{m|n}$ is defined as
follows: for $f,g\in\Z^{m|n}$, we say that $f\succ g$, if there
exists a sequence $f=h_1,   \ldots, h_r=g \in \Z^{m|n}$ such that
$h_1  \downarrow h_2, \cdots, h_{r-1} \downarrow h_r$. It can also
be described (cf. \cite{Br}) in terms of the $\epsilon$-weights,
which are defined by:
\begin{equation}\label{superweight}
\wt^{\epsilon} (f):=\sum_{i\in
I(m|n)}-\text{sgn}(i)\epsilon_{f(i)}, \qquad \text{for } f \in
\Z^{m|n}.
\end{equation}
The super Bruhat ordering is defined to be compatible with the one
on the set of integral weights for $\mathfrak g\mathfrak l(m|n)$.
{\em Here and further the superscript ${Br}$ stands for Brundan's
version \cite{Br}.} Our weight $f$ corresponds to $f^{Br}=-f$. Our
$d_i$ here differs from the one in \cite{Br} by a sign. Set (for a
finite $n$)
\begin{eqnarray*}
\Z_+^{m|n} &:=& \{f \in \Z^{m|n} \mid f(-m) > \cdots > f(-1), f(1)
< \cdots < f(n) \},   \\
\Z_{++}^{m|n} &:=& \{f \in \Z_+^{m|n} \mid f(n)\le n \},   \\
\Z_+^{m|\infty} &:=& \{f \in \Z^{m|\infty} \mid f(-m) > \cdots >
f(-1), \\
&& \qquad\qquad \quad f(1)  < f(2) < \cdots; f(i) =i \text{ for }
i\gg0 \}.
\end{eqnarray*}

An element in $\Z^{m|n}_+$ is called a {\em dominant weight}.
(Note that $\Z_+^{m|n} = - \Z_+^{m|n,{Br}}$.) For
$n\in\N\cup\infty$ and $f\in\Z^{m|n}$, conjugate under the action
of $S_{m|n}$ to an element in $\Z^{m|n}_+$, define the {\em degree
of atypicality} of $f$ to be the
$$\#f:=|\{f(-m),\cdots,f(-1)\}\cap\{f(1),f(2),\cdots\}|.$$
If $f, g \in  \Z^{m|n}_+$ are comparable under the super Bruhat
ordering, then $\#f =\#g.$ If $\#f=0$, we say that $f$ is {\em
typical}. For later use we introduce the following.
\begin{definition}
Let $f\in\Z^{m|k}_{+}$ (respectively $f\in\Z^{m+k}_{+}$), with
$k\in\N\cup\infty$, $n\le k$. Define $f^{(n)} \in\Z^{m|n}_{+}$
(respectively $f^{(n)} \in\Z^{m+n}_{+}$) to be the restriction of
$f$ to $I(m|n)$. We say that $n$ is {\em sufficiently large for
$f$}, if $\#f=\#f^{(n)}$ and $f(n')=n'$ (respectively $f(n')=
1-n'$) for $n' \ge n+1$. When $f$ is clear from the context we
shall simply write $n\gg0$.
\end{definition}
\subsection{The spaces $\wgw$, $\mathcal E^{m|n}$ and $\mathcal E^{m|\infty}$}

Recall that $\mathbb W =\mathbb V^*$ is the dual module of
$\mathcal U$ with basis $\{w_a\}_{a \in \Z}$. The space
$\mathbb T^{m|n} : = \mathbb V^{\otimes m} \otimes \mathbb W^{\otimes n}$
has the {\em monomial basis}
$M_f :=v_{f(-m)} \otimes \cdots  \otimes v_{f(-1)} \otimes
w_{f(1)} \otimes \cdots \otimes w_{f(n)}$, $f \in \Z^{m|n}$.
Let $\mathcal H_{m|n}$ act on $\mathbb T^{m|n}$ on the right (cf.
\cite{Br}) by:
\[
 M_f H_i = \left\{
 \begin{array}{ll}
 M_{f\cdot s_i}, & \text{if } f \prec f \cdot s_i \\
 q^{-1} M_{f}, & \text{if } f = f \cdot s_i \\
 M_{f \cdot s_i} - (q-q^{-1}) M_f, & \text{if } f \succ f \cdot s_i
 \end{array}
 \right.
\]
The algebra $\mathcal U$  acts on $\mathbb T^{m|n}$ via the
co-multiplication $\Delta$.

\begin{lemma} \label{lem:commute2}
The actions of $\mathcal U$ and $\mathcal H_{m|n}$ on $\mathbb
T^{m|n}$ commute with each other.
\end{lemma}

For a finite $n$, we can again define $\mathcal E^{m|n}$ to be
either the image of the operator $\widehat{H}_0$ in $\mathbb
T^{m|n}$, or the quotient of $\mathbb T^{m|n}$ by the kernel of
$\widehat{H}_0$. The following straightening relations hold in
$\mathcal E^{m|n}$ (note the conditions $a_{i-1} < a_i$ and $b_i >
b_{i+1}$ reflect the super Bruhat order)
\begin{eqnarray}
\cdots \wedge v_{a_{i-1}} \wedge v_{a_{i}} \wedge \cdots
 &=& -q^{-1} (\cdots  v_{a_{i}} \wedge v_{a_{i-1}} \wedge
\cdots), \quad \text{if } a_{i-1} < a_i, i<0; \nonumber \\
%
 %
 \cdots \wedge w_{b_i} \wedge w_{b_{i+1}} \wedge \cdots
 &=& -q^{-1} (\cdots  w_{b_{i+1}} \wedge w_{b_i} \wedge
\cdots), \; \text{if } b_i > b_{i+1}, i>0. \label{eq:straightensuper}
%
\end{eqnarray}

We can repeat the procedure in \cite{KMS}
to construct the space $\Lambda^\infty \mathbb W$ of
semi-infinite $q$-wedges
$w_{n_1} \wedge w_{n_2}  \wedge \cdots,$
where $n_i =i$ for $i\gg0$, which carries a $\mathcal
U$-module structure. Writing the conjugate partition of $\la$ as
$\la' =(\la'_1, \la'_2, \cdots)$, we set
$$ |\la'_*\rangle := w_{1 -\la'_1} \wedge w_{2 -\la'_2} \wedge w_{3
-\la'_3} \wedge \cdots.$$ The set $\{|\la'_*\rangle\}$ provides
another basis for $\Lambda^\infty \mathbb W$.
Note that for a finite $n$, we have $\mathcal
E^{m|n} = \Lambda^m \mathbb V \otimes \Lambda^n \mathbb W$ as
$\mathcal U$-modules. Denote
$$\mathcal E^{m|\infty} := \Lambda^m \mathbb V \otimes
\Lambda^\infty \mathbb W,$$
which is a $\mathcal U$-module via $\Delta$. $\mathcal E^{m|n}$
has the {\em monomial basis}
$$K_f := \left\{
 \begin{array}{ll}
 v_{f(-m)} \wedge \cdots  \wedge v_{f(-1)} \otimes
w_{f(1)} \wedge \cdots \wedge w_{f(n)}, & \text{for finite } n, \\
 v_{f(-m)} \wedge \cdots  \wedge v_{f(-1)} \otimes
w_{f(1)} \wedge w_{f(2)} \wedge \cdots , & \text{for } n=\infty,
 \end{array}
 \right.$$
 where $f$ runs over $\Z_+^{m|n}$.
\subsection{Bases for $\widehat{\mathcal E}^{m|n}$ and  $\widehat{\mathcal
E}^{m|\infty}$}\label{aux:base}

Let $n \in \N \cup \infty$. For $d\in\N$ let ${\mathcal
E}^{m|n}_{\ge -d}$ be the $\mathbb Q(q)$-subspace of $\mathcal
E^{m|n}$ spanned by $K_f$ with $f(i)\ge -d$, for all $i\in
\{1,\cdots,n\}$. We shall denote a certain topological completion
of ${\mathcal E}^{m|n}$ by $\widehat{\mathcal E}^{m|n}$ whose
elements may be viewed as infinite $\mathbb Q(q)$-linear
combinations of elements in $\mathcal E^{m|n}$, which under the
projection onto $\mathcal E^{m|n}_{\ge -d}$ are finite sums for
all $d\in\N$ (cf.~\cite[\S2-d]{Br}). The following proposition for
a finite $n$ is Theorem~3.5 \cite{Br}, which is similar to results
of Lusztig \cite{Lu}. The proof of $n=\infty$ case is similar and
in a sense even simpler, since $\wgw$ is a highest weight
$\mathcal U$-module.

\begin{proposition}\label{thm:involutionsuper}
Let $n \in \N \cup \infty$. There exists a unique continuous,
anti-linear bar map ${}^-: \widehat{\mathcal E}^{m|n} \rightarrow
\widehat{\mathcal E}^{m|n}$ such that
\begin{enumerate}
\item $\overline{K_f} = K_f$ for all typical $f \in
\Z_+^{m|n}$.

\item $\overline{X u} = \overline{X} \overline{u}$ for all $X \in
\mathcal U$ and $u \in \widehat{\mathcal E}^{m|n}$.

\item The bar map is an involution. \item $\overline{K_f} = K_f +
(*)$ where $(*)$ is (possibly infinite) $\Z [q,q^{-1}]$-linear
combination of $K_g$'s for $g \in \Z_+^{m|n}$ with $g  \prec f$.
\end{enumerate}
\end{proposition}

The next theorem now follows by standard arguments (cf. \cite{KL, Lu,
Du, Br}).

\begin{theorem}
Let $n \in \N \cup \infty$. There exist unique {\em canonical basis}
$\{U_f\}$ and {\em dual canonical basis} $ \{L_f \}$, where ${f \in \Z_+^{m|n}}$,  for
$\widehat{\mathcal E}^{m|n}$ such that
\begin{enumerate}
\item $\overline{U}_f =U_f$ and $\overline{L}_f =L_f$.
\item $U_f \in K_f + \widehat{\sum}_{g \in \Zmn} q \Z [q] K_g$
 and $L_f \in K_f + \widehat{\sum}_{g \in \Zmn} q^{-1} \Z [q^{-1}]
 K_g$.
\item $U_f = K_f + (*)$ and  $L_f = K_f + (**)$ where $(*)$ and
$(**)$ are (possibly infinite) $\Z [q,q^{-1}]$-linear combinations
of $K_g$'s for $g \in \Z_+^{m|n}$ such that $g  \prec f$.
\end{enumerate}
\end{theorem}

Let $n \in \N \cup \infty$. Following \cite{Br}, we define the
Kazhdan-Lusztig polynomials $u_{g,f}(q) \in \Z[q], \ell_{g,f}(q)
\in \Z[q^{-1}]$ associated to $f,g \in \Z_+^{m|n}$ by
\begin{eqnarray}  \label{eq:superul}
U_f =\sum_{g \in \Z_+^{m|n}} u_{g,f}(q) K_g, \qquad
L_f =\sum_{g \in \Z_+^{m|n}} \ell_{g,f}(q) K_g.
\end{eqnarray}
Note that $u_{g,f}(q) =\ell_{g,f}(q)=0$ unless $g \preccurlyeq f$,
$u_{f,f}(q) =\ell_{f,f}(q) =1$, and $u_{g,f}(q)\in q\Z[q]$,
$\ell_{g,f}(q) \in q^{-1} \Z[q^{-1}]$ for $g \ne f$.

\subsection{The truncation map}
\label{trunct-bar}

Let $n$ be finite. Denote by ${\mathcal E}^{m|n}_+$ the subspace
of ${\mathcal E}^{m|n}$ spanned by $K_f$ for $f \in
\Z_{++}^{m|n}$. For $n>0$, define the subalgebra $U_q({\mathfrak
s\mathfrak l}_{\le n})$ to be the subalgebra of $U_q({\mathfrak
s\mathfrak l}(\infty))$ generated by $E_{a-1}$, $F_{a-1}$ and
$K_{a-1,a}$, with $a\le n$.

\begin{lemma} \label{embed+}
${{\mathcal E}}^{m|n}_+$ is a bar-invariant subspace of
${{\mathcal E}}^{m|n}$, and it is a $U_q({\mathfrak s\mathfrak
l}_{\le n})$-module.
\end{lemma}

\begin{proof}
The first half follows from
Proposition~\ref{thm:involutionsuper}~(4), and the second half is
an easy consequence of Lemma~\ref{lem:uactsonk}.
\end{proof}

We define the {\em truncation map} to be the linear map
$ {\mathfrak{Tr}}_{n+1,n}:  \widehat{\mathcal E}^{m|n+1}_+
\longrightarrow \widehat{\mathcal E}^{m|n}_+ $
which sends $K_f$ to $K_{f^{(n)}}$ if $f(n+1) =n+1$, and to $0$
otherwise. For $n'>n$, we define the truncation map
${\mathfrak{Tr}}_{n',n}: \widehat{\mathcal E}^{m|n'}_+
\longrightarrow \widehat{\mathcal E}^{m|n}_+$ to be
${\mathfrak{Tr}}_{n',n}:= {\mathfrak{Tr}}_{n+1,n} \circ \cdots \circ
{\mathfrak{Tr}}_{n',n'-1}.$ Similarly, the {\em truncation map}
${\mathfrak{Tr}}_{n}: \widehat{\mathcal E}^{m|\infty}
\longrightarrow \widehat{\mathcal E}^{m|n}_+ $ is defined by sending
$K_f$ to $K_{f^{(n)}}$ if $f =(f^{(n)},n+1, n+2, \cdots)$, and to
$0$ otherwise.

\begin{proposition} \label{commutativity2}
The truncation map ${\mathfrak{Tr}}_{n+1, n}: \widehat{{\mathcal
E}}^{m|n+1}_+\rightarrow\widehat{{\mathcal E}}^{m|n}_+$ commutes
with the bar-involution.
\end{proposition}

\begin{proof}

Throughout the proof $\mathcal E^{m|n}$ will be regarded as a
subspace of ${\mathbb T}^{m|n}$ for all $m$ and $n$. Let
$\widehat{\mathbb T}^{m|n, 1}$ be a topological completion of
${{\mathcal E}}^{m|n}\otimes \mathbb W$ defined in an analogous
way as the completion $\widehat{{\mathcal E}}^{m|n}$ in Subsection
\ref{aux:base}. Thus $\widehat{\mathbb T}^{m|n, 1}$ contains
$\widehat{{\mathcal E}}^{m|n+1}$ as a subspace.  Now
$K_{f^{(n)}}\otimes w_{f(n+1)}$ with $f^{(n)}\in\Z^{m|n}_+$,
$f(n+1)\in\Z$, form a basis of $\widehat{\mathbb T}^{m|n, 1}$. One
can follow Lusztig (Section 24.1.1 in \cite{Lu}) to define a
bar-involution on $\widehat{\mathbb T}^{m|n, 1}$ using a quasi
$R$-matrix $\Theta^{(m|n+1)}$ such that \begin{eqnarray*}
\overline{K_{f^{(n)}}\otimes w_{f(n+1)}}:=
\Theta^{(m|n+1)}\left(\overline{K_{f^{(n)}}}\otimes
w_{f(n+1)}\right).
\end{eqnarray*}
Let $X_{a b}$ be endomorphisms of ${\mathbb W}$ such that $X_{a b}
w_c=\delta_{a c} w_b$, $a,b,c\in\Z$. By exploiting the structure
of the quasi $R$-matrix of ${\mathcal U}$ (Theorem 8.1 in
\cite{KT}, and Theorem 3 in \cite{KR}), we can show that there
exist elements $\hat{E}_{a b}$ in ${\mathcal U}$ constructed from
$E_{c, c+1}$ and $K_{c}^{\pm 1}$ only and with the property $ K_c
\hat{E}_{a b} K_c^{-1} = q^{\delta_{a c} - \delta_{b c}}
\hat{E}_{a b}$ such that
\begin{eqnarray}   \label{quasiR}
\Theta^{(m|n+1)}&=& 1\otimes 1 + (q-q^{-1}) \sum_{a<b} \hat{E}_{a b}\otimes X_{b a}. \end{eqnarray}
This is a variant of Jimbo's result \cite{Jim} on the $R$-matrix in the tensor product of an arbitrary
representation with the natural representation of ${U}_q({\mathfrak{gl}}(k))$ for finite $k$. We have
\begin{eqnarray} \overline{K_{f^{(n)}}\otimes w_{f(n+1)}} &=&\overline{K_{f^{(n)}}}\otimes w_{f(n+1)} \nonumber\\ && + (q-q^{-1}) \sum_{a< {f(n+1)}} \hat{E}_{a\, {f(n+1)}}\overline{K_{f^{(n)}}}\otimes w_a .\label{R-matrix} \end{eqnarray}

Denote by $\widehat{\mathbb T}^{m|n, 1}_+$ the subspace of
$\widehat{\mathbb T}^{m|n, 1}$ spanned by all $K_{f^{(n)}}\otimes
w_{f(n+1)}$ with $f^{(n)}\in\Z^{m|n}_+$ and $f(n+1)\le n+1$. By
\eqref{R-matrix}, $\widehat{\mathbb T}^{m|n, 1}_+$ is invariant
under the bar-involution. Introduce the following linear map
\begin{eqnarray*} \check{\mathfrak{Tr}}_{n+1, n}: \widehat{\mathbb
T}^{m|n, 1}_+\rightarrow \widehat{{\mathcal E}}^{m|n}, &&
K_{f^{(n)}}\otimes w_{f(n+1)} \mapsto \delta_{f(n+1), n+1}
K_{f^{(n)}}. \end{eqnarray*} Applying it to \eqref{R-matrix} we
see that $\check{\mathfrak{Tr}}_{n+1, n}$ commutes with the
bar-involution.

Now $\widehat{{\mathcal E}}^{m|n+1}_+=\widehat{\mathcal
E}^{m|n+1}\cap\widehat{\mathbb T}^{m|n, 1}_+$ and hence
$\widehat{{\mathcal E}}^{m|n+1}_+$ is bar-invariant. The
proposition follows since the restriction of
$\check{\mathfrak{Tr}}_{n+1, n}$ to $\widehat{{\mathcal
E}}^{m|n+1}_+$ coincides with ${\mathfrak{Tr}}_{n+1, n}$.
\end{proof}

\begin{remark}  \label{rem:trunbar}
It follows from Proposition~\ref{commutativity2} that
${\mathfrak{Tr}}_{n', n}: \widehat{{\mathcal
E}}^{m|n'}_+\rightarrow\widehat{{\mathcal E}}^{m|n}_+$ commutes
with the bar-involution for $n \le n' \le \infty$.
\end{remark}

\begin{corollary}\label{aux43}
\begin{enumerate}
\item $\{U_f \}_{f \in \Z_{++}^{m|n}}$ (respectively $\{L_f \}_{f
\in \Z_{++}^{m|n}}$) is a basis for $\widehat{\mathcal E}^{m|n}_+$.

\item ${\mathfrak{Tr}}_{n+1,n}$ sends $U_f \in \widehat{\mathcal
E}^{m|n+1}_+$ to $U_{f^{(n)}}$ if $f(n+1) =n+1$, and to $0$
otherwise.

\item ${\mathfrak{Tr}}_{n+1,n}$ sends $L_f \in \widehat{\mathcal
E}^{m|n+1}_+$ to $L_{f^{(n)}}$ if $f(n+1) =n+1$, and to $0$
otherwise.

\item For $f, g \in \Z_{++}^{m|n+1}$ such that $f(n+1) =g(n+1)=
n+1$, we have
$$u_{g,f} (q) =u_{g^{(n)}, f^{(n)}} (q), \qquad \ell_{g,f} (q) =
\ell_{g^{(n)}, f^{(n)}} (q).$$
\end{enumerate}
\end{corollary}

\begin{remark} For $n'>n$ the map ${\mathfrak{Tr}}_{n', n}$ is a
$U_q(\mathfrak{sl}_{\le n})$-module homomorphism.  Now for
$f\in\Z^{m|\infty}_+$ Procedure~3.20 in \cite{Br} is the same for
all $f^{(n)}$ with $n{\gg}0$, in the sense that it involves the
same Chevalley generators and the same sequence of weights. The
Chevalley generators lie in $U_q({\mathfrak sl}_{\le {n_0}})$, for
some fixed $n_0{\gg}0$. This together with Corollary \ref{aux43}
imply that Procedure~3.20 in \cite{Br} remains valid for
$n=\infty$ as well.
\end{remark}
\subsection{The transition matrices}
%
%
%

Let $n\in\N$ and $f \in \Z^{m|n}$ be $S_{m|n}$-conjugate to an
element in $\Z^{m|n}_+$. For $-m \le i <0 <j \le n$ with $f(i)
=f(j)$, define \cite{Br}
\begin{eqnarray*}
\texttt L_{i,j} (f) := f -a(d_i -d_j),
\end{eqnarray*}
where $a$ is the smallest positive integer such that $f -a(d_i
-d_j)$ and all $\texttt L_{k,l}(f) - a(d_i -d_j)$ for $-m\le
i<k<0<l<j\le n$ with $f(k) =f(l)$ are $S_{m|n}$-conjugate to
elements of $\Z_{+}^{m|n}$. Furthermore define \cite{JZ}
$$\texttt R_{i,j} (f) := f +b (d_i -d_j),$$
where $b$ is the smallest positive integer such that $f +b (d_i
-d_j)$ and all ${\texttt R}_{k,l}(f) + b(d_i -d_j)$ for $-m \le k
<i<0<j<l \le n$ with $f(k)=f(l)$ are $S_{m|n}$-conjugate to
elements of $\Z_{+}^{m|n}$.

Now let $f\in\Z^{m|n}_+$ and suppose $\#f=k$. Let $-m\le
i_i<i_2<\cdots<i_k\le -1$ and $1\le j_k<j_{k-1}<\cdots<j_1\le n$
be such that $f(i_l)=f(j_l)$, for $l=1,\cdots,k$.  For a $k$-tuple
$\theta=(\theta_1,\cdots,\theta_{k})\in\N^{k}$ we define \cite{Br}
\begin{align*}
{\texttt L}_\theta(f) =\Big{(}{\texttt L}^{\theta_k}_{i_k, j_k}
\circ\cdots
 \circ {\texttt L}^{\theta_1}_{i_1 ,j_1}(f)\Big{)}^+, \quad
{\texttt L}'_\theta(f) =\Big{(}{\texttt
L}^{\theta_1}_{i_1,j_1}\circ\cdots
 \circ {\texttt L}^{\theta_k}_{i_k ,j_k}(f)\Big{)}^+,\\
{\texttt R}_\theta(f) =\Big{(}{\texttt R}^{\theta_1}_{i_1
,j_1}\circ\cdots
 \circ {\texttt R}^{\theta_k}_{i_k ,j_k}(f)\Big{)}^+, \quad
{\texttt R}'_\theta(f) =\Big{(}{\texttt R}^{\theta_k}_{i_k,
j_k}\circ\cdots \circ {\texttt R}^{\theta_1}_{i_1,
j_1}(f)\Big{)}^+,
\end{align*}
where the superscript $+$ here and further stands for the unique
$S_{m|n}$-conjugate in $\Z^{m|n}_+$. We shall denote $\texttt
L_{(1,1,\cdots,1)}(f)$ and $\texttt R_{(1,1,\cdots,1)}(f)$ by
$\texttt L (f)$ and $\texttt R (f)$, or sometimes by $f^{\texttt
L}$ and $f^{\texttt R}$, respectively.

The corresponding operators ${\texttt L}_\theta$ and ${\texttt
L}'_\theta$ on $\Z^{m|\infty}_+$ are defined analogously, but it
takes extra care to make sense of the $\texttt R$ operators. Given
$f,g\in\Z^{m|\infty}_{+}$ and $\theta \in \mathbb N^{\# f}$, we
say
${\texttt R}_\theta(g)=f$
if there exists $n\gg0$ (for $f$ and $g$) so that ${\texttt
R}_\theta(g^{(n)})=f^{(n)}$. (Note the subtle point that ${\texttt
R}_\theta(g)$ is not defined for every $g\in\Z^{m|\infty}_{+}$.)

The next two lemmas follow from the definitions.

\begin{lemma}\label{lem:stabL}
Let $f=(f^{(n)},n+1)\in\Z^{m|n+1}_{++}$ and $\# f=\# f^{(n)}$, and
let $\theta\in\N^{\#f}$ be fixed. If ${\texttt L}_\theta(f)=g$,
then $g\in\Z^{m|n+1}_{++}$ and $g=({\texttt
L}_\theta(f^{(n)}),n+1)$. Similar statements hold for the operator
${\texttt L}'$ as well.
\end{lemma}

\begin{lemma} \label{lem:stableR}
Let $n\in\N$ and $f=(f^{(n)},n+1)\in\Z^{m|n+1}_{++}$ with
$\#f=\#f^{(n)}$.
\begin{itemize}
\item[(1)] Let $g\in\Z^{m|n+1}_+$ and $\theta\in\N^{\#g}$ such that ${\texttt
R}_\theta(g)=f$. Then $g=(g^{(n)},n+1)\in\Z^{m|n+1}_{++}$ with
$\#g=\#g^{(n)}$ and ${\texttt R}_\theta(g^{(n)})=f^{(n)}$.

\item[(2)] Let $\tilde{g} \in\Z^{m|n}_+$ and $\theta\in\N^{\#\tilde{g}}$ such that ${\texttt
R}_\theta(\tilde{g})=f^{(n)}$. Then
$\#(\tilde{g},n+1)=\#\tilde{g}$ and ${\texttt
R}_\theta((\tilde{g},n+1))=f$.
\end{itemize}
Similar statements hold for the operator ${\texttt R}'$ as well.
\end{lemma}

Lemmas \ref{lem:stabL}, \ref{lem:stableR} and Corollary
\ref{aux43} imply the following $n=\infty$ analogue of
\cite[Theorem~3.34, Corollary~3.36]{Br}.

\begin{theorem}\label{th:u-k}
For $f \in \Z^{m|\infty}_{+}$ we have
\begin{enumerate}
 \item $U_f = \sum_{ \theta \in \{0,1\}^{\#f}} q^{|\theta|}
 K_{\texttt L_\theta (f)}$;

 \item $K_f = \sum_{ \theta \in \mathbb N^{\#f}} (-q)^{|\theta|}
 U_{\texttt L_\theta '(f)}$;

 \item $K_f = \sum_g q^{-|\theta_g|} L_g$, where the sum is over
 all $g \in \Z^{m|\infty}_+$ such that $\texttt R_{\theta_g} (g)=f$ for some (unique)
 $\theta_g  \in \{0,1\}^{\# g}$;

 \item $L_f = \sum (-q)^{-|\theta|} K_g$, summed over
 $g \in \Z^{m|\infty}_+$ and $\theta\in \mathbb N^{\# g}$ with $\texttt R_{\theta}'
 (g)=f$.
\end{enumerate}
\end{theorem}


\begin{remark}  \label{rem:shift}
Let $n\in\N$ and let ${\bf 1}^{m|n} =(1, \cdots, 1| 1,\cdots, 1)$.
For $f, g \in \Zmn$, one has $(f -p{\bf 1}^{m|n}) (n) \le n, (g
-p{\bf 1}^{m|n}) (n) \le n$ for $p$ large enough. One sees that
$$u_{g,f} (q) = u_{g-p{\bf 1}^{m|n}, f -p{\bf 1}^{m|n}} (q),\quad
\ell_{g,f} (q) = \ell_{g-p{\bf 1}^{m|n}, f -p{\bf 1}^{m|n}} (q).$$
The right-hand sides can be computed first by computing in
$\widehat{{\mathcal E}}^{m|\infty}$ (Theorem~\ref{th:u-k}) and
then applying the truncation map ${\mathfrak{Tr}}_n$. Thus, the
structure of the (dual) canonical bases in $\widehat{{\mathcal
E}}^{m|\infty}$ completely controls those in $\widehat{{\mathcal
E}}^{m|n}$ via the truncation map.
\end{remark}

\section{Representation theory of $\mathfrak g\mathfrak l(m|n)$}
\label{sec:superRep}

The main goal of this section is to describe the relations between
tilting modules, Kac modules, and irreducibles of $\mathfrak
g\mathfrak l(m|\infty)$. To do this, we first study the
connections between representations of $\mathfrak g\mathfrak
l(m|n)$ and $\mathfrak g\mathfrak l(m|n+1)$.
\subsection{The categories $\FPn$ and $\mathcal O^{++}_{m|n}$}

For $n \in \N $ the Lie superalgebra $\mathfrak g =\glmn$ is
generated by $e_{ij}$, where $i,j \in I(m|n)$. For $i \in I(m|n)$,
let $\bar{i}=\bar{0}$ if $i>0$ and $i=\bar{1}$ if $i<0$. The
subalgebra $\mathfrak g_{\bar{0}}$ of $\mathfrak g$ is generated
by those $e_{ij}$ such that $\bar{i} +\bar{j} =\bar{0}$ and it is
isomorphic to ${\mathfrak g\mathfrak l}(m) \oplus {\mathfrak
g\mathfrak l}(n)$. Let $\mathfrak h$ be the standard Cartan
subalgebra of $\mathfrak g$ consisting of all diagonal matrices,
$\mathfrak b$ be the standard Borel subalgebra of all upper
triangular matrices, and let $\mathfrak p =\mathfrak g_{\bar{0}}
+\mathfrak b$. The Lie superalgebra $\mathfrak g$ is endowed with
a natural $\Z$-gradation
\begin{equation}\label{glmngrad}
\mathfrak g=\glmn_{-1}\oplus \glmn_0\oplus \glmn_{+1},
\end{equation}
consistent with its $\Z_2$-gradation.  Here $\glmn_{\pm 1}$ is the
subalgebra spanned by the odd positive/negative root vectors. We
let $\glmn_{\le 0}=\glmn_0\oplus\glmn_{-1}$. By means of the
natural inclusion ${\mathfrak{gl}}(m|n) \subseteq
{\mathfrak{gl}}(m|n+1)$ via $I(m|n) \subseteq I(m|n+1)$, we let
${\mathfrak{gl}}(m|\infty):=\lim\limits_{{\longrightarrow}}
{\mathfrak{gl}}(m|n)$.


Let $\{\delta_i| i \in I(m|n)\}$ be the basis of $\mathfrak h^*$
dual to $\{e_{ii}\vert i \in I(m|n)\}$. Let $X_{m|n}$ be the set
of integral weights $\la=\sum_{i \in I(m|n)} \la_i \delta_i$,
$\la_i \in\Z$. A symmetric bilinear form on $\mathfrak h^*$ is
defined by
$$(\delta_i|\delta_j)=-{\rm sgn}(i)\delta_{ij},\qquad i,j \in I(m|n).$$
(Our bilinear form differs from the one in \cite{Br} by a sign.)
Let $X^+_{m|n}$ be the set of all $\la\in X_{m|n}$ such that
$\la_{-m} \ge \cdots \ge \la_{-1}$, $\la_1 \ge \cdots \ge \la_n$.
Such a weight is called {\em dominant}. Let $X^{++}_{m|n}$ be the
set of all $\la\in X^+_{m|n}$ with $\la_n\geq 0$. We may regard an
element $\la$ in $X^{++}_{m|n}$ as an element in $X^{++}_{m|n+1}$
by letting $\la_{n+1} =0$ and let $X^{+}_{m|\infty}
:=\lim\limits_{{\longrightarrow}} X^{++}_{m|n}$.
For $n \in \N \cup \infty$ define
$$\rho = -\sum_{i \in I(m|n)} i \delta_i.$$
Define a bijection
\begin{eqnarray}
 X_{m|n} \longrightarrow \Z^{m|n}, \qquad \la \mapsto f_\la,
\end{eqnarray}
where $f_\la \in \Z^{m|n}$ is given by $f_\la (i) = (\la +\rho |
\delta_i)$, for $i \in I(m|n)$.
This map induces bijections $X^{+}_{m|n} \rightarrow \Zmn$ (for
$n$ possibly infinite), and $X^{++}_{m|n} \rightarrow
\Z_{++}^{m|n}$ (for $n$ finite). Using this bijection we define
the notions such as the degree of atypicality, the $\texttt{L}$
and $\texttt{R}$ operators, $\epsilon$-weight, partial order
$\preccurlyeq$, et cetera, for elements in $X^+_{m|n}$ by
requiring them to be compatible with those defined for elements in
$\Zmn$.

For $\la \in X_{m|n}$, we define the {\em Verma module} and the
{\em Kac module} to be
$$M_n(\la) :=U(\mathfrak g) \otimes_{U (\mathfrak b)}
\C_\la,\quad\text{and}\quad K_n(\la) :=U(\mathfrak g) \otimes_{U
(\mathfrak p)} L^0_n(\la),$$ respectively. The irreducible module
is denoted by $L_n(\la)$. Here $\C_\la$ is the standard
one-dimensional module over $\mathfrak h$ extended trivially to
$\mathfrak b$, and $L^0_n(\la)$ is the irreducible module of
$\mathfrak g_{\bar{0}}$ of highest weight $\la$. Let $[M:
L_{n}(\la)]$ denote the multiplicity of the irreducible module
$L_{n}(\la)$ in a $\mathfrak g\mathfrak l(m|n)$-module $M$. When
$n=\infty$ we will make it a convention to drop the subscript $n$.

For $n\in \mathbb N $, $\FPn$ is the category of
finite-dimensional $\mathfrak{gl}(m|n)$-modules $M$ with
\begin{equation*}
M=\bigoplus_{\g\in X_{m|n}}M_\g,
\end{equation*}
where as usual $M_\g$ denotes the $\g$-weight space of $M$ with
respect to $\mathfrak h$. We denote by $\FPPn$ the full
subcategory of $\FPn$ which consists of modules whose composition
factors are of the form $L_n(\la), \la \in \Xmn^{++}$.

We let $\mathcal O^{++}_{m|\infty}$ be the category of $\mathfrak
h$-semisimple finitely generated $\mathfrak{gl}(m|\infty)$-modules
that are locally finite over $ \mathfrak{gl}(m|N)$, for all finite
$N$, and such that the composition factors are of the form
$L(\la), \la \in X_{m|\infty}^{+}$.

\begin{remark}  \label{rem:shift2}
When $n$ is finite one can pass all the relevant information
between $\FPn$ and $\FPPn$ by tensoring with the one-dimensional
module $L_n(p{\delta}^{m|n})$ for suitable a $p \in\Z$, where
$\delta^{m|n}:= \sum_{i=-m}^{-1}\delta_i -\sum_{i=1}^n \delta_i$.
\end{remark}

\subsection{The truncation functor}

Let $\wt(v)$ denote the weight (or $\delta$-weight) of a weight
vector $v$ in a $\mathfrak g\mathfrak l(m|n)$-module.
\begin{lemma} \label{lem:weight} Let $n$ be finite.
\begin{enumerate}
\item If $\la \in \Xmn^{++}$ and $(\la| \delta_n) <0$ (respectively
$\le 0$), then for every weight vector $v$ in $K_n(\la)$ we have
$(\wt(v)| \delta_n) <0$ (respectively $\le 0$);

\item If $\la \in \Xmn^{++}$ and $(\la| \delta_n) <0$ (respectively
$\le 0$), then for every weight vector $v$ in $L_n(\la)$ we have
$(\wt(v)| \delta_n) <0$ (respectively $\le 0$);

\item For every weight vector $v$ in any module $M \in \FPPn$
we have $(\wt(v)| \delta_n) \le 0$.
\end{enumerate}
\end{lemma}

\begin{proof}
Since $K_n(\la)$ and $L_n(\la)$ are highest weight modules and
every negative root $\alpha$ of $\mathfrak g$ satisfies $(\alpha
|\delta_n) \le 0$, (1) and (2) are clear. Part (3) follows from
(2).
\end{proof}

\begin{corollary}\label{kacstab} For $n\in\N$ and $\la \in \Xmn^{++}$ we have $K_n(\la)\in\FPPn.$
\end{corollary}


\begin{definition}\label{def:trunc}
For $n<n' \le \infty$, the {\em truncation functor}
${\mathfrak{tr}}_{n',n}: \mathcal O_{m|n'}^{++} \longrightarrow
\FPPn$ is defined by sending an object $M$ to
$${\mathfrak{tr}}_{n',n}(M) := \text{span } \{ v \in M \mid (\wt(v)| \delta_{k}) =
0, \text{ for all }n+1 \le k \le n'\}.$$ When $n'$ is clear from
the context we will also write ${\mathfrak{tr}}_{n}$ for
$\mathfrak{tr}_{n',n}$.
\end{definition}
We have a system  of categories $\FPPn$ with a compatible sequence
of functors ${\mathfrak{tr}}_{k,n}$ in the sense that
${\mathfrak{tr}}_{n'',n} = {\mathfrak{tr}}_{n',n} \circ
{\mathfrak{tr}}_{n'',n'}$ for $n''>n'>n$.

%
%
%
\subsection{Kac and irreducible modules in $\FPPn$ and $\Fi$}

\begin{lemma}\label{super}
Let $n$ be finite.
For $\la \in X^{++}_{m|n}$, we have the following natural
inclusions of $\glmn$-modules:
\begin{eqnarray*} L_n(\la) \subseteq L_{n+1} (\la), \qquad
 K_n(\la) \subseteq K_{n+1} (\la).
\end{eqnarray*}
\end{lemma}

\begin{proof}
The eigenvalue of the operator $e_{n+1, n+1} \in \mathfrak{gl}
(m|n+1)$ provides an $\N$-gradation on the Kac module $K_{n+1}
(\la)$, whose degree zero subspace is isomorphic to $K_{n} (\la)$.
Similarly, we have $L_{n}(\la)\subseteq L_{n+1} (\la)$.
\end{proof}
%
%
%
%
%
By Lemma~\ref{super} and the natural inclusions
${\mathfrak g\mathfrak l}(m|1) \subset \cdots \subset \glsuper,$
$\cup_n K_n(\la)$ and $\cup_n L_n(\la)$ are naturally
$\glsuper$-modules. They are direct limits of $\{K_n(\la)\}$ and
$\{L_n(\la)\}$ and isomorphic to $K(\la)$ and $L(\la)$,
respectively. Similarly $\cup_n L_n^0(\la)$ is an irreducible
${\mathfrak g\mathfrak l}(m)\oplus {\mathfrak g\mathfrak
l}(\infty)$-module. The proof of Lemma~\ref{super} implies the
following.



%


\begin{corollary} \label{cor:easyTrun}
Let $n <n' \le \infty$. Let $\la =(\la_{-m}, \cdots, \la_{n'}) \in
X_{m|n'}^{++}$ if $n'$ is finite, and $\la\in X^+_{m|\infty}$
otherwise. Then, for $Y =L$ or $K$ we have
\begin{equation*}{\mathfrak{tr}}_{n',n}(Y_{n'} (\la)) =
 \left\{ \begin{array}{rr}
  Y_n(\la), & \quad \text{if }  \la_{i} =0,\forall i>n, \\
  0, & \quad \text{otherwise.}
  \end{array} \right.
 \end{equation*}
\end{corollary}


\begin{remark}\label{stabkac}
Fix $\la\in X^{++}_{m|n_0}$. Let $n \ge n_0$ and regard $\la \in
X^{++}_{m|n}$. Let $\mathfrak J_{n}$ denote the set of the highest
weights of the composition factors of $K_{n}(\la)$ and $r(n)$
denote the length of a composition series of $K_{n}(\la)$. By
\cite{Br}, $\mathfrak J_{n}$
consists precisely of those weights $\mu$ such that ${\texttt
R}_\theta(\mu)=\la$, for some $\theta\in\{0,1\}^{\#\la}$. As the
operator ${\texttt R}_\theta$ only affects the atypical parts of
$\mu$, we see that $r(n)$ and $\mathfrak J_n$ (with the tail of
zeros in a weight ignored) are independent of $n$, once we have
chosen $n$ large enough so that $\#\la$ is the same when $\la$ is
regarded as an element in $X^{++}_{m|n}$.
\end{remark}

\begin{lemma} \label{lem:series}
Let $\la\in X^{+}_{m|\infty}$ with $\la_i=0$ for $i \ge n_0$, and
regard $\la\in X^{++}_{m|n}$ for $n \ge {n_0}$. Suppose $n_0$ is
chosen so that every $K_n(\la)$, for $n\ge n_0$, has the same
number of composition factors (see Remark~\ref{stabkac}).  Then
${\mathfrak{tr}}_{n',n}$ maps bijectively the set of
Jordan-H\"older series of $K_{n'}(\la)$ onto the set of
Jordan-H\"older series of $K_n(\la)$, for all $n\le n'\le\infty$.
%
In particular $[K(\la): L(\mu)] =[K_{n}(\la): L_{n}(\mu)]$ for
every $n \ge n_0$.

\end{lemma}

\begin{proof}
We denote by ${U(\glmn)}_{+}$, respectively ${U(\glmn)}_{-}$, the
subalgebra of ${U({\mathfrak g\mathfrak l}(m|n))}$ generated by
the positive, respectively negative, root vectors of $\glmn$.

For $n$ finite let $0\subsetneq V^1_n \subsetneq V^2_n \subsetneq
\cdots \subsetneq V^r_n =K_n(\la)$ be a composition series of
$K_n(\la)$ with composition factors $V_n^i/V_n^{i-1}\cong
L_n(\la^i)$. Attached to it we have a set of weight vectors say
$\{v^1,\cdots,v^r\}$ determined in such a way that the irreducible
module $L_n(\la^i)\cong V^{i}_n/V^{i-1}_n$ is generated by the
highest weight vector $v^i+V^{i-1}_n$ of highest weight $\la^i$.
Viewing $v^k \in K_n(\la) \subseteq K_{n+1}(\la)$, we let
$V^i_{n+1} :=\sum_{k=1}^i {U({\mathfrak g\mathfrak
l}(m|n+1))_{-}}v^k$. Now by Lemma \ref{lem:weight},
$e_{l,n+1}v^k=0$ for $l\le n$, and hence $V^i_{n+1}$ is a
${\mathfrak g\mathfrak l}(m|n+1)$-submodule of $K_{n+1}(\la)$.
Also $V^i_{n+1}/V^{i-1}_{n+1}\not=0$,
%
%
and so $0\subsetneq V^1_{n+1} \subsetneq V^2_{n+1} \subsetneq
\cdots \subsetneq V^r_{n+1}=K_{n+1}(\la)$ is a composition series
for $K_{n+1}(\la)$. In this way we have defined a map
$\phi_{n,n+1}$, which takes a composition series of $K_n(\la)$ to
a composition series of $K_{n+1}(\la)$.

By construction $V^i_{n+1}=V^i_n\oplus C^i$, where $C^i$ is a
$\mathfrak{gl}(m|n)$-submodule on which $e_{n+1,n+1}$ acts
non-trivially, and thus $\mathfrak{tr}_{n+1,n}\circ\phi_{n,n+1}$
is the identity.  Now for a composition series of $K_{n+1}(\la)$
let $\{v^1,\cdots,v^r\}$ be defined as before.  It is clear that
$\{v^1,\cdots,v^r\}$ also gives rise to a composition series for
$\mathfrak{tr}_{n+1,n}(K_{n+1}(\la))=K_n(\la)$. By construction of
$\phi_{n,n+1}$ this composition series of $K_n(\la)$ lifts to the
original composition series of $K_{n+1}(\la)$, and hence
$\phi_{n,n+1}\circ\mathfrak{tr}_{n+1,n}$ is the identity.


Now let $V^{i}$ be the $\glsuper$-module $\cup_n V^{i}_n$. Since
$v^i\not\in V^{i-1}$, we have $0\subsetneq V^1 \subsetneq V^2
\subsetneq\cdots \subsetneq V^r =K(\la)$. It is straightforward to
verify that each factor module $V^i /V^{i-1}$ is irreducible and
hence isomorphic to $L(\la^i)$ by construction.  Bijection of
composition series now is established as before.
\end{proof}

\begin{corollary}
The Kac module $K(\la)$ for $\la \in X^{+}_{m|\infty}$ has a
finite composition series.  Furthermore its composition factors
are of the form $L(\mu)$ with $\mu\in X^{+}_{m|\infty}$, and thus
$K(\la)$ belongs to the category $\Fi$.
\end{corollary}
\subsection{Relating tilting modules in $\FPPn$ and $\mathcal O_{m|n+1}^{++}$}

Throughout this subsection we assume that $n$ is finite. An object
$M \in \FPn$ is said to have a {\em Kac flag} if it has a
filtration of $\mathfrak g\mathfrak l (m|n)$-modules:
$$0 =M_0 \subseteq \cdots \subseteq M_r =M$$
such that each $M_i/M_{i-1}$ is isomorphic to $K_n(\la^i)$ for
some $\la^i \in X_{m|n}^+$. We define $(M :K_n(\mu))$ for $\mu\in
X_{m|n}^+$ to be the number of subquotients of a Kac flag of $M$
that are isomorphic to $K_n(\mu)$.

\begin{definition} \label{def:tilt}
The tilting module associated to $\la\in X_{m|n}^+$ in the
category $\FPn$ is the unique indecomposable $\glmn$-module
$U_n(\la)$ satisfying:
\begin{itemize}
\item[(1)] $U_n(\la)$ has a Kac flag with $K_n(\la)$ at the
bottom.

\item[(2)] $\text{Ext}^1(K_n(\mu),U_n(\la))=0$, for all $\mu\in
X^+_{m|n}$.
\end{itemize}
\end{definition}

Let $n\in\N$ and $\la \in \Xmn^{++}$. By Corollary \ref{kacstab}
and Theorem 4.37 (i) of \cite{Br} we see that $U_n(\la)\in\FPPn$ .
Now let $\la =(\la_{-m},\cdots, \la_{n+1}) \in X_{m|n+1}^{++}$ and
let $U_{n+1}(\la)$ be the tilting module of ${\mathfrak g\mathfrak
l}(m|n+1)$ corresponding to $\la$. Let $\#\la$ denote the degree
of atypicality of $\la$. By \cite{Br} $U_{n+1}(\la)$ has a Kac
flag of the form:
\begin{equation}\label{kflagofu}
0 =U_{n+1}^0 \subsetneq U_{n+1}^1\subsetneq
U_{n+1}^2\subsetneq\cdots \subsetneq U_{n+1}^{a_{n+1}}
=U_{n+1}(\la),
\end{equation}
where
$U_{n+1}^{i}/U_{n+1}^{i-1}\cong K_{n+1}(\la^{i,n+1})$
for $i=1,\cdots, a_{n+1}$ and $\la^{1,n+1}=\la$. Furthermore
$a_{n+1} =2^{\#\la}$, and the $\la^{i,n+1}$'s are of the form
$\la^{{\texttt L}_\theta}$, where $\theta =(\theta_1,
\cdots,\theta_{\#\la})\in\{0,1\}^{\#\la}$.

If $(\la|\delta_{n+1}) <0$, then $(\la^{\texttt
L_\theta}|\delta_{n+1}) <0$ for every $\theta$, since
$\la^{\texttt L_\theta} \preccurlyeq \la$. Thus by
Lemma~\ref{lem:weight} ${\mathfrak{tr}}_n( U_{n+1}(\la)) =0$.

Now assume that $(\la|\delta_{n+1}) = 0$. Note that the degree of
atypicality of $\la$, regarded as an element in $X_{m|n}^{++}$, is
either $\#\la-1$ or $\#\la$, depending on whether $\la_{n+1}=0$
affects the atypicality. If it is $\#\la$, then $\la_{n+1}$ is not
an atypical part of $\la$.  Since the operators ${\texttt
L}_\theta$ only affects the atypical part of $\la$, it follows
that $(\la^{\texttt L_\theta} |
\delta_{n+1})=(\la|\delta_{n+1})=0$ for each $\theta$, and the two
sets $\{\la^{1,n},\cdots,\la^{a_n,n}\}$ and
$\{\la^{1,n+1},\cdots,\la^{a_{n+1},n+1}\}$ are identical.

On the other hand if the degree of atypicality is $\#\la-1$, then
for $\theta$ of the form $(1,\theta_2, \cdots, \theta_{\#\la})$,
we have $(\la^{\texttt L_\theta}| \delta_{n+1})<0$, and thus these
$2^{\#\la-1}$ Kac modules will not contribute to
${\mathfrak{tr}}_n(U_{n+1}(\la))$. The remaining $2^{\#\la-1}$ Kac
modules have highest weights satisfying $(\la^{\texttt L_\theta}|
\delta_{n+1}) =0$, and under ${\mathfrak{tr}}_n$ will contribute
$K_n(\la^{\texttt L_\theta})$ in a Kac flag of
${\mathfrak{tr}}_{n}(U_{n+1}(\la))$. But these are precisely the
factors in a Kac flag of $U_n(\la)$.

Denote by $J_n(\la) =\{\la^{1,n},\cdots,\la^{a_{n},n}\}$ the set
of the highest weights of the Kac modules of a Kac flag of
$U_n(\la)$. Summarizing we have the following.

\begin{proposition}\label{stabtilt}
Let $\la\in X^{++}_{m|n+1}$ with $(\la|\delta_{n+1})=0$,  and
$\#\la$ be its degree of atypicality. Then $J_n(\la)$ consists of
those $\mu$'s from $J_{n+1}(\la)$ with $(\mu|\delta_{n+1})=0$.
\end{proposition}

\begin{proposition} \label{tiltflag}
Let $\la\in X^{+}_{m|\infty}$ be such that $\la_i=0$ for $i \ge
n_0$, and regard $\la\in X^{++}_{m|n}$ for $n \ge {n_0}$. Suppose
that $n_0$ is chosen such that the degree of atypicality of $\la$
regarded as an element in $X^+_{m|n}$ is the same for $n\ge n_0$.
Then $J_n(\la)=J_{n_0}(\la)$.
\end{proposition}

\begin{proposition}\label{nestedtilt3}
For $\la\in X^{++}_{m|n+1}$ the map ${\mathfrak{tr}}_n$ sends
$U_{n+1} (\la)$ to $ U_n(\la)$, if $(\la|\delta_{n+1})=0$, and to
$0$, otherwise.
\end{proposition}

\begin{proof}
By \cite{Br} $U_{n+1}(\la)$ has a Kac flag with subquotients
isomorphic to $K_{n+1}(\mu)$, where $\la^{\texttt{L}} \preccurlyeq
\mu \preccurlyeq \la$. If $(\la|\delta_{n+1})>0$, then
Corollary~\ref{cor:easyTrun} implies
${\mathfrak{tr}}_n(U_{n+1}(\la))=0$.

Let $\la\in X^{++}_{m|n+1}$ with $(\la|\delta_{n+1})=0$ and
consider $U_{n+1}(\la)$. Let
$$\partial:=\hf\big{(}\sum_{i<0}e_{ii}-\sum_{j>0}e_{jj}\big{)}.$$
Then $\partial$ equips $\mathfrak{gl}(m|n+1)$ with the
$\Z$-gradation (\ref{glmngrad}), which allows us to regard $\mathcal
O_{m|n+1}^{++}$ as a $\Z$-graded module category with
gradation-preserving morphisms as in \cite{So3,Br2}. By the
construction of \cite{So3} one has a filtration of indecomposable
modules\footnote{It is asserted in \cite{So3} that the modules
$T_{\ge k}$ are indecomposable. A simple proof of this fact can be
read off from the proof of Proposition 3.1 in \cite{So4}. We learned
this proof from W.~Soergel through the anonymous expert, and we
thank both of them.} ($a=\la(\partial):=|\la|$)
\begin{equation}\label{nestedtilt}
K_{n+1}(\la)=T_{\ge a} \subseteq T_{\ge a-1}\subseteq T_{\ge a-2}
\subseteq \cdots\subseteq T_{\ge r} =U_{n+1}(\la).
\end{equation}
Furthermore since according to \cite{Br} tilting modules have
multiplicity-free Kac flags, we have $T_{\ge i-1}/T_{\ge i}\cong
\bigoplus_{|\mu|=i-1} K_{n+1}(\mu)$, where the summation is over
those $\mu$'s in $X^{++}_{m|n+1}$ with ${\rm
Ext}^1(K_{n+1}(\mu),T_{\ge i})\not=0$. Furthermore $T_{\ge i-1}$
contains a non-trivial extension of $K_{n+1}(\mu)$ by $T_{\ge i}$
for every such $\mu$. Applying the truncation functor
${\mathfrak{tr}}_n$ to (\ref{nestedtilt}) we obtain a filtration
\begin{equation}\label{nestedtilt2}
K_{n}(\la)\subseteq S_{\ge a-1}\subseteq S_{\ge
a-2}\subseteq\cdots\subseteq S_{\ge r},
\end{equation}
with $S_{\ge i-1}/S_{\ge i}\cong\bigoplus_{|\mu|=i-1}K_{n}(\mu)$,
and with summation over $\mu$ with $(\mu|\delta_{n+1})=0$.


We claim that the extension of each $K_n(\mu)$ by $S_{\ge i}$ in
(\ref{nestedtilt2}) is non-trivial, for $\mu$ satisfying
$(\mu|\delta_{n+1})=0$, $|\mu|=i-1$, and ${\rm
Ext}^1(K_{n+1}(\mu),T_{\ge i})\not=0$. Suppose it were trivial for
some $\mu$.  Let us denote the extension of $K_{n+1}(\mu)$ by
$T_{\ge i}$ by $T_{\ge i}^\mu$ and set $S^\mu_{\ge
i}={\mathfrak{tr}}_n(T_{\ge i}^\mu)$. Let $w\in T_{\ge i}^\mu$ of
weight $\mu$ be such that its image in $T_{\ge i}^\mu/T_{\ge i}$
generates over ${\mathfrak g\mathfrak l}(m|n+1)_{\le 0}$ a Kac
module $K_{n+1}(\mu)$. Since $S^\mu_{\ge i}$ is a trivial
extension of $K_n(\mu)$ by $S_{\ge i}$, we may assume (by adding
to $w$ an element from $S_{\ge i}$ if necessary) that $e_{ij}w=0$,
for all $-m\le i<j\le n$. By Lemma \ref{lem:weight}
$e_{l,n+1}w=0$, for $l\le n$, and hence $w$ is a genuine
${\mathfrak g\mathfrak l}(m|n+1)$-highest weight vector of highest
weight $\mu$ in $T_{\ge i-1}$.  This implies that the extension
$T^\mu_{\ge i}$ is split, which is a contradiction.

Now Proposition \ref{stabtilt} and Soergel's construction for
$U_n(\la)$ imply that (\ref{nestedtilt2}) is a construction of the
tilting module $U_n(\la)$ and hence $S_{\ge r}\cong U_n(\la)$.
\end{proof}

\begin{corollary} \label{cor:tilt99}
Let $\la\in X^{+}_{m|\infty}$ with $\la_i=0$ for $i \ge n_0$, and
regard $\la\in X^{++}_{m|n}$ for $n \ge {n_0}$. Suppose that $n_0$
is such that $J_n(\la)=J_{n_0}(\la)$, for $n\ge n_0$. Then
$U_{n_0}(\la)\subseteq U_{n}(\la)$ and
$(U_n(\la):K_n(\mu))=(U_{n_0}(\la):K_{n_0}(\mu))$. Furthermore a
Kac flag of $U_n(\la)$, for every $n\ge n_0$, can be chosen to
have the same ordered sequence of weights.
\end{corollary}

\begin{proof} Consider the construction of $U_{n_0}(\la)$ as in
(\ref{nestedtilt}).  We define a total order on the $\mu$'s that
appear in the construction by requiring that $\mu>\nu$ if
$|\mu|>|\nu|$, and among the $\mu$'s with the same $|\mu|$ we
choose an arbitrary total order.  Starting with $K_{n_0}(\la)$,
following this total order, we construct $U_{n_0}(\la)$. This
gives a Kac flag of $U_{n_0}(\la)$.  Regarding the same $\mu$'s as
weights of $\mathfrak{gl}(m|n_0+1)$ we construct in the same
fashion $U_{n_0+1}(\la)$ with the a Kac flag having the same
sequence of weights.
\end{proof}

\subsection{Tilting modules in $\Fi$}

In the same way as in Definition~\ref{def:tilt} one defines
tilting modules $U(\la)$ in the category $\Fi$, where $\la \in
X^{+}_{m|\infty}$. The following lemma is standard using induction
and the long exact sequence.

\begin{lemma}\label{lem:finiteflag}
Let $U$ be a ${\mathfrak g\mathfrak l}(m|\infty)$-module such that
${\rm Ext}^1(K(\mu),U)=0$, for all $\mu$. If $V$ is a ${\mathfrak
g\mathfrak l}(m|\infty)$-module possessing a finite Kac flag, then
we have ${\rm Ext}^1(V,U)=0$.
\end{lemma}

\begin{theorem} \label{tiltInf}
Let $\la\in X^{+}_{m|\infty}$.
\begin{enumerate}
\item There exists a unique (up to isomorphism) tilting module
associated to $\la$ in $\Fi$. Moreover, $U(\la) \cong \cup_n
U_n(\la)$.

\item The functor ${\mathfrak{tr}}_n$ sends $U (\la)$ to
 $ U_n(\la)$ if $(\la|\delta_{n+1})=0$ and to $0$ otherwise.

\item We have $(U(\la):K(\mu))=(U_n(\la):K_n(\mu))$ for $n\gg0$,
and hence $U(\la)$ has a Kac flag consisting of Kac modules of
highest weights $\mu$ with $\la^{\texttt L}\preccurlyeq \mu
\preccurlyeq \la$.
\end{enumerate}
\end{theorem}

\begin{proof}
Corollary \ref{cor:tilt99} allows us to define the ${\mathfrak
g\mathfrak l}(m|\infty)$-module $U(\la)$ to be
$$U(\la)=\cup_{n}U_n(\la).$$
Parts (2) and (3) follow from Corollary \ref{cor:tilt99}.

It remains to prove (1). First by Corollary \ref{cor:tilt99}
$U(\la)$ admits a Kac flag with the same ordered sequence of
weights as $U_n(\la)$, for $n\gg 0$. Now for $\la\in X^{++}_{m|n}$
the highest weights of the Kac modules that appear in a Kac flag
of $U_n(\la)$ lie in $X^{++}_{m|n}$. Also the highest weights of
the composition factors of a Kac module with highest weight lying
in $X^{++}_{m|n}$ also lie in $X^{++}_{m|n}$. Thus $U(\la)$ has a
finite composition series and all its composition factors have
highest weights lying in $X^{+}_{m|\infty}$. Thus $U(\la)$ belongs
to the category $\Fi$.

Suppose $U(\la)=M_1\oplus M_2$ as ${\mathfrak g\mathfrak
l}(m|\infty)$-modules. Choose $n$ so that $(\la_k|\delta_{n})=0$,
for all $\la_k$, where the $\la_k$'s are all the highest weights
of composition factors of $U(\la)$ (and hence of $M_1$ and $M_2$
as well). Now $U_n(\la)={\mathfrak{tr}}_n(M_1)\oplus
{\mathfrak{tr}}_n(M_2)$. But $M_j\not=0$ would imply
${\mathfrak{tr}}_n(M_j)\not=0$, which contradicts the
indecomposability of $U_n(\la)$. Thus, $U(\la)$ is indecomposable.

Given an exact sequence of ${\mathfrak g\mathfrak
l}(m|\infty)$-modules of the form
\begin{equation}\label{Uexact}
0\longrightarrow U(\la)\longrightarrow M\longrightarrow
K(\mu)\longrightarrow 0,
\end{equation}
where $\mu \in X_{m|\infty}^+$. We apply ${\mathfrak{tr}}_n$ to
(\ref{Uexact}) and obtain a split exact sequence
\begin{equation*}
0\longrightarrow U_n(\la)\longrightarrow {\mathfrak{tr}}_n(M)
\longrightarrow K_n(\mu) \longrightarrow 0
\end{equation*}
of $\mathfrak{gl}(m|n)$-modules. Thus, we can choose $v_\mu \in
{\mathfrak{tr}}_n(M)$ of weight $\mu$ such that $e_{ij}v_\mu=0$,
for all $-m\le i<0<j\le n$. Observe that $e_{kl}v_\mu=0$ for $k,l$
such that $k<l$ and $l>n$ by Lemma \ref{lem:weight}. Thus $v_\mu$
is a genuine ${\mathfrak g\mathfrak l}(m|\infty)$-highest weight
vector and hence (\ref{Uexact}) is split. Thus $U(\la)$ is a
tilting module.

A standard Fitting-type argument as in \cite{So3} proves the
uniqueness of $U(\la)$.
%
\end{proof}
%
%
%
%

Denote by $G(\mathcal O^{++}_{m|\infty})$ the Grothendieck group
of $\mathcal O^{++}_{m|\infty}$ and let $G(\mathcal
O^{++}_{m|\infty})_{\mathbb Q}:=G(\mathcal
O^{++}_{m|\infty})\otimes_\Z {\mathbb Q}$. For $M\in\mathcal
O^{++}_{m|\infty}$ let $[M]$ denote the corresponding element in
$G(\mathcal O^{++}_{m|\infty})_{\mathbb Q}$. Let $\mathcal
E^{m|\infty}_{q=1}$ denote the specialization of $\mathcal
E^{m|\infty}$ at $q\to 1$. The topological completion
$\widehat{\mathcal E}^{m|\infty}$ of $\mathcal E^{m|\infty}$
induces a topological completion $\widehat{\mathcal
E}^{m|\infty}_{q=1}$ of $\mathcal E^{m|\infty}_{q=1}$. Using the
bijection between ${\mathcal E}^{m|\infty}_{q=1}$ and $G(\mathcal
O^{++}_{m|\infty})_{\mathbb Q}$, induced by sending $K_{f_\la}$ to
$[K(\la)]$ for $\la\in X^+_{m|\infty}$, we may define a
topological completion $\widehat{G}(\mathcal
O^{++}_{m|\infty})_{\mathbb Q}$ of $G(\mathcal
O^{++}_{m|\infty})_{\mathbb Q}$.

\begin{proposition}  \label{transl=chev}
\begin{enumerate}
\item The linear map $j: \widehat{G}(\mathcal
O^{++}_{m|\infty})_{\mathbb Q} \rightarrow \widehat{\mathcal
E}^{m|\infty}|_{q=1}$ that sends $[K(\la)]$ to $K_{f_\la}(1)$ for
each $\la \in X^+_{m|\infty}$ is an isomorphism of vector spaces.

\item $j$ sends $[U(\la)]$ to $U_{f_\la}(1)$ for each $\la \in
X^+_{m|\infty}$.

\item $j$ sends $[L(\la)]$ to $L_{f_\la}(1)$ for each $\la \in
X^+_{m|\infty}$.
\end{enumerate}
\end{proposition}

\begin{proof}
Clearly $j$ is an isomorphism of vector spaces. Write
$$[K(\la)] =\sum_\mu a_{\la,\mu} [L(\mu)], \qquad [K_n(\la)]
=\sum_\mu a_{\la,\mu;n} [L_n(\mu)],$$
for some $a_{\la,\mu}, a_{\la,\mu;n}\in \Z$.
Lemma~\ref{lem:series} says that $a_{\la,\mu} =a_{\la,\mu;n}$ for
$n\gg0$ relative to $\la,\mu$. Corollary~\ref{aux43} implies that
$b_{\la,\mu} =b_{\la,\mu;n}$ for $n\gg0$ if we write
$$K_{f_\la} =\sum_\mu b_{\la,\mu} L_{f_\mu}, \qquad K_{f_\la^{(n)}}
=\sum_\mu b_{\la,\mu;n} L_{f_\mu^{(n)}}.$$
The finite $n$ analogue of (3) in \cite{Br} says that
$a_{\la,\mu;n} =b_{\la,\mu;n}.$ Thus (3) follows.

In a similar fashion (2) follows from Corollary~\ref{aux43},
Theorem \ref{tiltInf}, and Brundan's finite $n$ analogue of (2).
\end{proof}

\begin{remark}
We note that the operator $\texttt L: \Xmi^+ \rightarrow \Xmi^+$
is not surjective in contrast to the finite $n$ case. For example
the trivial weight does not lie in the image of $\texttt L$.
\end{remark}

\section{Kazhdan-Lusztig polynomials and canonical bases for ${\mathcal E}^{m+n}$}
\label{sec:bar}

The goal of this section is to establish closed formulas for the
(dual) canonical bases and Kazhdan-Lusztig polynomials in the Fock
space ${\mathcal E}^{m+n}$. We begin with some standard results on
canonical bases on tensor and exterior spaces, which seem to be
well known (cf. e.g. \cite{Lu, FKK, Br}).
\subsection{Bases for $ {\mathbb T}^{m+n}$}
\label{sec:tensormodule}

We say that $f \in \Z^{m+n}$ is {\em anti-dominant}, if  $f(-m)
\le \cdots \le f(-1)\le f(1) \le \cdots \le f(n)$. The following
is standard.

\begin{proposition}  \label{prop:classical}
Let $n$ be finite. There exists a unique anti-linear bar map
$\,\bar{} : {\mathbb T}^{m+n} \rightarrow {\mathbb T}^{m+n}$ such
that
\begin{enumerate}
 \item $\overline{\mathcal V_f} =\mathcal V_f$ for all anti-dominant $f \in
 \Z^{m+n}$;
 \item $\overline{XuH} =\overline{X} \overline{u} \overline{H}$
 for all $X \in \mathcal U, u \in {\mathbb T}^{m+n}, H \in
 \mathcal H_{m+n}$.

 \item The bar map is an involution;
 \item $\overline{\mathcal V_f} =\mathcal V_f + (*)$ where $(*)$ is a finite
 $\Z[q,q^{-1}]$-linear combination of $\mathcal V_g$'s for $g <f$.
\end{enumerate}
\end{proposition}

By standard arguments again, Proposition \ref{prop:classical}
implies the existence and uniqueness of the canonical basis
$\{\mathcal T_f\}$ and dual canonical basis $\{\mathcal L_f \}$
for $\mathbb T^{m+n}$.

\begin{theorem}
Let $n$ be finite. There exist unique bases $\{\mathcal T_f\},
\{\mathcal L_f \}$ for ${\mathbb T}^{m+n}$, where ${f \in
\Z^{m+n}}$, such that
\begin{enumerate}
\item $\overline{\mathcal T}_f =\mathcal T_f$ and
$\overline{\mathcal L}_f =\mathcal L_f$;

\item $\mathcal T_f \in  \mathcal V_f + {\sum}_{g \in \Z^{m+n}} q
\Z [q] \mathcal V_g$ and
$\mathcal L_f \in \mathcal  V_f +  {\sum}_{g \in \Z^{m+n}} q^{-1}
\Z [q^{-1}] \mathcal V_g$.

\item $\mathcal T_f = \mathcal V_f + (*)$ and  $\mathcal L_f =
\mathcal V_f + (**)$ where $(*)$ and $(**)$ are finite $\Z
[q,q^{-1}]$-linear combinations of $\mathcal V_g$'s for $g \in
\Z^{m+n}$ with $g <f$.
\end{enumerate}
\end{theorem}

We define $\mathfrak t_{g,f}(q) \in \Z[q], \mathfrak l_{g,f}(q)
\in \Z[q^{-1}]$ associated to $f,g \in \Z^{m+n}$ by
\begin{eqnarray} \label{eq:coeff-tl}
\mathcal T_f =\sum_{g \in \Z^{m+n}} \mathfrak t_{g,f}(q) \mathcal
V_g, \qquad
\mathcal L_f =\sum_{g \in \Z^{m+n}} \mathfrak l_{g,f}(q) \mathcal
V_g.
\end{eqnarray}
Note that $\mathfrak t_{g,f}(q) =\mathfrak l_{g,f}(q)=0$ unless $g
\le f$ and $\mathfrak t_{f,f}(q) =\mathfrak l_{f,f}(q) =1$. We
will refer to these polynomials as {\em Kazhdan-Lusztig
polynomials} (see Proposition~\ref{prop:KL=KL}).

\begin{remark}  \label{rem:duality}
One can show that the matrices $[\mathfrak t_{-f,-g}(q^{-1})]$ and
$[\mathfrak l_{f,g}(q)]$ are inverses of each other. In light of
Proposition~\ref{prop:KL=KL}, this can be regarded as a
reformulation of the inversion formula of the parabolic
Kazhdan-Lusztig polynomials. Alternatively, this follows from
similar arguments as in $\S$2-i, \cite{Br}. Thus
(\ref{eq:coeff-tl}) implies the following duality formulas:
\begin{eqnarray}   \label{eq:duality}
 \mathcal V_f = \sum_{g  \in \Z^{m+n}} \mathfrak
 t_{-f,-g}(q^{-1})\, \mathcal L_g
    = \sum_{g  \in \Z^{m+n}}
 \mathfrak l_{-f,-g}(q^{-1})\, \mathcal T_g, \quad f  \in \Z^{m+n}.
\end{eqnarray}
\end{remark}

\subsection{Bases for ${\mathcal E}^{m+n}$ and $\widehat{\mathcal E}^{m+\infty}$}

For $n\in\N$ let $\mathcal E^{m+n}_+$ denote the subspace of
$\mathcal E^{m+n}$ spanned by elements of the form $\mathcal K_f$,
$f\in\Z^{m+n}_{++}.$
Define the {\em truncation map}
$ {\textsf{Tr}}_{n+1,n}:  {\mathcal E}^{m+(n+1)}_+ \longrightarrow
{\mathcal E}^{m+n}_+ $
by sending $\mathcal K_f$ to $\mathcal K_{f^{(n)}}$ if $f(n+1)
=-n$, and to $0$ otherwise. For $n'>n$, ${\textsf{Tr}}_{n',n}:
{\mathcal E}^{m+n'}_+ \longrightarrow {\mathcal E}^{m+n}_+$ is
defined similarly. This gives rise to ${{\textsf{Tr}}_{n}}:
\mathcal E^{m+\infty}\rightarrow \mathcal E_+^{m+n}$, for all $n$,
which in turn allows us to define a topological completion
$\widehat{\mathcal
E}^{m+\infty}:=\lim\limits_{{\longleftarrow}}\mathcal E_+^{m+n}$,
similarly as in \cite[\S2-d]{Br}. We call the map
${\textsf{Tr}}_n: \widehat{\mathcal E}^{m+\infty}\rightarrow
\mathcal E^{m+n}_+$ also the truncation map.

\begin{lemma} \label{lem:simplemove}
Let $n\in \N \cup \infty$ and let $f,g \in \Z_+^{m+n}$. Then $f
> g$ in the Bruhat ordering if and only if there exists $f^0,\cdots,
f^s \in \Z_+^{m+n}$ such that $f=f^0 > f^1 > \cdots > f^s =g$ and
for each $0 \le a<s$ there exists $i_a<0<j_a$ with $f^{a+1} = (f^a
\cdot \tau_{i_aj_a})^+, f^a(i_a)> f^a(j_a),$ and
$f^a(i_a)\not=f^a(k)$, $\forall k>0$.
\end{lemma}

\begin{proof}[Sketch of a proof]
This follows from an equivalent formulation of the Bruhat ordering
in terms of the $\wt^\epsilon$ combined with an analogous proof of
the reductive analog of Brundan's Lemma~3.42 \cite{Br}.
\end{proof}

For $n\in \N \cup \infty$, a weight $f$ in $\Z_+^{m+n}$ is called
{\em $J$-typical} if it is minimal in the Bruhat ordering among
elements in $\Z_+^{m+n}$. The following proposition can be
established in the same way as a similar statement in
\cite[Theorem~3.5]{Br}. For a finite $n$ let $\widehat{\mathcal
E}^{m+n}:={\mathcal E}^{m+n}$.

\begin{proposition} \label{th:barinf}
Let $n\in \N \cup \infty$. There exists a unique anti-linear bar
map ${}^-: \widehat{\mathcal E}^{m+n} \rightarrow
\widehat{\mathcal E}^{m+n}$ such that
\begin{enumerate}
\item $\overline{\mathcal K_f} = \mathcal K_f$ for all $J$-typical
$f \in \Z_+^{m+n}$;

\item $\overline{X u} = \overline{X} \overline{u}$ for all $X \in
\mathcal U$ and $u \in \widehat{\mathcal E}^{m+n}$.

 \item The bar map is an involution;

\item $\overline{\mathcal K_f} = \mathcal K_f + (*)$ where $(*)$
is (possibly infinite when $n=\infty$) $\Z [q,q^{-1}]$-linear
combination of $\mathcal K_g$'s for $g \in \Z_+^{m+n}$ such that
$g  < f$ in the Bruhat ordering.
\end{enumerate}
\end{proposition}

The next theorem follows from Proposition~\ref{th:barinf}.

\begin{theorem} \label{th:canonical}
Let $n \in \N \cup \infty$. There exist unique topological bases
$\{\mathcal U_f\},$ $ \{\mathcal L_f \}$, where ${f \in
\Z_+^{m+n}}$,  for $\widehat{\mathcal E}^{m+n}$ such that
\begin{enumerate}
\item $\overline{\mathcal U}_f =\mathcal U_f$ and
$\overline{\mathcal L}_f =\mathcal L_f$;

\item $\mathcal U_f \in \mathcal K_f + \widehat{\sum}_{g \in
\Z_+^{m+n}} q \Z [q] \mathcal K_g$ and $\mathcal L_f \in \mathcal
K_f + \widehat{\sum}_{g \in \Z_+^{m+n}} q^{-1} \Z [q^{-1}]
 \mathcal K_g$.

\item $\mathcal U_f = \mathcal K_f + (*)$ and  $\mathcal L_f =
\mathcal K_f + (**)$ where $(*)$ and $(**)$ are (possibly infinite
when $n=\infty$) $\Z [q,q^{-1}]$-linear combinations of $\mathcal
K_g$'s for dominant $g  <f$. For $n$ finite, $(*)$ and $(**)$ are
always finite sums.
\end{enumerate}
\end{theorem}

We define $\mathfrak u_{g,f}(q) \in \Z[q], \mathfrak l_{g,f}(q)
\in \Z[q^{-1}]$ for $f,g \in \Z_+^{m+n}$ by
\begin{eqnarray} \label{ulpoly}
\mathcal U_f =\sum_{g \in \Z_+^{m+n}} \mathfrak u_{g,f}(q)
\mathcal K_g, \qquad
\mathcal L_f =\sum_{g \in \Z_+^{m+n}} \mathfrak l_{g,f}(q)
\mathcal K_g.
\end{eqnarray}
Note that $\mathfrak u_{g,f}(q) =\mathfrak l_{g,f}(q)=0$ unless $g
\le f$ and $\mathfrak u_{f,f}(q) =\mathfrak l_{f,f}(q) =1$. For
reasons that will become clear in Section \ref{sec:KLpoly} we will
refer to these polynomials as {\em (parabolic) Kazhdan-Lusztig
polynomials}.

\begin{remark}  \label{rem:samel}
In the same way as in \cite[(3.7)]{Br} we can show that our
notation of $\mathcal L_f$ and $\mathfrak l_{g,f}(q)$ here is
consistent with the same notation introduced in
Section~\ref{sec:tensormodule}. Also, similar to
\cite[Lemma~3.8]{Br}, we have $\mathcal U_f =\mathcal T_{f \cdot
w_0} \widehat{H}_0$, for $f \in \Z_+^{m+n}.$
\end{remark}

\begin{remark}  \label{dualUK}
By modifying the arguments in \S3-c of \cite{Br} we can introduce
a symmetric bilinear form $\langle\cdot,\cdot\rangle$ on
${\mathcal E}^{m+n}$ such that
%
$\langle\mathcal L_f,\mathcal U_{-g\cdot
w_0}\rangle=\delta_{f,g}$,
which readily implies that the matrices $[\mathfrak{u}_{-f \cdot
w_0,-g \cdot w_0}(q)]$ and $[\mathfrak l_{f ,g}(q^{-1})]$ are
inverses of each other. Equivalently, we have
\begin{equation*}
\mathcal K_f=\sum_{g\in\Z^{m+n}_+} \mathfrak{u}_{-f\cdot w_0,-g
\cdot w_0}(q^{-1})\mathcal L_g =
\sum_{g\in\Z^{m+n}_+}\mathfrak{l}_{-f \cdot w_0,-g \cdot
w_0}(q^{-1})\mathcal U_g, \quad  f\in\Z^{m+n}_+.
\end{equation*}
\end{remark}


\subsection{Degree of $J$-atypicality}

Let $n$ be finite and let $f\in\Z^{m+n}$, $S_{m|n}$-conjugate to
an element $f^+$ in $\Z^{m+n}_+$.

For a pair of integers $(i|j)$, with $-m\le i\le -1$ and $1\le
j\le n$, we define the {\em distance} of $(i|j)$ to be
$d(i|j):=f(i)-f(j).$
Such a pair is called an {\em atypical pair} for $f$, if the
following three conditions are satisfied:
\begin{itemize}
\item[(A1$+$)] $f(i)>f(j)$;

\item[(A2$+$)] For every $k$ with $1\le k\le n$ we have
$f(i)\not=f(k)$;

\item[(A3$+$)] For every $l$ with $-m\le l\le -1$ we have
$f(l)\not=f(j)$.
\end{itemize}

Two atypical pairs $(i_1|j_1)$ and $(i_2 |j_2)$ for $f$ are said
to be {\em disjoint}, if $i_1\not=i_2$ and $j_1\not=j_2$. Two
subsets $A_1$ and $A_2$ of atypical pairs of $f$ are said to be
disjoint, if any two atypical pairs $(i_1|j_1)\in A_1$ and $(i_2
|j_2)\in A_2$ are disjoint.
Let $A^+_f$ denote the set of all atypical pairs of $f$.  For
$k \ge 1$ we define the set $\Sigma^k_f$ recursively as follows:
\begin{align*}
%
&\Sigma^k_f:=\{(i|j)\in A^+_f\vert d(i|j)=k\text{ and }
(i|j)\text{ disjoint from }\bigcup_{s=1}^{k-1}\Sigma^s_f\},\\
&\Sigma^+_f:=\bigcup_{k\ge 1} \Sigma^k_f.
\end{align*}

\begin{definition}
An element in $\Sigma^+_f$ will be called a {\em positive pair} of
$f$.  The {\em degree of $J$-atypicality} of $f$, denoted by
$\divideontimes f$, is defined to be the cardinality of
$\Sigma^+_f$.
\end{definition}
In particular $\divideontimes f\le\text{min}(m,n)$. By
Lemma~\ref{lem:simplemove}, $\divideontimes  f=0$ if and only if
$f$ is minimal in the Bruhat ordering among elements in
$\Z_+^{m+n}$, i.e.~it is {\em $J$-typical}.

\begin{remark}
Our definition of $J$-atypicality was partly inspired by the
definition of Leclerc and Miyachi's function $\psi$ in \cite{LM},
2.5. Their $\psi$ was in turn inspired by Lusztig's notion of
``admissible involutions".
\end{remark}


We denote by $A^-_f$ the collection of all pairs $(i|j)$, $i<0$
and $j>0$, satisfying the following three conditions:
\begin{itemize}
\item[(A1$-$)] $f(i)<f(j)$;

\item[(A2$-$)] For every $k$ with $1\le k\le n$ we have
$f(i)\not=f(k)$;

\item[(A3$-$)] For every $l$ with $-m\le l\le -1$ we have
$f(l)\not=f(j)$.
\end{itemize}
Disjointness of subsets of $A^-_f$ is defined similarly.
We define the set $\Sigma^{-k}_f$ recursively:
\begin{align*}
&\Sigma^{-k}_f:=\{(i|j)\in A^-_f\vert d(i|j)=-k\text{ and }
(i|j)\text{ disjoint from }\cup_{s=1}^{k-1}{\Sigma}^{-s}_f\}, \text{ for } k \ge 1.
\end{align*}
\\
Set $\Sigma^-_f :=\bigcup_{k \ge 1} {\Sigma}^{-k}_f$. Elements in
$\Sigma^-_f$ will be called {\em negative pairs} of $f$.

\subsection{An algorithm for canonical basis}

Recall  $d_i \in \Z^{m+n}$ from Subsection \ref{def:di}. We give a
procedure to reduce the degree of $J$-atypicality in a controlled
manner (see Lemma~\ref{pro321lem}).

\begin{procedure}\label{pro321}  Let $f\in\Z^{m+n}_+$.  We
construct an element $h\in\Z^{m+n}_+$ and a Chevalley generator
$\mathcal X_a$ according to the following procedure.
\begin{itemize}
\item[(0)] Choose $i<0$ such that $i=\text{max}\{l\vert \text{
there exists } $k$\text{ with } (l|k)\in \Sigma^+_f\}$.

\item[(1)] If $i<-1$ and $f(i)-1\not=f(i+1)$ or $i=-1$, then put
${\mathcal X}_a=F_{f(i)-1}$, $h=f- d_i$ and the procedure stops.
Otherwise go to step (2).

\item[(2)] There exists an $s>0$ such that $f(i+1)=f(s)$ (see
Remark \ref{pro321re}). We let ${\mathcal X}_a=E_{f(i)-1}$, $h=f-
d_s$ and the procedure stops.
\end{itemize}
\end{procedure}

\begin{remark}\label{pro321re}
To justify step (2) above, we note that $(i|j)$ is an atypical
pair for some $j>0$.  Thus $f(j) \neq f(t)$ for all $t<0$, and in
addition $f(i)>f(j)$ which implies that $f(i+1) =f(i)-1\ge f(j)$.
If $f(i+1) \neq f(s)$ for all $s>0$, then $(i+1| j)$ is an
atypical pair, and so there exists $t>0$ with
$(i+1|t)\in\Sigma^+_f$, contradicting the choice of $i$ in
step~(0).
\end{remark}

\begin{remark}
Notice that step (1) produces an $h\in\Z^{m+n}_+$ with
$\sum_{i<0}f(i)-1=\sum_{i<0}h(i)$, while step (2) produces an
$h\in\Z^{m+n}_+$ with $\sum_{j>0}f(j)+1=\sum_{j>0}h(j)$. It
follows that a finite number of applications of the procedure
reduces the degree of $J$-atypicality, and thus a finite number of
repeated applications of Procedure \ref{pro321} will produce a
$J$-typical element. Procedure \ref{pro321} may be regarded as a
reductive analogue of Brundan's Procedure~3.20 \cite{Br}.
\end{remark}

We will write $h \rightsquigarrow f$ to denote that $h$ is
obtained from $f$ via Procedure \ref{pro321}.  Thus we have $h_k
\rightsquigarrow h_{k-1} \rightsquigarrow \cdots \rightsquigarrow
h_1 \rightsquigarrow f$, with $h_k$ $J$-typical. Denote the
Chevalley generator corresponding to $h_i \rightsquigarrow
h_{i-1}$ by ${\mathcal X}_{a_i}$. We have the following lemma.

\begin{lemma}\label{pro321lem}
Let $f\in\Z^{m+n}_+$ with $\divideontimes  f>0$ and let $h$ and
${\mathcal X}_a$ be defined as in Procedure \ref{pro321}.
\begin{itemize}
\item[(1)] If $\divideontimes   f=\divideontimes   h$, then
${\mathcal X}_a{\mathcal U}_h={\mathcal U}_f$ and ${\mathcal
X}_a{\mathcal K}_h={\mathcal K}_f$.

\item[(2)] If $\divideontimes   f=\divideontimes  h +1$, then
$h=f- d_i$, for some $i<0$.  In this case we have ${\mathcal
X}_a=F_{f(i)-1}$ and ${\mathcal X}_a{\mathcal U}_h={\mathcal U}_f$
and ${\mathcal X}_a {\mathcal K}_h={\mathcal K}_f + q{\mathcal
K}_{h'}$, where $h'=(h^{<0}|\tilde{e}_{f(i),f(i)-1}h^{>0})$. In
particular $f> h'$.

\item[(3)] Suppose $h_k \rightsquigarrow h_{k-1}\rightsquigarrow
\cdots \rightsquigarrow  h_1=f$ with $h_k$ $J$-typical.  Then we
have ${\mathcal U}_f={\mathcal X}_{a_1}{\mathcal X}_{a_2}\cdots
{\mathcal X}_{a_k}{\mathcal K}_{h_k}$.
\end{itemize}
\end{lemma}

\begin{proof}
As usual we shall write $f=(f^{<0}|f^{>0})$ and
$g=(g^{<0}|g^{>0})$ et cetera.

{}From the description of Procedure \ref{pro321}, there are the
following three cases to be considered depending on the forms of
$f$ and $h$, with $a+1=f(i)$. Here, $\cdots$ denotes entries
different from $a$ and $a+1$.

\vspace{.1cm}

(i) $f=(\cdots,a+1,a,\cdots|\cdots,a,\cdots),\quad
h=(\cdots,a+1,a,\cdots|\cdots,a+1,\cdots);$

(ii) $f=(\cdots,a+1,\cdots|\cdots),\quad
h=(\cdots,a,\cdots|\cdots);$

(iii) $f=(\cdots,a+1,\cdots|\cdots,a,\cdots), \quad
h=(\cdots,a,\cdots|\cdots,a,\cdots).$

\vspace{.1cm}
\noindent We have the following description of
$\Sigma^+_h$ in these cases:
\begin{itemize}
\item[(i)] $\Sigma^+_h=\Big{(}\Sigma^+_f\setminus\{(i|k)\}\Big{)}
\cup \{(i+1|k)\}$, with unique $k>0$ such that
$(i|k)\in\Sigma^+_f$;

\item[(ii)] $\Sigma^+_h=\Sigma^+_f$;

\item[(iii)] $\Sigma^+_h=\Sigma^+_f\setminus\{(i|s)\}$, with unique $s>0$
such that $f(s)=a$.
\end{itemize}
Thus in (i) and (ii) above $\divideontimes h=\divideontimes f$,
while in (iii) $\divideontimes h=\divideontimes f-1$.

The formulas for $\mathcal X_a \mathcal K_h$ in (1) and (2) are
verified case by case using Lemma~\ref{lem:uactsonk}. It remains
to verify the action of $\mathcal X_a$ on $\mathcal U_h$.

In (i), by Procedure \ref{pro321}, $\mathcal X_a$ is of the form
$E_a$. We have
$$\mathcal U_h=\mathcal K_h+ \sum_{g<h} \mathfrak
u_{g,h}(q) \mathcal K_g, \quad \text{for  } \mathfrak u_{g,h}(q)
\in q\Z[q].$$ Observe that any $g$ such that $g<h$  is of the form
\begin{equation*}
g_1=(\cdots,a+1,a,\cdots|\cdots,a+1,\cdots) \text{ or }
g_2=(\cdots,a+1,\cdots|\cdots,a+1,a,\cdots).
\end{equation*}
One checks that
$$\mathcal X_a\mathcal U_h =\mathcal
K_f +\sum_{g_1<h} \mathfrak u_{g_1,h}(q) \mathcal K_{g'_1}
+\sum_{g_2<h} \mathfrak u_{g_2,h}(q) \mathcal K_{g'_2}$$
where
$$g'_1=(g_1^{<0}|\tilde{e}_{a,a+1}g_1^{>0}),
\qquad g'_2=(\tilde{e}_{a,a+1} g_2^{<0}| g_2^{>0}).$$
Note that $g'_1<f$ and $g'_2<f$, for $g_1<h$ and $g_2<h$. Hence
$\mathcal X_a \mathcal U_h$ is of the form $\mathcal K_f
+\sum_{g<f} q\Z[q] \mathcal K_{g}$ and it is bar-invariant since
$$\overline{\mathcal X_a \mathcal U_h} =\overline{\mathcal
X_a}\cdot\overline{\mathcal U_h} = {\mathcal X_a \mathcal U_h}.$$
Thus $\mathcal X_a \mathcal U_h =\mathcal U_f$ by the uniqueness
of the canonical basis elment $\mathcal U_f$.

In case (ii), $\mathcal X_a=F_a$. Now if $g<h$, then $g$ can be
either $(\cdots,a,\cdots|\cdots)$ or $(\cdots|\cdots,a,\cdots)$.
Thus $F_a\mathcal K_{g}=\mathcal K_{g'}$, where $g'$ is either
${(\cdots,a+1,\cdots|\cdots)}$ or ${(\cdots|\cdots,a+1,\cdots)}$.
In either case $g'<f$. From a similar argument as in case~(i), we
conclude that  $\mathcal X_a \mathcal U_h =\mathcal U_f$.
In case (iii), $\mathcal X_a=F_a$ and if $g<h$ then $g$ is
of the same form as $h$.  Now this case is straightforward to check,
thus completing the proof of (1) and (2).
By (1) and (2), ${\mathcal X}_{a_1}{\mathcal
X}_{a_2}\cdots {\mathcal X}_{a_k}{\mathcal K}_{h_k}$ is a
bar-invariant element of the form ${\mathcal K}_f+\sum_{g<f}q\Z[q]
{\mathcal K}_g$ hence equals ${\mathcal U}_f$, proving (3).
\end{proof}

\begin{example}
Consider $f=(2,1,0|3,0,-2)$ with $\divideontimes f=1$.  We will
use Procedure~\ref{pro321} to compute ${\mathcal U}_f$.  We have
the following sequence:
\begin{equation*}
(2,1,-2|3,1,-2)\rightsquigarrow(2,1,-1|3,1,-2)\rightsquigarrow
(2,1,0|3,1,-2)\rightsquigarrow(2,1,0|3,0,-2).
\end{equation*}
Hence, ${\mathcal U}_f= E_0F_{-1}F_{-2}{\mathcal
K}_{(2,1,-2|3,1,-2)}={\mathcal K}_{(2,1,0|3,0,-2)}+q{\mathcal
K}_{(2,0,-2|3,1,0)}.$
\end{example}

Let $\Sigma=\{(i_1|j_1),\cdots,(i_k|j_k)\}$ be a subset of
$\{-m,\cdots,-1\}\times\{1,\cdots,n\}$ such that all the $i_s$'s
and also all the $j_s$'s are distinct. Let $f\in\Z^{m+n}_+$.
Define an element $f_\Sigma\in\Z^{m+n}_+$ obtained from $f$ as
follows. It is enough to specify the values of $f^{<0}_\Sigma$ and
$f^{>0}_\Sigma$. The values of $f_\Sigma^{<0}$ are obtained from
the values of $f^{<0}$ with $f(j_s)$ replacing $f(i_s)$, and the
values of $f_\Sigma^{>0}$ are obtained from the values of $f^{>0}$
with $f(i_s)$ replacing $f(j_s)$, for all $s=1,\cdots,k$.

\begin{theorem}\label{prop:aux84}
Let $f\in\Z^{m+n}_+$.  We have
\begin{equation*}
\mathcal U_f=\sum_{\Sigma\subseteq \Sigma^+_f}q^{|\Sigma|}\mathcal
K_{f_\Sigma}.
\end{equation*}
where the summation is over all subsets $\Sigma$ of $\Sigma^+_f$.
\end{theorem}

\begin{proof}
If $f$ itself is $J$-typical, then $\mathcal K_f=\mathcal U_f$ and
the theorem is true in this case, since $\Sigma^+_f=\emptyset$.

Suppose that $\divideontimes  f>0$. Using the notation of
Procedure~\ref{pro321} and Lemma~\ref{pro321lem} we have $\mathcal
U_f=\mathcal X_{a}\mathcal U_h$. We may assume by induction based
on Procedure~\ref{pro321} that $\mathcal U_h$ is of the form
$\mathcal U_h=\sum_{\Sigma\subseteq
\Sigma^+_h}q^{|\Sigma|}\mathcal K_{h_\Sigma}.$
There are three possibilities for $\mathcal X_a$, where
$a+1=f(i)$.

In the first case $\mathcal X_a=E_{a}$ and the $J$-atypicality is
not changed.  Here we must have
$f=(\cdots,a+1,a,\cdots|\cdots,a,\cdots)$ and
$h=(\cdots,a+1,a,\cdots|\cdots,a+1,\cdots)$.  In this case if
$(i|k)\in\Sigma^+_f$, then
$\Sigma^+_h=\Big{(}\Sigma^+_f\setminus\{(i|k)\}\Big{)}\cup\{(i+1,k)\}$.
Now one checks that
\begin{equation*}
\mathcal U_f=E_{a}\mathcal U_h=\sum_{\Sigma\subseteq
\Sigma^+_h}q^{|\Sigma|}E_{a}\mathcal K_{h_\Sigma}=
\sum_{\Sigma\subseteq \Sigma^+_f}q^{|\Sigma|}\mathcal
K_{f_\Sigma}.
\end{equation*}

In the second case $\mathcal X_a=F_a$ and the $J$-atypicality is
not changed.  In this case $f=(\cdots,a+1,\cdots|\cdots)$ and $h$
is of the form $(\cdots,a,\cdots|\cdots)$, and we have
$\Sigma^+_h=\Sigma^+_f$. Thus
\begin{equation*}
\mathcal U_f=F_{a}\mathcal U_h=\sum_{\Sigma\subseteq
\Sigma^+_h}q^{|\Sigma|}F_{a}\mathcal K_{h_\Sigma}=
\sum_{\Sigma\subseteq \Sigma^+_f}q^{|\Sigma|}\mathcal
K_{f_\Sigma}.
\end{equation*}

Finally suppose that $\mathcal X_a=F_{a}$ and the atypicality is
changed.  In this case $f=(\cdots,a+1,\cdots|\cdots,a,\cdots)$ and
$h=(\cdots,a,\cdots|\cdots,a,\cdots)$.  Let $s>0$ be such that
$f(s)=a$.  Then
$\Sigma^+_h=\Sigma^+_f\setminus\{(i|s)\}$. Thus any $h_\Sigma$ is
of the form $(\cdots,a,\cdots|\cdots,a,\cdots)$ and we have
\begin{align*}
\mathcal U_f=F_{a}\mathcal U_h
=&\sum_{\Sigma\subseteq \Sigma^+_h}q^{|\Sigma|}F_{a}\mathcal
K_{h_\Sigma} \\
=& \sum_{\Sigma\subseteq \Sigma^+_h}q^{|\Sigma|}\mathcal
K_{f_\Sigma}
+\sum_{\Sigma\subseteq \Sigma^+_h}q^{|\Sigma|+1}\mathcal
K_{f_{\Sigma \cup \{(i|s)\}}}
=\sum_{\Sigma\subseteq \Sigma^+_f}q^{|\Sigma|}\mathcal
K_{f_\Sigma}.
\end{align*}
This completes the proof.
\end{proof}

\begin{remark}
Theorem~\ref{prop:aux84} is a generalization of Theorem 3 of
\cite{LM}, where the canonical basis of the irreducible $\mathcal
U$-module generated by the vacuum vector of the Fock space was
computed. The same combinatorics was used in \cite{Lec} to
describe the composition factors of tensor products of evaluation
modules for quantum affine algebras.
\end{remark}

Let $f\in\Z^{m+n}_+$ and let $w_0$ be the longest element in
$S_{m|n}$. From the definitions the negative pairs of $f$ are
precisely the positive pairs of $-f\cdot w_0$ via the
correspondence $$(i|j)\leftrightarrow(-m-i-1|n-j+1).$$ We have the
following consequence of Remark \ref{dualUK} and
Theorem~\ref{prop:aux84}.

\begin{corollary}\label{auxaux111}
Let $f\in\Z^{m+n}_+$.  We have
\begin{equation*}
\mathcal K_f=\sum_{g} q^{-|\Sigma|}\mathcal L_{g},
\end{equation*}
summed over all $g\in\Z^{m+n}_+$ such that there exists
$\Sigma\subseteq \Sigma^-_g$ with $g_\Sigma=f$.
\end{corollary}

\subsection{The operators $\mathbb L$ and $\mathbb R$}
Theorem \ref{prop:aux84} computes $\mathfrak{u}_{g,f}$ and
Corollary \ref{auxaux111} computes the polynomials
$\mathfrak{l}_{g,f}$ upon inversion. In the remainder of this
section we will introduce some combinatorics and establish  the
reductive counterparts of Theorem 3.34 and Corollary 3.36 of
\cite{Br} that allow us to describe this inversion explicitly. We
start with the following lemma which follows from the definition.

\begin{lemma}\label{lem:aux81}
Let $f\in\Z^{m+n}$ be $S_{m|n}$-conjugate to an element in
$\Z^{m+n}_+$ and suppose that $(i|j)$ is an atypical pair of $f$.
Then $(i|j)$ and $\Sigma^+_f$ are not disjoint.
\end{lemma}

\begin{corollary}\label{cor:atypical}  Let $f\in\Z^{m+n}$, $S_{m|n}$-conjugate to an element in
$\Z^{m+n}_+$, and suppose that $(i|j)$ is an atypical pair of $f$.
Then there exists $i'<0$ with $f(i')\le f(i)$ such that
$(i'|j)\in\Sigma^+_f$ or there exists $j'>0$ with $f(j')\ge f(j)$
such that $(i|j')\in\Sigma^+_f$.
\end{corollary}



Next we define the reductive analogue of the ${\texttt
L}$-operators and ${\texttt R}$-operators. Let $f\in\Z^{m+n}$ be
$S_{m|n}$-conjugate to an element in $\Z^{m+n}_+$, with positive
pairs $\Sigma^+_f$ and negative pairs $\Sigma^-_f$. Recalling
$\tau_{ij}\in S_{m+n}$ from Subsection \ref{tauij} we define for
each $i=-m,\cdots,-1$
\begin{eqnarray*}
{\mathbb L}_{i}(f):= \left\{
\begin{array}{ll}
f\cdot\tau_{ij},\quad\text{if there exists $j$ with $(i|j)\in\Sigma^+_f$}, \\
 \emptyset,\quad\text{otherwise}. \\
 \end{array}
 \right.
\end{eqnarray*}
\begin{eqnarray*}
{\mathbb R}_{i}(f):= \left\{
\begin{array}{ll}
f\cdot\tau_{ij},\quad\text{if there exists $j$ with $(i|j)\in\Sigma^-_f$}, \\
 \emptyset,\quad\text{otherwise}. \\
 \end{array}
 \right.
\end{eqnarray*}

Suppose now that $f\in\Z^{m+n}_+$. For
$\theta=(\theta_{-m},\cdots,\theta_{-1})\in\N^{m}$ we define
\begin{align*}
&{\mathbb L_\theta}(f) :=\Big{(}\mathbb L_{-1}^{\theta_{-1}}\circ
\mathbb L_{{-2}}^{\theta_{-2}}\circ\cdots\circ\mathbb
L_{-m}^{\theta_{-m}}(f)\Big{)}^+, &{\mathbb L'_{\theta}}(f)
:=\Big{(}\mathbb L_{-m}^{\theta_{-m}}\circ\cdots\circ\mathbb
L_{-1}^{\theta_{-1}}(f)\Big{)}^+,\\
&{\mathbb R'_{\theta}}(f)
:=\Big{(}\mathbb R_{-1}^{\theta_{-1}}\circ\mathbb
R_{-2}^{\theta_{-2}}\circ\cdots\circ\mathbb
R_{-m}^{\theta_{-m}}(f)\Big{)}^+,
&{\mathbb R_\theta}(f) :=\Big{(}\mathbb
R_{-m}^{\theta_{-m}} \circ\cdots\circ\mathbb
R_{-1}^{\theta_{-1}}(f)\Big{)}^+.
\end{align*}
Write $|\theta|:= \sum_{i=-m}^{-1}\theta_{i}$. Let
$\Sigma^+_f=\{(i_1|j_1),\cdots,(i_k|j_k)\}$, where
$i_1<i_2<\cdots<i_k<0$.
We shall denote by $\mathbb L(f)$ (or $f^{\mathbb L}$) the
expression
\begin{equation*}
\mathbb L(f):=\Big{(}{\mathbb L}_{i_k}\circ{\mathbb
L}_{i_{k-1}}\cdots\circ{\mathbb L}_{i_1}(f)\Big{)}^+.
\end{equation*}
Similarly, if $\Sigma^-_f=\{(i_1|j_1),\cdots,(i_k|j_k)\}$, with
$i_1<i_2<\cdots<i_k<0$, we shall write $\mathbb R(f)$ (or
$f^{\mathbb R}$) for the expression
\begin{equation*}
\mathbb R(f):=\Big{(}{\mathbb R}_{i_1}\circ{\mathbb
R}_{i_{2}}\cdots\circ{\mathbb R}_{i_k}(f)\Big{)}^+.
\end{equation*}

The following follows from the definitions.

\begin{lemma}\label{aux112}
Let $f\in\Z^{m+n}_+$ and let $w_0$ be the longest element in
$S_{m|n}$. For $\theta=(\theta_{-m},\cdots,\theta_{-1})\in\N^{m}$
let $\varphi=(\theta_{-1},\cdots,\theta_{-m})$.  We have
$$-{\mathbb L}_\theta(-f\cdot w_0)\cdot w_0={\mathbb R}_\varphi(f),\quad
-{\mathbb L}'_\theta(-f\cdot w_0)\cdot w_0={\mathbb
R}'_\varphi(f).$$
\end{lemma}

\subsection{Transition matrices between monomial and
canonical bases in ${\mathcal E}^{m+n}$}
We start with several
simple lemmas whose proofs are straightforward.

\begin{lemma}\label{pingpong}
Let $f\in\Z^{m+n}$ be $S_{m|n}$-conjugate to an element in
$\Z^{m+n}_+$. Suppose that $(i_0|j_0)\in\Sigma^+_f$.
\begin{itemize}
\item[(i)] If there exists $b>0$ such that $f(b)\not\in\{f(i)\vert
i<0\}$ and $f(i_0)>f(b)>f(j_0)$, then there exists $a<0$ with
$f(i_0)>f(a)>f(b)$ and $(a|b)\in\Sigma^+_f$.

\item[(ii)] If there exists $a<0$ with $f(a)\not\in\{f(j)\vert
j>0\}$ such that $f(i_0)>f(a)>f(j_0)$, then there exists $b>0$
with $f(a)>f(b)>f(j_0)$ and $(a|b)\in\Sigma^+_f$.
\end{itemize}
\end{lemma}

\begin{lemma}\label{lem:aux82}
Let $f\in\Z^{m+n}$ be $S_{m|n}$-conjugate to an element in
$\Z^{m+n}_+$, and let $(i|k),(j|l)\in\Sigma^+_f$ with $f(i)>f(j)$.
Then $(j|l)\in\Sigma^+_{\mathbb L_{i}(f)}$.
\end{lemma}

\begin{lemma}\label{lem:aux83}
Let $f\in\Z^{m+n}$ be $S_{m|n}$-conjugate to an element in
$\Z^{m+n}_+$, and let $(i|k)\in\Sigma^+_f$. Suppose $i<0, j<0$ and
$f(i)>f(j)$, and there is no $l>0$ such that $(j|l)\in\Sigma^+_f$.
Then there is no $l$ such that $(j|l)\in\Sigma^+_{\mathbb
L_i(f)}$.
\end{lemma}

\begin{proof} The lemma is obvious if there exists $b>0$ such that
$f(j)=f(b)$.  Hence we may assume that this is not the case.

Suppose there exists an $l>0$ such that $f(l)\not=f(c)$, for every
$c<0$, and $f(j)>f(l)$. Then there exists an $a<0$ with
$f(j)>f(a)>f(l)$ and $(a|l)\in\Sigma^+_f$.  But then by Lemma
\ref{lem:aux82} we have $(a|l)\in\Sigma^+_{\mathbb L_i(f)}$, and
hence $(j|l)\not\in\Sigma^+_{\mathbb L_i(f)}$.
\end{proof}

\begin{theorem}\label{canonical-reductive}
Let $f\in\Z^{m+n}_+$.  We have
\begin{itemize}
\item[(1)] $\mathcal
U_f=\sum_{\theta\in\{0,1\}^m}q^{|\theta|}\mathcal K_{\mathbb
L_\theta(f)}$; \item[(2)] $\mathcal
K_f=\sum_{\theta\in\N^m}(-q)^{|\theta|}\mathcal U_{\mathbb
L'_\theta(f)}$.
\end{itemize}
\end{theorem}

\begin{proof} Part (1) follows from Theorem~\ref{prop:aux84},
Lemmas \ref{lem:aux82} and \ref{lem:aux83}.

Part (2) can be proved using  (1) following the strategy of the
proof of Theorem 3.34  (ii) in \cite{Br} as follows. We let
$M^{m+n}$ be the free $\Z[q,q^{-1}]$-module with basis $\{[f]\}$,
where the $f$'s are all the $S_{m|n}$-conjugates of elements in
$\Z^{m+n}_+$.  We define linear maps $\mathcal U,\mathcal
K:M^{m+n}\rightarrow\mathcal E^{m+n}$ by
$\mathcal K([f])=\mathcal K_{f^+},  \mathcal U([f])=\mathcal
U_{f^+}.$
Let $\mathcal K^{-1},\mathcal U^{-1}$ denote the right inverse
functions of $\mathcal K$ and $\mathcal U$, respectively, given by
$\mathcal K^{-1}(\mathcal K_{f})=[f],  \mathcal U^{-1}(\mathcal
U_{f})=[f].$
Now for $-m\le i\le -1$ define $\la_i:M^{m+n}\rightarrow M^{m+n}$ by
\begin{eqnarray}
{\mathbb \la}_{i}([f]):= \left\{
\begin{array}{ll}
[f\cdot\tau_{ij}],\quad\text{if there exists $j$ with $(i|j)\in\Sigma^+_f$}, \\
 0,\quad\text{otherwise}. \\
 \end{array}
 \right.
\end{eqnarray}
Consider the maps
\begin{align*}
&\mathcal P:=\mathcal K\circ\Big{(}(1+q\la_{-1})\cdots(1+q\la_{-m})
\Big{)}\circ \mathcal U^{-1}:\mathcal E^{m+n}\rightarrow\mathcal
E^{m+n},\\
&\mathcal Q:=\mathcal U\circ\Big{(}  \frac{1}{1+q\la_{-m}} \cdots
\frac{1}{1+q\la_{-1}}  \Big{)}\circ \mathcal K^{-1}:\mathcal
E^{m+n}\rightarrow\mathcal E^{m+n}.
\end{align*}
Note that every $\la_i$ is nilpotent, so the expression ${1+q\la_i}$
is invertible.

Now  (1) says that $\mathcal P(\mathcal U_f)=\mathcal U_f$ and hence
$\mathcal P$ is the identity map.  Part (2) is equivalent to saying
that $\mathcal Q$ is the identity map.  So it suffices to show that
$\mathcal P\circ\mathcal Q$ is the identity map, which amounts to
removing $\mathcal U^{-1}\circ \mathcal U$ in $\mathcal
P\circ\mathcal Q$. Now this would follow, once we can show that,
given $f\in\Z^{m+n}_+$ and a nonzero summand $[g]$ of
$\frac{1}{1+q\la_{-m}} \cdots \frac{1}{1+q\la_{-1}} [f]$,
\begin{equation}\label{aux822}
\mathcal K\circ (1+q\la_{-1})\cdots(1+q\la_{-m})  [g]= \mathcal
K\circ (1+q\la_{-1})\cdots(1+q\la_{-m}) [g^+].
\end{equation}
We make the following three observations for [g] a summand of
$\frac{1}{1+q\la_{-m}} \cdots \frac{1}{1+q\la_{-1}}[f]$, for
$f\in\Z^{m+n}_+$:
\begin{itemize}
\item[(i)] If $g$ is of the form $(\cdots g(i)\cdots
g(j)\cdots|\cdots)$ with $g(i)<g(j)$, $i<j$ and $(i|k),
(j|l)\in\Sigma^+_g$, for some $k,l>0$, then $g(l)>g(i)$.

\item[(ii)] For $g$ having
the form of (i), $\mathbb L_i\mathbb L_j(g)$ and $\mathbb L_j\mathbb
L_i(g)$ are $S_{m|n}$-conjugates.

\item [(iii)] For $g$ of the form (i) and $-m\le i_1<i_2<\cdots <i_k<j<0$ with
$i\not=i_s$ for $s=1,\cdots,k$, the element $\mathbb
L_{i_k}\cdots\mathbb L_{i_1}(g)$ is also of the form (i).
\end{itemize}
Now let $g$ be such a summand.  For $-m\le i_1<i_2<\cdots <i_k<0$
consider $h=\mathbb L_{i_k}\cdots\mathbb L_{i_1}(g)$.  Of course
there exist distinct $-m\le j_1,\cdots,j_k<0$ such that $\mathbb
L_{j_k}\cdots\mathbb L_{j_1}(g^+)$ is $S_{m|n}$-conjugate to $h$. In
order to establish (\ref{aux822}) we only need to show that if we
rearrange the $j_1,\cdots,j_k$ in decreasing order and apply the
corresponding $\mathbb L$-operators to $g^+$ we get an element that
is also $S_{m|n}$-conjugate to $h$.  This follows from (ii) and
(iii).  Below is an example to illustrate this. Suppose
\begin{align*}
&g=(\cdots,g(i_1),\cdots,g(i_2),\cdots,g(i_3)\cdots|\cdots),\\
&g^+=(\cdots,g(i_3),\cdots,g(i_2),\cdots,g(i_1),\cdots|\cdots),
\end{align*}
with $g(i_s)=g^+(j_s)$ with $s=1,2,3$.  We write $h\sim h'$, if $h$
and $h'$ are $S_{m|n}$-conjugates. We have by (ii) and (iii)
\begin{align*}
\mathbb L_{i_3}\mathbb L_{i_2}\mathbb L_{i_1}g\sim \mathbb
L_{i_3}\mathbb L_{i_1}\mathbb L_{i_2}g\sim \mathbb L_{i_1}\mathbb
L_{i_3}\mathbb L_{i_2}g\sim \mathbb L_{i_1}\mathbb L_{i_2}\mathbb
L_{i_3}g.
\end{align*}
But $\mathbb L_{i_1}\mathbb L_{i_2}\mathbb L_{i_3}g$ is
$S_{m|n}$-conjugate to $\mathbb L_{j_3}\mathbb L_{j_2}\mathbb
L_{j_1}g^+$.
\end{proof}

\begin{remark}
The formula for $\mathcal U_f$ in Theorem
\ref{canonical-reductive} (or in Theorem~\ref{prop:aux84}) implies
that $\mathcal U_f$ is a sum of precisely $2^{\divideontimes f}$
monomial basis elements.
\end{remark}

\subsection{Transition matrices with dual
canonical basis in ${\mathcal E}^{m+n}$}

\begin{corollary}  \label{cor:redKL}
Let $f\in\Z^{m+n}_+$.  We have
\begin{itemize}
\item[(1)] $\mathcal K_f=\sum_{\mathbb R_\theta(g)=f}
q^{-|\theta|}\mathcal L_{g}$,

\item[(2)] $\mathcal L_f=\sum_{\mathbb R'_\theta(g)=f}
(-q)^{-|\theta|}\mathcal K_{g},$
\end{itemize}
where the summation in (1) is over all $g\in\Z^{m+n}_+$ such that
there exists $\theta\in\{0,1\}^m$ with ${\mathbb R}_\theta(g)=f$,
and the summation in (2) is over all $g\in\Z^{m+n}_+$ and
$\theta\in\N^m$ such that ${\mathbb R}'_\theta(g)=f$.
\end{corollary}

\begin{proof}
In order to establish part (1) we need to compute ${\mathfrak
u}_{-fw_0,-gw_0}(q^{-1})$ according to Remark \ref{dualUK}. By
Theorem \ref{canonical-reductive}~(1) it is equal to
$q^{-|\theta|}$ if $-f\cdot w_0=\mathbb L_\theta(-g\cdot w_0)$,
where $\theta\in\{0,1\}^m$, which is equivalent to saying that
$f=-\mathbb L_\theta(-g\cdot w_0)\cdot w_0$. Now by Lemma
\ref{aux112} we have $-\mathbb L_\theta(-g\cdot w_0)\cdot
w_0={\mathbb R}_\varphi(g)$, where we recall that the expression
$\varphi$ means $\theta$ in reversed order. Thus we obtain
$\mathbb R_\varphi(g)=f$, as desired.

Part (2) is proved in a similar way using the second identity in
Remark \ref{dualUK} together with Theorem
\ref{canonical-reductive}~(2) and Lemma \ref{aux112}.
\end{proof}

\begin{corollary}
Both $\mathfrak u_{g,f}(q)$ and $\mathfrak l_{g,f}(-q^{-1})$ are
polynomials in $q$ with positive integer coefficients. More
explicitly,
\begin{itemize}
\item[(1)] $\mathfrak u_{g,f}(q)=q^{|\theta|}$, if $g=\mathbb
L_{\theta}(f)$ for some $\theta\in\{0,1\}^m$, and ${\mathfrak
u}_{g,f}(q)=0$ otherwise; \item[(2)] $\mathfrak
l_{g,f}(-q^{-1})=\sum_{\theta}q^{|\theta|}$, where the sum is over
all $\theta\in\N^{m}$ with $\mathbb R'_\theta(g)=f$.
\end{itemize}
\end{corollary}

\subsection{Transition matrices for bases in $\widehat{\mathcal E}^{m+\infty}$}

\begin{proposition} \label{commutativity-reduct}
The truncation map ${\textsf{Tr}}_{n+1, n}:  {\mathcal
E}^{m+n+1}_+\rightarrow{\mathcal E}^{m+n}_+$ commutes with the
bar-involution.
\end{proposition}

\begin{proof}
The proposition can be proved by an almost identical argument as
Proposition~\ref{commutativity2}.  We will not give the details
except for pointing out that we now work over the space ${\mathcal
E}^{m+n}\otimes\mathbb V$ and we have the following counterparts
of (\ref{quasiR}) and (\ref{R-matrix}):
\begin{eqnarray*}
\Theta^{(m+n, 1)}&=& 1\otimes 1 + (q-q^{-1}) \sum_{a, b: a<b}
\tilde{E}_{a b}\otimes e_{b a},
\\
\overline{{\mathcal{K}}_{f^{(n)}}\otimes v_{f(n+1)}}
&=&\overline{{\mathcal{K}}_{f^{(n)}}}\otimes v_{f(n+1)} \nonumber\\
&& + (q-q^{-1}) \sum_{b>f(n+1)} \tilde{E}_{f(n+1),\,  b}
\overline{{\mathcal{K}}_{f^{(n)}}}\otimes v_b .\label{R-matrix1}
\end{eqnarray*}
\end{proof}

Similarly as in Subsection~\ref{trunct-bar},
Proposition~\ref{commutativity-reduct} implies the following.

\begin{corollary}
\begin{enumerate}
\item $\textsf{Tr}_{n+1,n}$ sends $\mathcal U_f$ (respectively
$\mathcal L_f$) to $\mathcal U_{f^{(n)}}$ (respectively $\mathcal
L_{f^{(n)}}$) if $f(n+1) =-n$ and to $0$ otherwise.

\item
For $f, g \in \Z_{++}^{m+(n+1)}$ such that $f(n+1) =g(n+1)= -n$,
we have
$$\mathfrak u_{g,f} (q) =\mathfrak u_{g^{(n)}, f^{(n)}} (q),\quad
\mathfrak l_{g,f} (q) =\mathfrak l_{g^{(n)}, f^{(n)}} (q).$$
\end{enumerate}
\end{corollary}

\begin{remark} \label{aux1000}
By Proposition~\ref{commutativity-reduct} and stabilizing its
finite $n$ counterparts (Theorem~\ref{canonical-reductive},
Corollary~\ref{cor:redKL}) we see that
Theorem~\ref{canonical-reductive} and Corollary~\ref{cor:redKL}
remain to be valid for $n=\infty$.
\end{remark}

%
%
%
\section{Representation theory of $\mathfrak g\mathfrak l(m+n)$}
\label{sec:reductiveRep}
\subsection{The categories $\mathcal O_{m+n}$, $\OPn$ and $\OPPn$}
\label{sec:categoryO}

Let $n \in \N $. We shall think of ${\mathfrak g\mathfrak l}(m+n)$
as the Lie algebra of complex matrices whose rows and columns are
parameterized by $I(m|n)$. We denote by $\mathfrak h_c$
(respectively $\mathfrak b_c$) the standard Cartan (respectively
Borel) subalgebra of ${\mathfrak g\mathfrak l}(m+n)$, which
consists of all diagonal (respectively upper triangular) matrices.
Let $\{\delta'_i, i \in I(m|n)\}$ be the basis of $\mathfrak
h_c^*$ dual to $\{e_{ii}, i \in I(m|n)\}$. Set $\mathfrak g_0=
{\mathfrak g\mathfrak l}(m)\oplus {\mathfrak g\mathfrak l}(n)$.
Let $\mathfrak q$ be a parabolic subalgebra $\mathfrak g_{0}
+\mathfrak b_c$. Similarly as in Section \ref{sec:superRep} the
Lie algebra $\mathfrak{gl}(m+\infty)$ is defined as
$\lim\limits_{\stackrel{\longrightarrow}{n}}\mathfrak{gl}(m+n)$.

Define the symmetric bilinear form $(\cdot\vert\cdot)_c$ on
$\mathfrak h_c^*$ by
$$(\delta'_i\vert\delta'_j)_c=\delta_{ij}, \qquad i,j \in I(m|n),$$
and let
$$\rho_c =-\sum_{i=-m}^{-1}i\delta'_i+\sum_{j=1}^n(1-j)\delta'_j.$$

Let $X_{m+n}$ be the set of integral weights $\la=\sum_{i \in
I(m|n)} \la_i \delta'_i$, $\la_i \in\Z$. let $X^+_{m+n}$ be the
set of all $\la \in X_{m+n}$ such that $\la_{-m} \ge \cdots \ge
\la_{-1}$, $\la_1 \ge \cdots \ge \la_n$ (such a weight will be
referred to as {\em dominant}), and let $X^{++}_{m+n}$ be the set
of all $\la \in X^+_{m+n}$ such that $\la_n \geq 0$. We may regard
an element $\la$ in $X^{++}_{m+n}$ as an element in $X^{++}_{m+
n+1}$ by letting $\la_{n+1} =0$. Analogously we denote by
$X^{+}_{m+\infty}$ the set of all dominant integral weights $\la
=\sum_{i \in I(m|\infty)} \la_i \delta'_i \in X_{m+\infty}$ such
that $\la_{-m} \ge \cdots \ge \la_{-1}$, $\la_1 \ge \la_2 \ge
\cdots$, and $\la_i =0$ for $i \gg0$.

Given $\la \in X_{m+n}^+$, $n\in\N\cup\infty$, we denote the Verma
and generalized Verma modules by
$${\mathcal V}_n(\la):= U({\mathfrak g\mathfrak l}(m+n))
\otimes_{U(\mathfrak b_c)} \C_\la,\quad{\mathcal K}_n(\la)
:=U({\mathfrak g\mathfrak l}(m+n)) \otimes_{U(\mathfrak q)}
L_n^0(\la), $$
respectively. Here $L_n^0(\la)$ is extended trivially to a
$\mathfrak q$-module, and $\C_\la$ is the standard one-dimensional
$\mathfrak h_c$-module extended trivially to a $\mathfrak
b_c$-module. Let $\mathcal L_n(\la)$ be the irreducible
${\mathfrak g\mathfrak l}(m+n)$-module of highest weight $\la$.

Let $n \in \N \cup \infty$. Define a bijection
\begin{eqnarray}
 X_{m+n} \longrightarrow \Z^{m+n}, \qquad \la \mapsto f_\la,
\end{eqnarray}
where $f_\la \in \Z^{m+n}$ is given by
$f_\la (i) = (\la +\rho_c \vert \delta'_i)_c$ for $i \in
I(m|n).$
This map induces bijections $X^{+}_{m+n} \rightarrow \Z_+^{m+n}$
(for $n$ possibly infinite) and $X^{++}_{m+n} \rightarrow
\Z_{++}^{m+n}$ (for finite $n$). Notions, such as Bruhat ordering
and degree of $J$-atypicality etc., are defined on elements of
$X_{m+n}$ in a way that is compatible with the ones defined on
elements of $\Z^{m+n}$.

Let $n\in\N$. Denote by $\mathcal O_{m+n}$ the category of
finitely generated $\glmpn$-modules $M$ that are locally finite
over $\mathfrak b_c$, semisimple over $\mathfrak h_c$ and
\begin{equation*}
M=\bigoplus_{\g\in X_{m+n}}M_\g,
\end{equation*}
where as usual $M_\g$ denotes the $\g$-weight space of $M$ with
respect to $\mathfrak h_c$. The objects in $\mathcal O_{m+n}$
include the Verma module $\mathcal V_n(\la)$ and the irreducible
module $\mathcal L_n(\la)$ for $\la \in X_{m+n}$. Denote by $\OPn$
the full subcategory of finitely generated $\glmpn$-modules $M$
which are locally finite over $\mathfrak q$. The generalized Verma
module $\mathcal K_n(\la)$ and the irreducible module $\mathcal
L_n(\la)$ for $\la \in X_{m+n}^+$ belong to $\OPn$. Denote by
$\OPPn$ the full subcategory of $\OPn$ which consists of
$\glmpn$-modules $M$ whose composition factors are isomorphic to
$\mathcal L_n(\la)$ with $\la \in X^{++}_{m+n}$.

Given $M\in\OPn$, we endow the restricted dual $M^*$ with the
usual $\mathfrak g\mathfrak l(m+n)$-module structure. Further
twisting the ${\mathfrak g\mathfrak l}(m+n)$-action on $M^*$ by
the automorphism given by the negative transpose of $\mathfrak
g\mathfrak l(m+n)$, we obtain another $\mathfrak g$-module denoted
by $M^\tau$. For $M$ with finite-dimensional weight spaces we have
$(M^\tau)^\tau \cong M$.

The category  $\Oi$ consists of finitely generated
$\mathfrak{gl}(m+\infty)$-modules that are locally finite over
$\mathfrak q\cap \mathfrak{gl}(m+N)$, for every $N$, and such that
the composition factors are of the form $\mathcal L(\la)$, $\la\in
X_{m+\infty}^+$.

A $\glmpn$-module $M$ is said to have {\em a generalized Verma
flag} if there exists a filtration of $\glmpn$-modules
$0 =M_0 \subseteq M_1 \subseteq \cdots \subseteq M_r =M$
such that $M_i/M_{i-1}$ is isomorphic to some generalized Verma
module for each $i$. Denote by $(M: \mathcal K_n(\la))$ the number
of subquotients of $M$ which are isomorphic to $\mathcal
K_n(\la)$.

Tilting modules in $\mathcal O_{m+n}$ and $\OPn$ are defined
similarly as before, cf.~Definition~\ref{def:tilt}.  The role of
standard modules there is played by Kac modules, while here it is
played by generalized Verma modules. For a finite $n$, the notion
of the tilting modules in $\mathcal O_{m+n}$ and $\OPn$ first
appeared in \cite{CoI} in a somewhat different formulation (see
\cite{So3} for generalizations based on the work of Ringel and
others; also cf. \cite{Br2}).

The character formula of the tilting module $\mathcal T_n (\mu)$
in $\mathcal O_{m+n}$ is given as follows \cite{So3} (also see
\cite{Br2}).
\begin{eqnarray} \label{eq:tiltT}
(\mathcal T_n (\mu): \mathcal V_n(\nu))
 =[\mathcal V_n (-\nu  -2 \rho_c) : \mathcal L_n (-\mu -2 \rho_c)],
 \quad \mu, \nu \in X_{m+n}.
\end{eqnarray}

Let
$\beta =-n(\delta_{-m}'+\cdots +\delta_{-1}')+
m(\delta_{1}'+\cdots +\delta_{n}')$
be the sum of all the negative roots of $\mathfrak g\mathfrak
l(m+n)$ that are not roots of $\mathfrak g\mathfrak l(m)\oplus
\mathfrak g\mathfrak l(n)$. Then Theorem~6.7 (and Lemma~7.4) of
\cite{So3} imply the following character formula for the tilting
module $\mathcal U_n (\mu)$ in $\OPn$.
\begin{eqnarray} \label{eq:tiltU}
(\mathcal U_n (\mu): \mathcal  K_n(\nu))
 =[\mathcal K_n (\beta -w_0 \nu) : \mathcal L_n(\beta -w_0\mu)],
 \quad \mu, \nu \in X_{m+n}^+.
\end{eqnarray}
The formula is best understood by the following formula
\begin{eqnarray}  \label{invol}
f_{\beta -w_0\mu} = (m-n+1) {\bf 1}^{m+n} -f_\mu\cdot w_0,
\end{eqnarray}
where ${\bf 1}^{m+n} =(1,\cdots,1|1,\cdots,1)\in \Z^{m+n}$. So the
difference of the parameters on both sides of (\ref{eq:tiltU}) is
essentially given by the symmetry $f \mapsto -f\cdot w_0$ (similar
phenomenon has been observed in (4.17) of \cite{Br} for $\mathfrak
g\mathfrak l(m|n)$.)

\subsection{The Jantzen irreducibility criterion}

The following is obtained by applying Jantzen's irreducibility
criterion (M+) \cite[Satz 4]{J} to ${\mathcal K}_n(\la)$.

\begin{proposition}\label{jantzen}
Let $n \in \N$. The $\mathfrak g\mathfrak l (m+n)$-module
${\mathcal K}_n(\la)$ is irreducible if and only if for any
$\beta$ of the form $\delta'_i-\delta'_j$, $i<0< j$, with
$(\la+\rho_c\vert\beta)_c$ being a positive integer, there exists
an $\alpha$ of the form $\delta'_i-\delta'_l$, $0< l$, or
$\delta'_k-\delta'_j$, $k<0$, with $(\la+\rho_c\vert\alpha)_c=0$.
\end{proposition}

Jantzen's criterion can be simplified for $\gl$.
\begin{proposition}\label{prop:jantzen}
The $\gl$-module ${\mathcal K}(\la)$ is irreducible if and only if
for every $-m \le i <0$, there exists an $\ell >0$ with
$(\la+\rho_c\vert\delta'_i -\delta'_\ell)_c=0$.
\end{proposition}

\begin{proof}
As Jantzen's criterion is obtained by the non-degeneracy of the
contravariant bilinear form on the generalized Verma module
${\mathcal K}(\la)$, it works for $n=\infty$ as well. Let us call
the criterion for $n=\infty$ in Proposition~\ref{jantzen}
Condition (A) and the one in Proposition~\ref{prop:jantzen}
Condition (B). Apparently Condition (B) implies (A).

On the other hand, assume that (A) is satisfied. Let $i$ be such
that $-m \le i<0$. Since $(\la+\rho_c\vert
\delta'_k-\delta'_{J+1})_c
>(\la+\rho_c\vert \delta'_k-\delta'_J)_c$, we can take $J\gg0$
such that $(\la+\rho_c\vert \delta'_k-\delta'_J)_c >0$ for every
$-m \le k \le -1$. In particular, for the root $ {\beta}
=\delta'_i-\delta'_J$ the number $(\la +\rho_c \vert {\beta})_c$
is a positive integer. Applying Condition (A) to ${\beta}$ now
gives us Condition (B).
\end{proof}

\subsection{Kazhdan-Lusztig polynomials and (dual) canonical bases}
\label{sec:KLpoly}

In \cite{KL} the Kazhdan-Lusztig polynomials were introduced and a
parabolic version was later defined by Deodhar \cite{Deo} (also
\cite{CC}). Here we follow the presentation in \cite{So2}.

Throughout this subsection we assume that $f$ is anti-dominant
(see Subsection \ref{sec:tensormodule}). Let $S_f$ denote the
stabilizer subgroup (of $S_{m+n}$) of $f$. Denote by $\mathcal
H_f$ the Iwahori-Hecke algebra associated to $S_f$. Denote by
${\bf 1}_{\mathcal H_f}$ the one-dimensional right $\mathcal
H_f$-module with basis element $1$ such that $1 \cdot H_i =q^{-1}
1$. Set
\begin{eqnarray*}
 \mathcal M^f = {\bf 1}_{\mathcal H_f} \otimes_{\mathcal H_f}
\mathcal H_{m+n}, \text{ and }
 M_x = 1 \otimes H_x, \quad x \in S_{m+n}.
\end{eqnarray*}
$\mathcal M^f$ carries a natural right $\mathcal H_{m+n}$-module
structure. Denote by $S^f$ the set of minimal length
representatives for the right cosets $S_f \backslash S_{m+n}.$
Then $\mathcal M^f$ has a basis $M_x, x \in S^f$. There exists a
unique bar involution on $\mathcal M^f$ which satisfies
$\overline{M_x} = 1\otimes \overline{H_x}$ and $\overline{MH}
=\overline M \cdot\overline{H}$ for all $M \in \mathcal M^f, H \in
\mathcal H_{m+n}$. According to \cite{KL, Deo, So2}, $\mathcal
M^f$ admits the Kazhdan-Lusztig bases $\{\underline{M}_x \}, \{
\tilde{\underline{M}}_x \}$, where $x \in S^f$, which are
characterized by the following properties:

(a) every $\underline{M}_x, \tilde{\underline{M}}_x$ with $x\in
S^f$ is bar-invariant;

(b) 
$
 {\underline{M}}_x = M_x +\sum_{y } m_{y,x}(q) M_y,\;
 \tilde{\underline{M}}_x = M_x +\sum_{y} \tilde{m}_{y,x}(q) M_y,
$
where $m_{y,x} \in q\Z[q]$ and $\tilde{m}_{y,x} \in q^{-1}
\Z[q^{-1}].$ (Furthermore, it is known that $m_{y,x} =0
=\tilde{m}_{y,x}$ unless $y < x$. By convention we set $m_{x,x}(q)
=\tilde{m}_{x,x}(q) =1$.)

By examining the $\mathcal H_{m+n}$-module homomorphism $\phi:
\mathcal M^f \rightarrow \mathbb T^{m+n}$, $M_x \mapsto \mathcal
V_{f\cdot x}$, one has the following well-known identification
(cf. \cite{So2, Br, FKK}).

%
\begin{proposition} \label{prop:KL=KL}
Suppose that $f \in\Z^{m+n}$ is anti-dominant. Then in $\mathbb
T^{m+n}$ we have
\begin{eqnarray*}
 \mathcal T_{f\cdot x} =  \mathcal V_{f\cdot x} +\sum_{y < x} m_{y,x}(q) \mathcal V_{f\cdot y} \quad
 \mathcal L_{f\cdot x} = \mathcal V_{f\cdot x} +\sum_{y < x} \tilde{m}_{y,x}(q)
 \mathcal V_{f\cdot y}.
 \end{eqnarray*}
 Equivalently, we have
 \begin{eqnarray*}
 {\mathfrak  t}_{g, f\cdot x}(q) &=& m_{y,x}(q), \text{ if } g=f \cdot y
 \text{ for some } y \in S^f,  \\
  {\mathfrak  l}_{g, f\cdot x}(q) &=& \tilde{m}_{y,x}(q), \text{ if } g=f \cdot y
 \text{ for some } y \in S^f,
 \end{eqnarray*}
 and ${\mathfrak  t}_{g, f\cdot x}$, ${\mathfrak  l}_{g, f\cdot x}$ are zero otherwise.
\end{proposition}

\subsection{Kazhdan-Lusztig conjecture and (dual) canonical bases}

We write $\mathfrak l_{g,f}(q), \mathfrak t_{g,f}(q)$ for
$\mathfrak l_{\mu,\la}(q), \mathfrak t_{\mu,\la}(q)$, where $f$,
$g$ correspond to $\la, \mu$, respectively, under the bijection
$X_{m+n}^{+} \rightarrow \Z_+^{m+n}$.

The following theorem is a reformulation in terms of dual
canonical and canonical bases of the Kazhdan-Lusztig conjecture
proved in \cite{BB, BK} combined with the translation principle of
Jantzen (cf. \cite{CC, So1, BGS, J2}) and the character formula of
tilting modules \cite{So3}. The proof here is inspired by a
similar argument as in the proof of Theorem 4.31, \cite{Br}. Such
a reformulation is known to experts.

\begin{theorem} \label{th:multiCan}
\begin{enumerate}
\item In the Grothendieck group $G(\mathcal O_{m+n})$, for $\nu
\in X_{m+n}$, we have
\begin{eqnarray*}
[\mathcal T_n(\nu)] = \sum_{\mu \in X_{m+n}} \mathfrak
t_{\mu,\nu}(1) [\mathcal V_n(\mu)],\quad
[\mathcal L_n(\nu)] =\sum_{\mu \in X_{m+n}} \mathfrak
l_{\mu,\nu}(1) [\mathcal V_n(\mu)].
\end{eqnarray*}

\item In the Grothendieck group $G(\OPn)$, for $\nu
\in X_{m+n}^{+}$, we have
\begin{eqnarray*}
[\mathcal U_n(\nu)] =\sum_{\mu \in X_{m+n}^{+}} \mathfrak
u_{\mu,\nu}(1) [\mathcal K_n(\mu)],\quad
[\mathcal L_n(\nu)] =\sum_{\mu \in X_{m+n}^{+}} \mathfrak
l_{\mu,\nu}(1) [\mathcal K_n(\mu)].
\end{eqnarray*}
\end{enumerate}
\end{theorem}

\begin{proof}
Let $\la \in X_{m+n}$ be such that $f_{-\la -2\rho_c} \in
\Z^{m+n}$ satisfies that $f_{-\la -2\rho_c}(-m) \le \cdots  \le
f_{-\la -2\rho_c}(n)$. The dot action of the Weyl group $S_{m+n}$
on $X_{m+n}$ is given by: $\sigma \cdot \mu = \sigma (\mu +\rho_c)
-\rho_c$, where $\sigma \in S_{m+n}, \mu \in X_{m+n}$. Let $W_\la$
be the stabilizer (in $S_{m+n}$) of $\la$ under the dot action,
$W^\la$ the set of maximal length representatives of the left
cosets $S_{m+n}/W_\la$, and $w_\la$ the longest element in
$W_\la$. By the Kazhdan-Lusztig conjecture, combined with the
translation principle, we have
\begin{eqnarray} \label{KLformula}
[\mathcal V_n (\sigma \cdot \la) : \mathcal L_n(\tau \cdot \la) ]
=P_{\sigma,\tau} (1) =m_{w_\la \sigma^{-1}, w_\la \tau^{-1}}(1),
\quad \sigma, \tau \in W^\la.
\end{eqnarray}
Here $P_{\sigma,\tau}(1)$ is the value at $1$ of the usual
Kazhdan-Lusztig polynomial \cite{KL}, and it is equal to $m_{w_\la
\sigma^{-1}, w_\la \tau^{-1}}(1)$ in the notation of the previous
subsection, according to \cite{So2}. By
Proposition~\ref{prop:KL=KL} we have $m_{w_\la \sigma^{-1}, w_\la
\tau^{-1}}(1) = \mathfrak t_{-\sigma \cdot \la -2 \rho_c, -\tau
\cdot \la -2 \rho_c} (1)$ for $\sigma, \tau \in W^\la$.
%
Thus by (\ref{KLformula}) we have
\begin{eqnarray} \label{eq:KLtype}
[\mathcal V_n (\sigma \cdot \la) : \mathcal L_n(\tau \cdot \la) ]
=\mathfrak t_{-\sigma \cdot \la -2 \rho_c, -\tau \cdot \la -2
\rho_c} (1).
\end{eqnarray}
Combining (\ref{eq:tiltT}) and (\ref{eq:KLtype}) we have proved
the first identity in part~(1).

By Remark~\ref{rem:duality}, the matrices $[\mathfrak t_{-\sigma
\cdot \la -2 \rho_c, -\tau \cdot \la -2 \rho_c} (1)]$ and
$[\mathfrak l_{\sigma \cdot \la,\tau \cdot \la}(1)]$ are inverses
to each other.
%
The second identity in (1) follows from this.

The second identity in (2) can be derived from the second identity
in (1) and the Weyl character formula (applied to
$\mathfrak{gl}(m) \oplus \mathfrak{gl}(n)$)
(cf.~\cite[Sect.~7]{So3}).
By Remark~\ref{dualUK}, the matrices $[\mathfrak{u}_{-f_\mu\cdot
w_0,-f_\nu\cdot w_0}(1)]$ and $[\mathfrak l_{f_\mu,f_\nu}(1)]$ are
inverses to each other. Now the remaining identity in (2) follows
by applying (\ref{eq:tiltU}) and (\ref{invol}).
\end{proof}
Note that a weight $\la \in X_{m+n}^+$ is $J$-typical (i.e.~$\la$
is minimal in $X_{m+n}^+$ in the Bruhat ordering) if and only if
the generalized Verma module $\mathcal K_n(\la)$ is irreducible.

\subsection{The categories $\OPPn$ and $\Oi$}

In what follows by $\text{wt}^\epsilon$ we mean the definition
given in (\ref{nonsuperweight}). Denote by $\chi_\la$ the integral
central character associated to $\la \in X_{m+n}$. It is known
that $\chi_\la =\chi_\mu$ for $\la, \mu \in X_{m+n}$ if and only
if $\la =\sigma \cdot \mu$ for some $\sigma \in S_{m+n}$, or
equivalently $\text{wt}^\epsilon(\la)=\text{wt}^\epsilon(\mu)\in
P$. We denote by $\mathcal O^+_\g$ the block in $\OPn$ associated
to $\g \in P$.

Let $V$ be the natural $\mathfrak g\mathfrak l (m+n)$-module and
$V^*$ its dual. For $a\in \Z$ we define the translation functors
$E_a, F_a: \OPn \longrightarrow \OPn$ by sending $M \in \mathcal
O^+_\g$ to
\begin{eqnarray*}
 F_a M := \text{pr}_{\g - (\epsilon_a
 -\epsilon_{a+1})} (M \otimes V),\quad
 E_a M := \text{pr}_{\g + (\epsilon_a
 -\epsilon_{a+1})} (M \otimes V^*).
\end{eqnarray*}
Let $G(\OPn)_{\mathbb Q}:=G(\OPn)\otimes_\Z{\mathbb Q}$ and let
$\mathcal E^{m+n}|_{q=1}$ be the specialization of $\mathcal
E^{m+n}$ at $q\to 1$.
\begin{theorem} \label{red-KL-finite}
Let $n \in \N$.
\begin{enumerate}
\item Sending the Chevalley generators $E_a, F_a$ to the
translation functors $E_a, F_a$ defines a $\mathcal
U_{q=1}$-module structure on $G(\OPn)_{\mathbb Q}$.

\item The linear map $i: G(\OPn)_{\mathbb Q} \rightarrow \mathcal
E^{m+n}|_{q=1}$ which sends $[\mathcal K_n(\la)]$ to $\mathcal
K_{f_\la}(1)$ for each $\la \in X^+_{m+n}$ is an isomorphism of
$\mathcal U_{q=1}$-modules.

\item The map $i$ sends $[\mathcal U_n(\la)]$ to $\mathcal
U_{f_\la}(1)$ and $[\mathcal L_n(\la)]$ to $\mathcal
L_{f_\la}(1)$, for each $\la \in X^+_{m+n}$.
\end{enumerate}
\end{theorem}

\begin{proof} The map $i$ is certainly a vector space isomorphism.
One checks that the action of the translation functors on the
generalized Verma modules is compatible with the action of the
Chevalley generators of $\mathcal U_{q=1}$ on the monomial basis.
Thus (1) and (2) follow.  Now (3) follows from Theorem
\ref{th:multiCan} and the definition of KL polynomials $\mathfrak
u_{\mu,\nu}$ and $\mathfrak l_{\mu,\nu}$.
\end{proof}

\begin{corollary}
Let $\la \in X_{m+n}^+$. The subquotients of a generalized Verma
flag of $\mathcal U_n(\la)$ are precisely $\mathcal K_n(\mu)$ with
$\mu =  {\mathbb L}_\theta (\la)$ associated to $\theta \in
\{0,1\}^{\divideontimes f_\la}$. Furthermore for $\mathcal X_a \in
\{E_a, F_a \}$ corresponding to the Chevalley operators in
Procedure~\ref{pro321}, $\mathcal X_a \mathcal U_n(\la)$ is a
tilting module. Also we have $\mathcal U_n(\la)^\tau \cong
\mathcal U_n(\la)$.
\end{corollary}

The corollary above can be proved using induction based on
Procedure~\ref{pro321} in an analogous way as Theorem~4.37 in
\cite{Br} is proved. The induction procedure also shows that the
indexing set $\{ \mu = {\mathbb L}_\theta (\la) \mid \theta \in
\{0,1\}^{\divideontimes f_\la}\}$ above is compatible with the
indexing set in Theorem~\ref{prop:aux84} for $f =f_\la$. Finally
$\mathcal U_n(\la)^\tau \cong \mathcal U_n(\la)$ is a consequence
of Procedure \ref{pro321} and the fact that $\tau$ commutes with
the translation functors.

Theorem \ref{red-KL-finite} and results from Section~\ref{sec:bar}
imply the stability of the composition factors in a generalized
Verma module, and the stability of generalized Verma flags in a
tilting module. This allows us to apply the machinery of
Section~\ref{sec:superRep} to the reductive setting.  We will
state these results below.

Given $\la, \mu \in X^{+}_{m+\infty}$, we may regard $\la, \mu \in
X^{++}_{m+n}$ for $n\gg0$. Then,
$$[\mathcal K(\la): \mathcal L(\mu)] =[\mathcal K_{n}(\la):
\mathcal L_{n}(\mu)], \quad \text{for } n \gg0.$$

For $n<n' \le \infty$, the {\em truncation functor}
${\textsf{tr}}_{n',n}: \mathcal O_{m+n'}^{++} \longrightarrow
\OPPn$
is defined by sending an object $M$ to the $\mathfrak g\mathfrak
l(m+n)$-module
$${\textsf{tr}}_{n',n}(M) := \text{span } \{ v \in M \mid (\wt(v)| \delta'_{k})_c = 0
\text{ for all }n+1 \le k \le n'\}$$
When $n'$ is clear from the context we will also write
${\textsf{tr}}_{n}$ for $\textsf{tr}_{n',n}$.
Theorem \ref{red-KL-finite} allows us to apply the construction of
Section \ref{sec:superRep} to prove the existence of a unique
tilting module $\mathcal U(\la)$ associated to $\la$ in $\Fi$,
isomorphic to $\cup_n \mathcal U_n(\la)$. Moreover,
$$(\mathcal U(\la):\mathcal K(\mu))
=(\mathcal U_n(\la):\mathcal K_n(\mu)), \quad \text{ for } n\gg 0,
$$
and $\mathcal U(\la)$ has a generalized Verma flag whose
subquotients are precisely $\mathcal K(\mu)$ with $\mu =  {\mathbb
L}_\theta (\la)$ associated to $\theta \in \{0,1\}^{\divideontimes
f_\la}$. Furthermore, $U(\la)^\tau \cong U(\la)$.

Let $n < n'\le \infty$ and $\la\in X^{++}_{m+n'}$. Then
$\textsf{tr}_n$ sends $\mathcal Y_{n'} (\la))$ to $\mathcal
Y_n(\la)$ if $(\la|\delta'_{n+1})_c =0$ and to $ 0$ otherwise, for
$\mathcal Y=\mathcal L, \mathcal K$ or $\mathcal U$.

Similarly one proves the reductive analogue of Proposition
\ref{transl=chev}.

\begin{proposition} \label{rem:isom}
For a suitable topological completion $\widehat{G}(\Oi)_{\mathbb
Q}$ of $G(\Oi)_{\mathbb Q}$, the linear map $i:
\widehat{G}(\Oi)_{\mathbb Q} \rightarrow \widehat{\mathcal
E}^{m+\infty}|_{q=1}$ which sends $[\mathcal K(\la)]$ to $\mathcal
K_{f_\la}(1)$ for each $\la\in X^+_{m+\infty}$ is an isomorphism
of vector spaces. The map $i$ further identifies $[\mathcal
U(\la)]$ with $\mathcal U_{f_\la}(1)$ and $[\mathcal L(\la)]$ with
$\mathcal L_{f_\la}(1)$.
\end{proposition}
\section{Super duality}\label{sec:isom}
\subsection{An isomorphism of Fock spaces}

\begin{proposition}
\begin{enumerate}
 \item The $\mathcal U$-module $\wgv$ is isomorphic to the basic
 representation of $\mathcal U$ with highest weight vector
 $\vac = v_0 \wedge v_{-1}\wedge v_{-2} \wedge  \cdots.$ More
 explicitly, $E_a \vac =0$ and $K_{a, a+1} \vac =q^{\delta_{a,0}}
 \vac, $ for all $a \in\Z;$
 \item The $\mathcal U$-module $\wgw$ is isomorphic  to the basic
 representation of $\mathcal U$ with highest weight vector
 $\vacc = w_1\wedge w_2 \wedge w_3 \wedge  \cdots.$
\end{enumerate}
\end{proposition}

\begin{proof}
We prove  (2). Using the formulas of the $\mathcal U$-action on
$\mathbb W$, we verify that:

a) $E_a \vacc =0$ for every $a \in \Z$;

b) $K_a \vacc =\vacc, a \le 0$ and $K_b \vacc =q^{-1}\vacc, b >0$.

This implies that the weights of $\vacc$ by $K_{a,a+1} =K_a
K_{a+1}^{-1} (= q^{h_a})$ are $1$ for $a \ne 0$ and $q$ for $a=0$.
This says that the weight is the fundamental weight $\Lambda_0$.
The claim now follows from the fact that $\wgw$ and the basic
representation of $\mathcal U$ have the same character.

Part (1) can be proved similarly and can be found in \cite{MM,
KMS}.
\end{proof}

Thus we have obtained a natural isomorphism of $\mathcal
U$-modules $$C: \wgv \stackrel{\cong}{\longrightarrow} \wgw.$$

Consider the lattice on the fourth quadrant of the $xy$-plane. Let
us label unit intervals on the $x$-axis $1, 2, 3, \cdots$ and
those on the $y$-axis $0, -1, -2, \cdots$. (See the example
below.) A lattice path has one end going down vertically along the
$y$-axis and the other end going to the right horizontally along
the $x$-axis. Clearly, there is a bijective correspondence between
such lattice paths and Young diagrams.

The edges of a lattice path are further labeled by integers; and
the labels are uniquely determined by two requirements: a) the set
of all labels coincides with $\Z$; b) the labels on the edges
which lie on the $x$- and $y$- axis coincide with the pre-fixed
labels therein. By abuse of notation, we will use $\la$ to refer
to both a partition and its associated lattice path. Note that our
labeling differs from the one in \cite{MM}, pp.~81, by a shift.
Denote by $\text{vL}(\la)$ the set of vertical labels of $\la$ and
by $\text{hL}(\la)$ the set of horizontal labels of $\la$. For
example, associated to the partition $\la =(5,3,2,2)$, the lattice
path is the bold line segments with labels attached.

\begin{equation}\label{diagram}
{\beginpicture \setcoordinatesystem units <1.5pc,1.5pc> point at 0
2 \setplotarea x from 0 to 1.5, y from -3 to 4 \plot 0 0 0 4 1 4 1
0 0 0 / \plot 1 1 2 1 2 4 1 4 / \plot 2 2 3 2 3 4 2 4 / \plot 3 3
5 3 5 4 3 4 / \plot 0 1 1 1 / \plot 0 2 2 2 / \plot 0 3 3 3 /
\plot 4 3 4 4 / \plot 1 0 2 0 2 1 / \plot 0 0 0 -2 /  \plot 5 4 8
4 / \plot 0.02 0 0.02 -2 / \plot 0 -0.02 2 -0.02 / \plot 1.98 0
1.98 2 / \plot 1.98 2.02 3 2.02 / \plot 2.98 2 2.98 3 / \plot 2.98
2.98 5 2.98 / \plot 4.98 3 4.98 4 / \plot 5 3.98 8 3.98 / \plot 0
0.01 2 0.01 / \plot 1.97 0 1.97 2 / \plot 1.97 2.03 3 2.03 / \plot
2.97 2 2.97 3 / \plot 2.97 2.97 5 2.97 / \plot 4.97 3 4.97 4 /
\plot 5 3.97 8 3.97 / \put{$0$} at -.5 3.5 \put{$-1$} at -.76 2.5
\put{$-2$} at -.76 1.5 \put{$-3$} at -.76 0.5 \put{$\vdots$} at
-.5 -1.5 \put{$1$} at 0.5 4.5 \put{$2$} at 1.5 4.5 \put{$3$} at
2.5 4.5 \put{$4$} at 3.5 4.5 \put{$5$} at 4.5 4.5 \put{$\cdots$}
at 7.7 4.5 \put{\Small {-3}} at .5 .2 \put{\Small {-2}} at 1.5 .2
\put{\Small {-1}} at 2.2 0.5 \put{\Small {0}} at 2.15 1.5
\put{\Small {1}} at 2.5 2.2 \put{\Small {2}} at 3.15 2.5
\put{\Small {3}} at 3.5 3.2 \put{\Small {4}} at 4.5 3.2
\put{\Small {5}} at 5.15 3.5 \put{$\shortmid$} at 6 4
\put{$\shortmid$} at 7 4 \put{-} at 0.03 -1 \put{$-4$} at -.76
-0.5 \put{$6$} at 5.5 4.5 \put{$7$} at 6.5 4.5
\endpicture}
\end{equation}

The following lemma seems to be well known.

\begin{lemma} \label{lem:label}
 Let $\la=(\la_1, \la_2, \cdots)$ be a partition, and let $\la'= (\la'_1, \la'_2, \cdots)$
 be its conjugate partition. Then,
 \begin{enumerate}
 \item the set of horizontal labels $\text{hL}(\la)$ is $\{i-\la'_i, i=1,2, \cdots\};$
 \item the set of vertical labels $\text{vL}(\la)$ is $\{\la_i -i+1, i=1,2, \cdots\}.$
 \end{enumerate}
\end{lemma}

The following observation plays an important role in this paper.
\begin{theorem}  \label{th:iso}
The $\mathcal U$-module isomorphism  $C: \wgv \rightarrow \wgw$ is
explicitly given by sending $|\la \rangle$ to $|\la'_* \rangle$
for each partition $\la$.
\end{theorem}

\begin{proof}
Define the notions of concave and convex corners in a lattice path
as follows:
\begin{equation}\label{diagram1}
{\beginpicture \setcoordinatesystem units <1.5pc,1.5pc> point at 0
2 \setplotarea x from 0 to 1.5, y from -2 to 1.5 \plot -3 0 -3 1
-2 1 / \plot 3 1 3 0 2 0 / \put{\Small{$i$}} at -3.2 0.5 \put
{\Small{$i+1$}} at -2.4 1.3 \put {\Small{$i$}} at 2.5 -.3 \put{
\Small{$i+1$}} at 3.6 .5 \put{{ Concave corner}} at -3.2 -1 \put{{
Convex corner}} at 3.2 -1
\endpicture}
\end{equation}
With our convention of labeling, the formulas in \cite{MM},
pp.~81--82, can be adjusted as follows, e.g.~for $F_i$.
We have $F_i |\la \rangle = |\nu \rangle$ if $\la$ has the concave
corner labeled by $i, i+1$ (see \ref{diagram1}) ({\em instead of
the $i-1, i$ in \cite{MM}}) and $\nu$ is the same as $\la$ except
this corner becomes convex labeled again by $i,i+1$, i.e., with
one particular cell added (see \ref{diagram1}). Otherwise $F_i
|\la \rangle$ equals $0$. In other words, the vertical labels of
the lattice path of $\la$, after changing the label $i$ to $i+1$
(when allowed), become the vertical labels of the lattice path of
$\nu$. By Lemma~\ref{lem:label}, the above formula is equivalent
to the statement that $F_i |\la \rangle$ is the semi-infinite
wedge obtained by replacing the $v_{i}$ appearing in $|\la
\rangle$ by $v_{i+1}$. By the way, this is consistent with the
action of $\mathcal U$ on $\mathbb V$.

Now we work out explicitly the action of $F_i$ on the basis
element $|\la'_* \rangle$ of $\wgw$. By the action of $\mathcal U$
on $\mathbb W$, $F_i |\la'_* \rangle$ is the semi-infinite wedge
in $\wgw$ which is obtained by replacing the $w_{i+1}$ appearing
in $|\la'_* \rangle$ by $w_{i}$. In light of
Lemma~\ref{lem:label}, the horizontal labels of the lattice path
associated to $F_i |\la'_* \rangle$ are obtained from those for
$\la$ by changing the label $i+1$ to $i$ (when allowed). This is
exactly the $\nu$ above. Thus we conclude that $F_i |\la'_*
\rangle = |\nu'_* \rangle$ if $\la$ has the concave corner labeled
by $i, i+1$, otherwise it is $0$.

One checks similarly that sending $|\la \rangle$ to $|\la'_*
\rangle$ respects the actions of $E_i, K_{i,i+1}$ for each
$i\in\Z$.
\end{proof}

\subsection{The match of typicality and $J$-typicality}

Given $\la = \sum_{i \in I(m|\infty)} \la_i\delta'_i \in
X_{m+\infty}^{+}$ so that by definition $\la^{>0} :=(\la_1, \la_2,
\cdots)$ is a partition. Denoting by $(\la'_1, \la'_2, \cdots)$
the conjugate partition of $\la^{>0}$, we define a weight in
$X_{m|\infty}^{+}$
$$\la^{\natural} := \sum_{i=-m}^{-1} \la_i \delta_i
+\sum_{j=1}^\infty \la_j' \delta_j.$$ This actually defines a
bijection
$X_{m+\infty}^{+} \stackrel{\natural}{\longleftrightarrow}
X_{m|\infty}^{+}.$


\begin{theorem} \label{th:typical}
\begin{enumerate}
\item For every $\la \in X_{m+\infty}^{+}$, the degree of
$J$-atypicality of the weight $\la \in X_{m+\infty}^{+}$ is equal
to the degree of atypicality of $\la^\natural$;

\item A weight $\la \in X_{m+\infty}^{+}$ is $J$-typical if and
only if $\la^\natural \in \Xmi^{+}$ is typical.
\end{enumerate}
\end{theorem}

\begin{proof}
Part (2) is a special case of (1), and thus it suffices to prove
(1).

Let us write $\la=(\la^{<0}|\la^{>0})$ as usual.  For the proof it
will be convenient to put $\mu=\la^{>0}$. By definition, the
degree of atypicality of $\la^\natural = (\la^{<0}| \mu')$ is the
number of $i$'s in $\{-m, \cdots, -1\}$ such that
$$(\lambda^\natural +\rho| \delta_i-\delta_j)= \la_i -i
+\mu_j' -j =0, \quad \text{ for some } j>0,$$
or by Lemma~\ref{lem:label}, it is the number of $i$'s in $\{-m,
\cdots, -1\}$ such that $\la_i -i$ is a horizontal label of the
lattice path $\mu$.

By definition,
the degree of $J$-atypicality of $\lambda = (\la^{<0}| \mu)$ is
the number of $i$'s in $\{-m, \cdots, -1\}$ such that
$$(\lambda +\rho_c| \delta'_i-\delta'_\ell)_c = \la_i -i -\mu_\ell
+\ell-1 \neq 0,\; \text{  for every } \ell>0,$$
or by Lemma~\ref{lem:label}, it is the number of $i$'s in $\{-m,
\cdots, -1\}$ such that $\la_i -i$ is not a vertical label of the
lattice path $\mu$.

Now the theorem follows because the set of horizontal labels of
$\mu$ and the set of vertical labels of $\mu$ are disjoint and
their union is $\Z$.
\end{proof}

\begin{remark}  \label{rem:index}
{}From the above proof, we observe that the subset of $\{-m,
\cdots, -1\}$ that contributes to the $J$-atypicality of $\la$
coincides with the one that contributes to the atypicality of
$\la^\natural$.
\end{remark}

\subsection{The canonical isomorphism}
\label{subsec:isom}

By coupling with the two bijections $X_{m|\infty}^{+} \rightarrow
\Z_+^{m|\infty}$ and $X_{m+\infty}^{+} \rightarrow
\Z_+^{m+\infty}$, the bijection
$X_{m+\infty}^{+}\stackrel{\natural}{\longleftrightarrow}
X_{m|\infty}^{+}$ induces a bijection
$\Z_+^{m+\infty}\stackrel{\natural}{\longleftrightarrow}\Z_+^{m|\infty}$,
which will be also denoted by $\natural.$

\begin{lemma}  \label{lem:move}
  For $f,g \in \Z_+^{m+\infty}$, $f \ge g$ in the Bruhat ordering if and
  only if $f^\natural \succcurlyeq g^\natural$ in the super Bruhat
  ordering.
\end{lemma}
\begin{proof}
First we note that $h \succ e$ in super Bruhat ordering for $e, h
\in \Z_+^{m|\infty}$ if and only if there exists a sequence $h_0,
h_1,\ldots, h_t \in \Z_+^{m|\infty}$ such that $h =h_0 \succ h_1
\succ \cdots \succ h_t =e$ and for each $0\le a <t$, $h_{a+1} =
(h_{a} -r(d_{i_a} -d_{j_a}))^+$ for some $r>0$ and $i_a<0<j_a$
with $h_{a} (i_a) =h_a(j_a)$. This claim can be deduced from the
proof of Lemma~3.42 in \cite{Br} using the equivalent formulation
of the super Bruhat ordering in terms of $\wt^\epsilon$ and some
simple inequalities (cf. \S 2-b, (2.4) \cite{Br}). We shall say
that $h_a$ is reduced to $h_{a+1}$ by a {\em super move}. Denote
by $(\la_<^a| \la^a) \in X^{+}_{m|\infty}$ the weight
corresponding to $h_a$. The horizontal label set $\text{hL}
(\la^{a+1})$ is obtained from $\text{hL}(\la^a)$ by substituting
$h_a(j_a)$ with the smaller number $h_a(j_a)-r$, and this
statement in turn characterizes a super move.

On the other hand, recall from Lemma~\ref{lem:simplemove} that $f
> g$ in the Bruhat ordering if and only if there exists $f^0,\cdots,
f^s \in \Z_+^{m+\infty}$ such that $f=f^0 > f^1 > \cdots > f^s =g$
and for each $0 \le a<s$ there exists $i_a<0<j_a$ with $f^{a+1} =
(f^a \cdot \tau_{i_aj_a})^+, f^a(i_a)> f^a(j_a),$ and
$f^a(i_a)\not=f^a(k)$, $\forall k>0$. We shall say $f^a$ is
reduced to $f^{a+1}$ by a {\em simple move.} Note that $f^a(i_a)$
does not belong to the vertical label set $\text{vL} (\mu^a)$ if
we denote by $(\mu_<^a|\mu^a) \in X^{+}_{m+\infty}$ the weight
corresponding to $f^a$. The set $\text{vL} (\mu^{a+1})$ is
obtained from $\text{vL}(\mu^a)$ by substituting $f^a(j_a)$ with
the larger number $f^a(i_a)$, and this statement in turn
characterizes a simple move. By Remark~\ref{rem:index},
Lemma~\ref{lem:label}, and the fact that $\Z = \text{hL} (\nu)
\bigsqcup \text{vL} (\nu)$ for any partition $\nu$, we conclude
that if $f$ is reduced to $g$ by a simple move then $f^\natural$
is reduced to $g^\natural$ by a super move, and vice versa. An
induction on the number of simple/super moves completes the proof.
\end{proof}

The $\mathcal U$-module isomorphism $C: \wgv \rightarrow \wgw$
gives rise an isomorphism of $\mathcal U$-modules ${\mathcal
E}^{m+\infty} \stackrel{\cong}{\longrightarrow} {\mathcal
E}^{m|\infty}.$ Recall that $\widehat{\mathcal E}^{m+\infty}$ and
$\widehat{\mathcal E}^{m|\infty}$ are topological completions of
$\Lambda^m \mathbb V \otimes \La^\infty \mathbb V$ and $\Lambda^m
\mathbb V \otimes \La^\infty \mathbb V^*$, respectively. Using the
fact that $\Z = \text{hL} (\la) \bigsqcup \text{vL} (\la)$ for a
partition $\la$, one can show that these two completions are
indeed compatible under the above map.  Hence we obtain a natural
isomorphism of $\mathcal U$-modules (which will also be
conveniently denoted by $\natural$ by abuse of notation):
$$\natural = 1\otimes C: \widehat{\mathcal E}^{m+\infty}
\stackrel{\cong}{\longrightarrow} \widehat{\mathcal
E}^{m|\infty}.$$

\begin{theorem} \label{correspondence}
The isomorphism $\natural: \widehat{\mathcal E}^{m+\infty}
\longrightarrow \widehat{\mathcal E}^{m|\infty}$ has the following
properties:
\begin{enumerate}

\item $\natural (\mathcal K_f) = K_{f^\natural}$ for each
$f\in\Z_+^{m+\infty}$;

 \item $\natural$ is compatible
with the bar involutions, i.e., $\natural (\bar{u})
=\overline{\natural (u)}$ for each $u \in \widehat{\mathcal
E}^{m+\infty}$;

\item $\natural (\mathcal L_f) = L_{f^\natural}$ for each
$f\in\Z_{+}^{m+\infty}$;

\item $\natural (\mathcal U_f) =U_{f^\natural}$ for each
$f\in\Z_{+}^{m+\infty}$.

\end{enumerate}
\end{theorem}

\begin{proof}
(1) follows from the definitions and Theorem~\ref{th:iso}. (2)
follows from Theorems \ref{th:iso} and \ref{th:typical}, Lemma
\ref{lem:move} and the characterizations of the bar involutions.
(3) and (4) follow from (1), (2), (Lemma \ref{lem:move}) and the
characterizations of these bases.
\end{proof}

Combining the isomorphisms $j$, $i$, $\natural$ of
Proposition~\ref{transl=chev}, Proposition~\ref{rem:isom},
Theorem~\ref{correspondence}, respectively, we conclude the
following. Denote $\sharp := j^{-1} \circ \natural \circ i.$

\begin{theorem} \label{sameGr}
\begin{enumerate}
\item There exists a linear isomorphism of Grothendieck groups
$\sharp: \widehat{G} (\Oi)_{\mathbb Q} \rightarrow
\widehat{G}(\Fi)_{\mathbb Q}$, which sends $[\mathcal K(\la)],
[\mathcal U(\la)]$ and $[\mathcal L(\la)]$ to $[K(\la^\natural)],$
$[U(\la^\natural)]$ and $[L(\la^\natural)]$ respectively, where
$\la \in X_{m+\infty}^+.$

\item
For $\la, \mu \in X_{m+\infty}^+$, we have
$$
 \mathfrak u_{\mu,\la} (q) = u_{\mu^\natural,\la^\natural} (q),
 \quad
 \mathfrak l_{\mu,\la} (q) = \ell_{\mu^\natural,\la^\natural} (q).
$$
\end{enumerate}
\end{theorem}

\begin{remark}
{}From Vogan's interpretation of the Kazhdan-Lusztig conjecture in
terms of the $\mathfrak u$-cohomology groups (\cite{Vo},
Conjecture 3.4), Brundan-Serganova's description of the
polynomials $\ell_{\mu,\la}(q)$ \cite{Se1, Br, Zou}, together with
Kostant's $\mathfrak u$-cohomology formula for finite-dimensional
irreducible modules of $\mathfrak g\mathfrak l(N)$, the
computation of the cohomology groups ${H}^i(\mathfrak g\mathfrak
l(m|n)_{+1};L_n(\la^\natural))$ in \cite{CZ} can be shown to be
equivalent to the identity $\mathfrak l_{\mu,\la} (q) =
\ell_{\mu^\natural,\la^\natural}(q)$ when $\la$ and $\mu$ are
partitions.  In the approach of \cite{CZ} the unitarity of the
module $L_n(\la^\natural)$ when $\la$ is a partition is used in a
crucial way. As noted in {\em loc.~cit.}, the Weyl character
formula for $\mathcal L(\la)$ translates into a Weyl-type
character formula for the $\glmn$-module $L_n(\la^\natural)$.
\end{remark}

Theorem \ref{sameGr} strongly suggests the following.

\begin{conjecture}
The categories $\Fi$ and $\Oi$ are equivalent.
\end{conjecture}

\begin{remark} \label{rem:KLn=N}
We summarize the precise relations between the Kazhdan-Lusztig
polynomials for $\glmn$ and the usual Kazhdan-Lusztig polynomials
for $\mathfrak{gl}(m+N)$ for {\em finite} $n$ and $N$. To compute
$u_{\mu,\la}(q)$ and $\ell_{\mu,\la}(q)$ for fixed $\la, \mu \in
\Xmn^{+}$, we may assume by Remark~\ref{rem:shift} that $\la, \mu
\in \Xmn^{++}$. Denote by $\la_\infty \in X_{m|\infty}^+$ the
extension of $\la$ by zeros, and we have $\la_\infty^\natural \in
X_{m+\infty}^+$. Write $\la_\infty^\natural =
((\la_\infty^\natural)^{<0} | (\la_\infty^\natural)^{>0})$.
Assuming the lengths of the partitions
$(\mu_\infty^\natural)^{>0}$ and $(\la_\infty^\natural)^{>0}$ are
no larger than $N$, we have $\la_\infty^{\natural,(N)},
\mu_\infty^{\natural,(N)} \in X_{m+N}^{++}.$ Then,
\begin{eqnarray*}
 u_{\mu,\la}(q)
 &=& u_{\mu_\infty,\la_\infty}(q) \qquad \text{by Theorem~\ref{th:u-k}}, \\
 &=& \mathfrak u_{\mu_\infty^\natural,\la_\infty^\natural}(q)
 \qquad \text{by Theorem~\ref{sameGr}}, \\
 &=& \mathfrak u_{\mu_\infty^{\natural, (N)},\la_\infty^{\natural,
 (N)}}(q)  \qquad \text{by Remark~\ref{aux1000}}.
\end{eqnarray*}
Similarly, we have
$$
 \ell_{\mu,\la}(q)
 = \mathfrak l_{\mu_\infty^{\natural, (N)},\la_\infty^{\natural,(N)}}(q).
 $$
Observe that, depending on $\la$ and $\mu$, the integer $N$ might
be arbitrarily large (respectively small) even if $n$ is small
(respectively large). To capture all the Kazhdan-Lusztig
polynomials for $\glmn$ with a fixed $n>0$, we use Kazhdan-Lusztig
polynomials for $\mathfrak{gl}(m+N)$ with arbitrarily large $N$.

One can also reverse the roles of Kazhdan-Lusztig polynomials for
Lie algebras and Lie superalgebras above.
\end{remark}
\subsection{Examples}

\begin{example}\label{ex:U}
Consider the element in $\Z^{4|\infty}_+$ given by
\begin{equation*}
f=(0,-1,-3,-4|-2,-1,0,4,5,6,7,\cdots).
\end{equation*}
Then $\#f=2$ and we have
\begin{align*}
&\big{(}{\texttt L}_{-4,3}(f)\big{)}^+=(-1,-3,-4,-6|-6,-2,-1,4,5,6,\cdots),\\
&\big{(}\texttt L_{-3,2}(f)\big{)}^+=(0,-3,-4,-5|-5,-2,0,4,5,6,\cdots),\\
&\big{(}\texttt L_{-3,2} \circ\texttt
L_{-4,3}(f)\big{)}^+=(-3,-4,-5,-6|-6,-5,-2,4,5,6,\cdots).
\end{align*}
Now consider
$$f^{(3)}=(0,-1,-3,-4|-2,-1,0)\in\Z^{4|3}_+
$$
with $\#f^{(3)}=2$. This is Example 3.22 in \cite{Br} in our
notational convention. Following Procedure~3.20 \cite{Br} one
obtains the following formula for $U_{f^{(3)}}$:
\begin{align*}
U_{(0,-1,-3,-4|-2,-1,0)}&=F_{-5}F_{-4}E_{-3}F_{-2}F_{-6}F_{-5}E_{-4}F_{-3}E_{-2}F_{-1}
K_{(-1,-3,-5,-6|-4,-2,0)}\\
&=K_{(0,-1,-3,-4|-2,-1,0)}+qK_{(-1,-3,-4,-6|-6,-2,-1)}\\&\
+qK_{(0,-3,-4,-5|-5,-2,0)}+q^2K_{(-3,-4,-5,-6|-6,-5,-2)}.
\end{align*}
One sees that $U_f$ is just $U_{f^{(3)}}$ with tails added
everywhere. Under $\natural$ we have
$$f^\natural=(0,-1,-3,-4|3,2,1,-3,-4,-5,-6,\cdots).$$
The truncation map from $\Z^{4+\infty}_+\rightarrow\Z^{4+3}_+$
takes $f^\natural$ to
\begin{equation*}
f^{\natural(3)}=(0,-1,-3,-4|3,2,1).
\end{equation*}
Note that $f^{\natural(3)}$ is $J$-typical, and hence
$$\mathcal U_{(0,-1,-3,-4|3,2,1)}=\mathcal K_{(0,-1,-3,-4|3,2,1)}.$$
This certainly does not correspond to $U_{(0,-1,-3,-4|-2,-1,0)}$
under $\natural$. So canonical bases of $\widehat{\mathcal
E}^{m|n}$ and $\mathcal E^{m+n}$ do not correspond for finite $n$
in general.

Now consider the $n=\infty$ case.  Apply Procedure \ref{pro321} to
$f^\natural$ repeatedly until we get a $J$-typical element. A
straightforward calculation shows that
\begin{align*}
\mathcal
U_{f^\natural}&=F_{-2}E_{-3}E_{-4}F_{-5}F_{-1}E_{-2}E_{-3}F_{-4}E_{-5}F_{-6}
\mathcal K_{(-1,-2,-4,-6|3,2,1-1,-2,-4,-6,-7,\cdots)}\\
&=\mathcal K_{(0,-1,-3,-4|3,2,1,-3,-4,-5,-6,\cdots)}+q\mathcal
K_{(-1,-3,-4,-6|3,2,1,0,-3,-4,-5,-7,\cdots)}\\&\ +q\mathcal
K_{(0,-3,-4,-5|3,2,1,-1,-3,-4,-6,-7,\cdots)}+q^2\mathcal
K_{(-3,-4,-5,-6|3,2,1,0,-1,-3,-4,-7,-8\cdots)}.
\end{align*}
We have
\begin{eqnarray*}
\big{(}{\texttt L}_{-4,3}(f)\big{)}^{+\natural}
&=&(-1,-3,-4,-6|3,2,1,0,-3,-4,-5,-7,\cdots)=\big{(}\mathbb
L_{-4}(f^\natural)\big{)}^+,\cr \big{(}\texttt
L_{-3,2}(f)\big{)}^{+\natural}
&=&(0,-3,-4,-5|3,2,1,-1,-3,-4,-6,-7,\cdots)=\big{(}\mathbb
L_{-3}(f^\natural)\big{)}^+,\cr \big{(}\texttt
L_{-3,2}\circ\texttt L_{-4,3}(f)\big{)}^{+\natural}
&=&(-3,-4,-5,-6|3,2,1,0,-1,-3,-4,-7,\cdots) \cr &=&\big{(}\mathbb
L_{-3}\circ\mathbb L_{-4}(f^\natural)\big{)}^+.
\end{eqnarray*}
Therefore $U_f$ corresponds to $\mathcal U_{f^\natural}$ under the
map $\natural$.
\end{example}

\begin{remark}
{}From Example \ref{ex:U} we see that in the Procedure~3.20 in
\cite{Br} and Procedure \ref{pro321} the functions do not
correspond at each step under the map $\natural$, because the
Chevalley generators in general are different at each step. Indeed
even in the final reduction the typical function and the
$J$-typical function do not correspond under $\natural$. So the
two procedures are different.
\end{remark}

\begin{example}
Consider the elements $f=(0,-2|-2,0,3,4,5,\cdots)$ and
$g=(-2,-4|-4,-2,3,4,5,\cdots)$ in $\Z^{2|\infty}$. Let us compute
${\ell}_{g,f}$.  The case of ${\ell}_{g^{(2)},f^{(2)}}$ is a
special case of Example 3.40 in \cite{Br}.  One sees that
only for $\theta=(2,2)$ and $\theta=(0,2)$ is it possible to have
$\texttt{R}'_\theta(g)=f$, i.e.~ we have
\begin{equation*}
\big{(}{\texttt R}_{-1,1}^2(g)\big{)}^+=f,\quad \big{(}{\texttt
R}_{-1,1}^2\circ {\texttt R}_{-2,2}^2(g)\big{)}^+=f,
\end{equation*}
so that $\ell_{g,f}(-q^{-1})=q^2+q^4$.

We have
\begin{align*}
f^\natural=(0,-2|2,1,-1,-3,-4,-5,-6\cdots),\
g^\natural=(-2,-4|2,1,0,-1,-3,-5,-6,\cdots).
\end{align*}
It is straightforward to verify that we have
\begin{equation*}
\big{(}{\mathbb R}_{-1}^2(g^\natural)\big{)}^+=f^\natural,\quad
\big{(}{\mathbb R}_{-1}^2\circ {\mathbb
R}_{-2}^2(g^\natural)\big{)}^+=f^\natural,
\end{equation*}
and furthermore these are the only $\theta$'s for which ${\mathbb
R}'_\theta(g^\natural)=f^\natural$.  Thus we have ${\mathfrak
l}_{g^\natural, f^\natural}(-q^{-1})=q^2+q^4$.  So $\ell_{g,f}=
{\mathfrak l}_{g^\natural, f^\natural}$.

On the other hand $f^{\natural(2)}$ is $J$-typical, thus $\mathcal
L_{f^{\natural(2)}}=\mathcal K_{f^\natural(2)}$. So we do not have
a correspondence between the dual canonical basis elements
$L_{f^{(2)}}$ and $\mathcal L_{f^{\natural(2)}}$ under $\natural$.
\end{example}

\subsection{Application to combinatorial characters and tensor products}

Denote by $\omega_+$ the involution of the ring of symmetric
functions in the variables $x_1,x_2,\cdots$, which in terms of
generating series in the indeterminate $t$ is given by
\begin{equation*}
\omega_+ \left({\prod_{i>0} (1-tx_i)}^{-1} \right)=\prod_{i>0}(1+t
x_i).
\end{equation*}
For a partition $\la$ we let $s_\la(x_1,x_2,\cdots)$ denote the
Schur function in the variables ${\bf x_+} := \{x_1,x_2,\cdots\}$.
The $s_\nu(x_{-m},\cdots,x_{-1})$ for finitely many variables
${\bf x_-} := \{x_{-m},\cdots,x_{-1}\}$ makes sense as a Laurent
polynomial for a generalized partition $\nu =(\nu_{-m},\cdots,
\nu_{-1})$ which by definition satisfies $\nu_{-m} \ge \cdots \ge
\nu_{-1}$ and $\nu_i \in\Z$.

Let $\mu =(\mu^{<0}|\mu^{>0}) \in X^+_{m|\infty}$. Setting
$e^{\delta_i} =x_i$, the character of $K(\mu)$ as a power series
in $\Z[[{\bf x_-}^{\pm 1}, {\bf x_+}]]$ is
\begin{equation*}
\text{ch}_{K(\mu)}={s_{\mu^{<0}}({\bf x_{-}})s_{\mu^{>0}}({\bf
x_{+})}} {\prod_{i<0<j}(1+x^{-1}_ix_j)}.
\end{equation*}
Hence by Proposition~\ref{transl=chev} the character of $L(\la)$,
where $\la\in X^+_{m|\infty}$, is
\begin{equation}\label{superchar}
\text{ch}_{L(\la)}= \sum_\mu \ell_{\mu,\la}(1){s_{\mu^{<0}}({\bf
x_{-}})s_{\mu^{>0}}({\bf x_{+})}} {\prod_{-m \le
i<0<j}(1+x^{-1}_ix_j)}.
\end{equation}

Similarly, setting $e^{\delta'_i} =x_i$ and using
Proposition~\ref{rem:isom}, the characters of $\mathcal K(\mu)$
and $\mathcal L(\la)$, where $\la, \mu \in X^+_{m+\infty}$, are
\begin{eqnarray}
\text{ch}_{\mathcal K(\mu)} &=&  {s_{\mu^{<0}}({\bf
x_{-}})s_{\mu^{>0}}({\bf x_{+})}} {\prod_{-m \le i<0<j}(1-x^{-1}_i
x_j)}^{-1},
\nonumber \\
 \text{ch}_{\mathcal L(\la)} &=& \sum_{\mu} \mathfrak l_{\mu,\la}(1)
 {s_{\mu^{<0}}({\bf x_{-}})s_{\mu^{>0}}({\bf x_{+})}}
{\prod_{-m \le i<0<j}(1-x^{-1}_ix_j)}^{-1}. \label{reductivechar}
\end{eqnarray}
Comparing (\ref{superchar}) and (\ref{reductivechar}), we have by
Theorem~\ref{sameGr} (2) for $\la\in X^+_{m|\infty}$ and $\mu\in
X^+_{m+\infty}$
\begin{equation}\label{dualchar}
\text{ch}_{L(\la)} =\omega_+(\text{ch}_{\mathcal
L(\la^\natural)}),\quad
\text{ch}_{\mathcal L(\mu)}
=\omega_+(\text{ch}_{L(\mu^\natural)}).
\end{equation}

Let $n$ be finite. Let $\la\in X^{++}_{m|n}$ and regard $\la\in
X^+_{m|\infty}$. By Corollary~\ref{cor:easyTrun} we have
\begin{equation*}
\text{ch}_{L_n(\la)}(x_{-m},\cdots,x_{-1},x_1,\cdots,x_n)=
\text{ch}_{L(\la)}(x_{-m},\cdots,x_{-1},x_1,\cdots,x_n,0,0,\cdots)
\end{equation*}
by setting the variables $x_{n+1},x_{n+2},\cdots$ in
$\text{ch}_{L(\la)}$ to $0$. (Similar remark applies to
$\text{ch}_{\mathcal L_n(\mu)}$ for $\mu\in X^{++}_{m+n}$.) This
together with (\ref{dualchar}) implies the following.

\begin{corollary}\label{involchar}
For $\la\in X^{++}_{m|n}$ and $\mu\in X^{++}_{m+n}$, we have
\begin{itemize}
\item[(1)] $\text{ch}_{L_n(\la)}(x_{-m},\cdots,x_n)=
\omega_+(\text{ch}_{\mathcal
L(\la^\natural)})(x_{-m},\cdots,x_n,0,0,\cdots)$;

\item[(2)] $\text{ch}_{\mathcal
L_n(\mu)}(x_{-m},\cdots,x_n)=\omega_+(\text{ch}_{
L(\mu^\natural)})(x_{-m},\cdots,x_n,0,0,\cdots)$.
\end{itemize}
\end{corollary}

\begin{remark}
When $\la^\natural$ is a partition, $\text{ch}_{\mathcal
L(\la^\natural)}$ is the Schur function $s_{\la^\natural}$, and
thus, $\omega_+(s_{\la^\natural})$ is the hook Schur function
associated to the partition $\la^\natural$ (cf.~\cite{BR}).
Corollary \ref{involchar}~(1) recovers the character formula for
the irreducible representation of $\mathfrak g\mathfrak l(m|n)$
appearing in the tensor powers of the natural module $\C^{m|n}$
\cite{Sv, BR}.
\end{remark}


\subsection{Isomorphism of Grothendieck rings}
We shall define product structures in the completed Grothendieck
groups $\widehat{G}(\Oi)_{\mathbb Q}$ and $\widehat{G}(\Fi)_{\mathbb
Q}$ induced by the tensor products of modules in the respective
categories. For $\la,\mu\in X^+_{m|\infty}$ consider the tensor
product of two Kac modules $K(\la)$ and $K(\mu)$.  In order to
compute their product we need to determine the Kac modules appearing
in the Kac flag of $K(\la)\otimes K(\mu)$. We have
\begin{align*}
K(\la)\otimes K(\mu)=&\Big{(}U(\glsuper \otimes_{U(\mathfrak
p)}L^0(\la)\Big{)}\otimes \Big{(}U(\glsuper)\otimes_{U(\mathfrak
p)}L^0(\mu)\Big{)}\\
\cong& U(\glsuper)\otimes_{U(\mathfrak p)}\Big{(}L^0(\la)\otimes
U(\mathfrak{gl}(m|\infty)_{-1})\otimes L^0(\mu)\Big{)},
\end{align*}
where we recall that $\mathfrak p=\mathfrak{gl}(m|\infty)_{\ge 0}$.
Since
$U(\mathfrak{gl}(m|\infty)_{-1})\cong\Lambda(\C^{m*}\otimes\C^\infty)$
as a $\mathfrak{gl}(m)\oplus\mathfrak{gl}(\infty)$-module,
$K(\la)\otimes K(\mu)$ has a Kac flag parameterized by the
$\mathfrak{gl}(m)\oplus \mathfrak{gl}(\infty)$-highest weights
appearing in the decomposition of the module
$\Lambda(\C^{m*}\otimes\C^\infty)\otimes L^0(\la)\otimes L^0(\mu)$.
By the skew-symmetric $(\mathfrak{gl},\mathfrak{gl})$-Howe duality
\cite{Ho} we have, as $\mathfrak{gl}(m)\oplus
\mathfrak{gl}(\infty)$-module,

$$\Lambda(\C^{m*}\otimes\C^\infty) \cong\sum_{\g}L^{\mathfrak{gl}(m)}(\g^*)\otimes
L^{\mathfrak{gl}(\infty)}(\g'),$$ where $\g$ is a partition with
$\ell(\g)\le m$, and $L^{\mathfrak{gl}(m)}(\g^*)$ denotes the
$\mathfrak{gl}(m)$-module dual to $L^{\mathfrak{gl}(m)}(\g)$. Since
the operator $-\sum_{i=-m}^{-1}{e_{ii}}$ provides an $\N$-gradation
on $\Lambda(\C^{m*}\otimes\C^\infty)$ with each graded component
consisting of only finitely many irreducible components, each
irreducible $\mathfrak{gl}(m)\oplus\mathfrak{gl}(\infty)$-component
in $\Lambda(\C^{m*}\otimes\C^\infty)\otimes L^0(\la)\otimes
L^0(\mu)$ appears with finite multiplicity. Thus $K(\la)\otimes
K(\mu)$ admits an infinite filtration of Kac modules such that each
Kac module appears with finite multiplicity.  This implies that the
product of $[K(\la) \otimes K(\mu)]$ in $\widehat{G}(\Fi)_{\mathbb
Q}$ is well-defined.

Now consider the tensor product of two generalized Verma modules
$\mathcal K(\la^\natural)$ and $\mathcal K(\mu^\natural)$.  An
analogous argument shows that $\mathcal
K(\la^\natural)\otimes\mathcal K(\mu^\natural)$ has a filtration of
generalized Verma modules parameterized by the irreducible
$\mathfrak{gl}(m)\oplus\mathfrak{gl}(\infty)$-components appearing
in the decomposition of $S(\C^{m*}\otimes\C^\infty)\otimes
L^0(\la^\natural)\otimes L^0(\mu^\natural)$.  By the symmetric
$(\mathfrak{gl},\mathfrak{gl})$-Howe duality \cite{Ho}, as
$\mathfrak{gl}(m)\oplus \mathfrak{gl}(\infty)$-module,
$$S(\C^{m*}\otimes\C^\infty)\cong\sum_{\g}L^{\mathfrak{gl}(m)}(\g^*)\otimes
L^{\mathfrak{gl}(\infty)}(\g),$$ where again the sum is over all
partitions $\g$ with $\ell(\g)\le m$. Thus the tensor product admits
an infinite generalized Verma flag with each generalized Verma
module appearing with finite multiplicity.

From the description of the
$\mathfrak{gl}(m)\oplus\mathfrak{gl}(\infty)$-modules
$\Lambda(\C^{m*}\otimes\C^\infty)\otimes L^0(\la)\otimes L^0(\mu)$
and $S(\C^{m*}\otimes\C^\infty)\otimes L^0(\la^\natural)\otimes
L^0(\mu^\natural)$ above, it is clear that $[K(\la)\otimes K(\mu)]$
is mapped to $[\mathcal K(\la^\natural)\otimes\mathcal
K(\mu^\natural)]$ under the map $\sharp$ given in
Theorem~\ref{sameGr}. Thus, we have established the following.

\begin{theorem} \label{ringisom}
The natural isomorphism $\sharp: \widehat{G}(\mathcal
O_{m+\infty}^{++})_{\mathbb Q} \longrightarrow \widehat{G} (\mathcal
O_{m|\infty}^{++})_{\mathbb Q}$ is an isomorphism of Grothendieck
rings. In particular, if $\la,\mu\in X^{+}_{m|\infty}$, then the
composition factors of $L(\la)\otimes L(\mu)$ are in one-to-one
correspondence with the composition factors of $\mathcal
L(\la^\natural)\otimes \mathcal L(\mu^\natural)$ via $\sharp$.
\end{theorem}


\begin{thebibliography}{ABCD}

\bibitem[BB]{BB} A. Beilinson and J. Bernstein, {\em Localisation de
$\mathfrak g$-modules}, C.R. Acad. Sci. Paris Ser. I Math. {\bf
292} (1981), 15--18.

\bibitem[BGS]{BGS} A. Beilinson, V. Ginzburg and W. Soergel,
{\em Koszul duality patterns in representation theory}, J. Amer.
Math. Soc. {\bf 9} (1996), 473--527.

\bibitem[BK]{BK} J.L. Brylinski and M. Kashiwara,
{\em Kazhdan-Lusztig conjecture and holonomic systerms}, Invent.
Math. {\bf 64} (1981), 387--410.

\bibitem[BL]{BL} J. Bernstein and D. Leites,
{\em Character formulae for irreducible representations of Lie
superalgebras of series $\mathfrak g\mathfrak l$ and $\mathfrak
s\mathfrak l$}, C.R. Acad. Bulg. Sci. {\bf 33} (1980), 1049--1051.

\bibitem[BR]{BR} A. Berele  and A. Regev, {\em Hook Young Diagrams with
Applications to Combinatorics and to Representations of Lie
Superalgebras}, Adv. Math. {\bf 64} (1987), 118--175.

\bibitem[Br1]{Br} J. Brundan,
{\em Kazhdan-Lusztig polynomials and character formulae for the
Lie superalgebra ${\mathfrak g\mathfrak l}(m|n)$}, J. Amer. Math.
Soc.  {\bf 16} (2003), 185--231.

\bibitem[Br2]{Br2} ------,
{\em Tilting modules for Lie superalgebras}, Commun.~Algebra {\bf
32} (2004), 2251-2268.



\bibitem[CC]{CC} L. Casian and D. Collingwood,
{\em The Kazhdan-Lusztig conjecture for generalized Verma
modules}, Math. Z. {\bf 195} (1987), 581--600.

\bibitem[CoI]{CoI}  D. Collingwood and R. Irving,
{\em A decomposition theorem for certain self-dual modules in the
category $\mathcal O$}, Duke Math. J. {\bf 58} (1989), 89--102.

\bibitem[CPS]{CPS} E. Cline, B. Parshall and L. Scott,
{\em Abstract Kazhdan-Lusztig theories}, Tohoku Math. J. {\bf 45}
(1993), 511--534.

\bibitem[CW]{CW} S.-J. Cheng and W. Wang, {\em Howe Duality for Lie
Superalgebras}, Compositio Math. {\bf 128} (2001), 55--94.

\bibitem[CZ]{CZ} S.-J. Cheng and R.B. Zhang,
{\em Analogue of Kostant's $\mathfrak u$-cohomology formula for
the general linear superalgebra}, Internat. Math. Res. Not.  {\bf
1} (2004), 31--53.


\bibitem[Deo]{Deo} V. Deodhar,
{\em On some geometric aspects of Bruhat orderings II: the
parabilic analogue of Kazhdan-Lusztig polynomials}, J. Algebra
{\bf 111} (1987), 483--506.

\bibitem[Don]{Don} S. Donkin,
{\em On tilting modules for algebraic groups}, Math. Z. {\bf 212}
(1993), 39--60.

\bibitem[Dr]{Dr} V.~Drinfeld, {\em Quantum groups},
Proceedings of the ICM (Berkeley, 1986), 798--820, Amer. Math. Soc.,
Providence, RI, 1987.

\bibitem[Du]{Du} J. Du,
{\em IC bases and quantum linear groups}, Proc. Symp. Pure Math.
{\bf 56} (1994), 135--148.

\bibitem[ES]{ES} T. Enright and B. Shelton,
{\em Categories of highest weight modules: applications to classical
Hermitian symmetric pairs}, Mem. Amer. Math. Soc. {\bf 67}, no. 367 (1987).

\bibitem[FKK]{FKK} I. Frenkel, M. Khovanov and A. Kirillov, Jr.,
{\em Kazhdan-Lusztig polynomials and canonical basis}, Transform.
Groups {\bf 3} (1998), 321--336.

\bibitem[Ho]{Ho} R. Howe,
{\it Perspectives on Invariant Theory: Schur Duality,
Multiplicity-free Actions and Beyond}, The Schur Lectures, Israel
Math. Conf. Proc. {\bf 8}, Tel Aviv (1992), 1--182.

\bibitem[Ir]{Ir} R.S. Irving, {\em Singular Blocks of the category $\mathcal
O$}, Math. Z. {\bf 204} (1990), 209--224.

\bibitem[Ja1]{J} J. Jantzen, {\em Kontravariante Formen auf induzierten Darstellungen
halbeinfacher Lie-Algebren}, Math. Ann. {\bf 226} (1977), 53--65.

\bibitem[Ja2]{J2} ------, {\em Moduln mit einem h\"ochsten Gewicht}, Lect. Notes in Math.
{\bf 750}, Springer Verlag, 1983.

\bibitem[Ji]{Ji} M.~Jimbo, {\em Quantum $R$ matrix for the generalized Toda
system}, Commun. Math. Phys. {\bf 102} (1986), 537--547.

\bibitem[JHKT]{JHKT} J. van der Jeugt, J. Hughes, R. King and J.
Thierry-Mieg, {\em A character formula for singly atypical modules
of the Lie superalgebra $\mathfrak sl (m|n)$}, Commun. Algebra
{\bf 18} (1990), 3453--3480.

\bibitem[Jim]{Jim} M. Jimbo,
{\em A $q$-analogue of $U({\mathfrak g\mathfrak l}(N+1))$, Hecke
algebra, and the Yang-Baxter equation}, Lett. Math. Phys. {\bf 11}
(1986), 247--252.

\bibitem[JZ]{JZ} J. van der Jeugt and R.B. Zhang, {\em Characters and
composition factor multiplicities for the Lie superalgebra
$\mathfrak{gl}(m\vert n)$}, Lett.~Math.~Phys.~{\bf 47} (1999),
49--61.

\bibitem[K1]{K1} V. Kac, {\em Lie Superalgebras}, Adv. Math. {\bf 16}
(1977), 8--96.

\bibitem[K2]{K2} ------, {\em Representations of classical
Lie Superalgebras}, Lect. Notes in Math. {\bf 676}, pp.~597--626,
Springer Verlag, 1978.

\bibitem[Kas]{Kas} M.~Kashiwara,
{\em On crystal bases of the $Q$-analogue of universal enveloping
algebras}, Duke Math. J. {\bf 63} (1991), 465--516.

\bibitem[KL]{KL} D. Kazhdan and G. Lusztig,
{\em Representations of Coxeter groups and Hecke algebras},
Invent. Math. {\bf 53} (1979), 165--184.

\bibitem[KMS]{KMS}
M.~Kashiwara, T.~Miwa, and E.~Stern, {\em Decomposition of
$q$-deformed Fock spaces}, Selecta Math. (N.S.) {\bf 1} (1995),
787--805.

\bibitem[KT]{KT} S. Khoroshkin and V. Tolstoy, {\em Universal $R$-matrix for
quantized (super)algebras}, Commun. Math. Phys. {\bf 141} (1991),
599--617.

\bibitem[KR]{KR} A.N. Kirillov and N. Reshetikhin, {\em $q$-Weyl group and a
multiplicative formula for universal $R$-matrices}, Commun. Math.
Phys. {\bf 134} (1990), 421-431.

\bibitem[Ku]{Ku} J. Kujawa,
{\em Crystal structures arising from representations of
$GL(m|n)$}, preprint, math.RT/0311251.

\bibitem[LLT]{LLT} A. Lascoux, B. Leclerc and J.-Y. Thibon,
{\em Hecke algebras at roots of unity and crystal bases of quantum
affine algebras}, Commun. Math. Phys. {\bf 181} (1996), 205--263.

\bibitem[LS]{LS} A. Lascoux and M.-P. Sch\" utzenberger,
{\em Polyn\^ omes de Kazhdan et Lusztig pour les grassmanniennes},
(French). Young tableaux and Schur functors in algebra and geometry
(Toru\'n, 1980), Ast\' erisque {\bf 87--88} (1981), 249--266.

\bibitem[Lec]{Lec} B. Leclerc,
{\em A Littlewood-Richardson rule for evaluation representations
of $U_q(\widehat{\mathfrak{sl}}_n)$}, S\'eminaire Lotharingien de
Combinatoire, {\bf B50e}, www.emis.ams.org/journals/SLC/.

\bibitem[LM]{LM} B. Leclerc and H. Miyachi, {\em Constructible characters and
canonical bases}, J.~Algebra {\bf 77} (2004), 298-317.

\bibitem[Lu1]{Lu1} G. Lusztig,
{\em Canonical bases arising from quantized enveloping algebras},
J. Amer. Math. Soc. {\bf 3} (1990), 447--498.

\bibitem[Lu2]{Lu} ------,
{\em Introduction to quantum groups}, Progress in Math. {\bf 110},
Birkh\"auser, 1993.

\bibitem[MM]{MM} K.~Misra and T.~Miwa,
{\em Crystal base for the basic representation of
$U_q(\widehat{\mathfrak{sl}}(n))$}, Commun. Math. Phys. {\bf 134}
(1990), 79--88.

\bibitem[PS]{PS} I. Penkov and V. Serganova,
{\em Representations of classical Lie superalgebras of type I},
Indag. Math. (N.S.) {\bf 3} (1992), 419--466.


\bibitem[Se1]{Se1} V. Serganova,
{\em Kazhdan-Lusztig polynomials for the Lie superalgebra
$GL(m\vert n)$}, Adv. Soviet Math. {\bf 16} (1993), 151--165.

\bibitem[Se2]{Se2} ------,
{\em Kazhdan-Lusztig polynomials and character formula for the Lie
superalgebra $\mathfrak g\mathfrak l(m\vert n)$}, Selecta Math.
(N.S.) {\bf 2} (1996), 607--651.

\bibitem[So1]{So1} W. Soergel,
{\em $\mathfrak n$-cohomology of simple highest weight modules on
walls and purity}, Invent. Math. {\bf 98} (1989), 565--580.

\bibitem[So2]{So2}  ------,
{\em Kazhdan-Lusztig polynomials and a combinatoric for tilting
modules}, Represent. Theory {\bf 1} (1997), 83--114.

\bibitem[So3]{So3}  ------,
{\em Character formulas for tilting modules over Kac-Moody
algebras}, Represent. Theory (electronic) {\bf 2} (1998),
432--448.

\bibitem[So4]{So4}  ------,
{\em Character formulas for tilting modules over quantum groups at
roots of one}. Current developments in mathematics, 1997
(Cambridge, MA), 161--172, Int. Press, Boston, MA, 1999.

\bibitem[St]{St} E. Stern,
{\em Semi-infinite wedges and vertex operators}, Internat. Math.
Res. Notices, no. {\bf 4} (1995) 201--219.

\bibitem[Sv]{Sv} A. Sergeev,
{\em The tensor algebra of the identity representation as a module
over the Lie superalgebras $\mathfrak{gl}(n|m)$ and $Q(n)$}, Math.
USSR Sbornik {\bf 51} (1985), 419--427.


\bibitem[Vo]{Vo} D. Vogan,
{\em Irreducible representations of semsimple Lie groups II: the
Kazhdan-Lusztig conjectures}, Duke Math. J. {\bf 46} (1979),
805--859.

\bibitem[Wa]{Wak} M.~Wakimoto,
{\em Infinite-dimensional Lie algebras}, Transl. Math. Monographs
{\bf 195}, Amer. Math. Soc., 2001.

\bibitem[Zou]{Zou} Y.M. Zou,
{\em Categories of finite dimensional weight modules over type I
classical Lie superalgebras}, J. Algebra {\bf 180} (1996),
459--482.

\end{thebibliography}
\end{document}